\documentclass[aop,MSNbibl,secthm,seceqn,citesort,dvips]{arximspdf}

\doi{10.1214/10-AOP538}
\volume{38}
\issue{6}
\pubyear{2010}
\firstpage{2170}
\lastpage{2223}

\makeatletter
\newcommand{\eqref}[1]{(\ref{#1})}
\DeclareMathAlphabet\mathcaligr{OMS}{cmsy}{m}{n}
\renewcommand{\mathcal}{\mathcaligr}
\renewcommand{\emptyset}{\varnothing}

\newtheorem{theorem}{Theorem}[section]
\newproclaim{definition}[theorem]{Definition}
\newtheorem{lemma}[theorem]{Lemma}
\newproclaim{remark}[theorem]{Remark}
\makeatother

\begin{document}
\begin{frontmatter}

\title{The Skorohod oblique reflection problem in time-dependent domains}
\runtitle{Skorohod problem in time-dependent domains}

\begin{aug}
\author[A]{\fnms{Kaj} \snm{Nystr\"{o}m}\thanksref{t1}\ead[label=e1]{kaj.nystrom@math.umu.se}}
\and
\author[A]{\fnms{Thomas} \snm{\"{O}nskog}\corref{}\thanksref{t2}\ead[label=e2]{thomas.onskog@math.umu.se}}
\runauthor{K. Nystr\"{o}m and T. \"{o}nskog}
\thankstext{t1}{Supported in part by Grant
VR-70629701 from the Swedish research council, VR.}
\thankstext{t2}{Supported in part by
the Swedish Defense Research Agency, FOI.}
\affiliation{Ume{\aa} University}
\address[A]{Department of Mathematics\\
\quad and Mathematical Statistics\\
Ume\aa\ University\\
SE-901 87 Ume\aa\\
Sweden\\
\printead{e1}\\
\phantom{\textsc{E-mail: }}\printead*{e2}} 
\end{aug}

\received{\smonth{5} \syear{2009}}
\revised{\smonth{12} \syear{2009}}

%
\begin{abstract}
The deterministic Skorohod problem plays an important role in the
construction and analysis of diffusion processes with reflection. In the
form studied here, the multidimensional Skorohod problem was
introduced, in
time-independent domains, by H. Tanaka \cite{Tanaka1979} and further
investigated by P.-L. Lions and A.-S. Sznitman \cite
{LionsSznitman1984} in
their celebrated article. Subsequent results of several researchers have
resulted in a large literature on the Skorohod problem in time-independent
domains. In this article we conduct a thorough study of the
multidimensional Skorohod problem in time-dependent domains. In particular,
we prove the existence of c\`{a}dl\`{a}g solutions $(x,\lambda)$ to the
Skorohod problem, with oblique reflection, for $(D,\Gamma,w)$
assuming, in
particular, that~$D$ is a time-dependent domain (Theorem \ref{Theorem 1}).
In addition, we prove that if $w$ is continuous, then $x$ is continuous as
well (Theorem \ref{Theorem 2}). Subsequently, we use the established
existence results to construct solutions to stochastic differential
equations with oblique reflection (Theorem \ref{Theorem 3}) in
time-dependent domains. In the process of proving these results we establish
a number of estimates for solutions to the Skorohod problem with bounded
jumps and, in addition, several results concerning the convergence of
sequences of solutions to Skorohod problems in the setting of time-dependent
domains.
\end{abstract}

%



\end{frontmatter}

\section{Introduction}\label{intro}

In time-independent domains the Skorohod problem, in the form studied in
this article, goes back to Tanaka \cite{Tanaka1979}, who established
existence and uniqueness of solutions to the Skorohod problem in convex
domains with normal reflection. These results were subsequently generalized
to wider classes of time-independent domains by, in particular, Lions and
Sznitman \cite{LionsSznitman1984} and Saisho \cite{Saisho1987}. By imposing
an admissibility condition on the domain, Lions and Sznitman \cite{LionsSznitman1984} proved existence and uniqueness of solutions to the
Skorohod problem in two different cases. The first of the two cases
considered normal reflection on domains satisfying a uniform exterior sphere
condition, meaning that the domain is smooth except for ``convex corners.''
Moreover, the second case considered smoothly varying (possibly oblique)
directions of reflection on smooth domains. In addition, for smoothly
varying directions of reflection on domains satisfying a uniform exterior
sphere condition, existence and uniqueness results were obtained in the
special case when the oblique reflection cone can be transformed into the
normal cone by multiplication by a smooth matrix function. Saisho \cite
{Saisho1987} later showed that in the first case considered in \cite{LionsSznitman1984}, that is, for normal reflection, the admissibility condition
is not necessary and can be removed. Moreover, concerning oblique
reflection, that is, when the cone of reflection differs from the cone
of inward
normals, we note that in the case of an orthant with constant
directions of
reflection on the sides, Harrison and Reiman \cite{HarrisonReiman1981} found
sufficient conditions for the existence and uniqueness of solutions to the
Skorohod problem as well as for continuity of the reflection map. In this
context we also mention that Bernard and El Kharroubi \cite{BernardElKharroubi1991} provided necessary and sufficient conditions for
the existence of solutions to the Skorohod problem in an orthant with
constant directions of reflection on each face. The most general
results so
far concerning the existence of solutions to the Skorohod problem with
oblique reflection in time-independent domains were derived by Costantini
\cite{Costantini1992}. Costantini \cite{Costantini1992} proved
existence of
solutions to the Skorohod problem for domains satisfying a uniform exterior
sphere condition with a nontangential reflection cone given as a continuous
transformation of the normal cone. Note that this allowed for discontinuous
directions of reflection at the corners. The question of uniqueness of
solutions to the Skorohod problem with oblique reflection is, in general,
still an open question and has been settled only in some specific
cases. For
example, Dupuis and Ishii \cite{DupuisIshii1991} obtained uniqueness
for a
convex polyhedron with constant directions of reflection on the faces
assuming the existence of a certain convex set, defined in terms of the
normal directions and the directions of reflection. Dupuis and Ishii
\cite{DupuisIshii1993,DupuisIshii2008} later extended this result to
piecewise smooth domains with smoothly varying directions of reflection on
each face. In addition, we here also mention the work of Dupuis and Ramanan
\cite{DupuisRamanan1999b,DupuisRamanan1999a} based on convex duality
techniques. In particular, in \cite{DupuisRamanan1999a} convex duality is
used to transform the condition of Dupuis and Ishii \cite{DupuisIshii1991}
into one that is much easier to verify. Before we proceed, we here note that
the outline above is an attempt to briefly discuss relevant previous
developments concerning the Skorohod problem in the form studied in this
article. In particular, the study of reflected diffusion based on Skorohod
problems was first introduced by Skorohod \cite{Skorohod1961a} and this
approach has, as briefly described, subsequently been developed in many
articles, including \cite{Costantini1992,DupuisIshii1991,HarrisonReiman1981,LionsSznitman1984,Saisho1987} and
\cite{Tanaka1979}. However, we emphasize that the literature devoted to Skorohod
problems, their extensions and applications is much larger than what is
conveyed above and, in fact, many more researcher have contributed to this
rich field. In particular, applied areas where Skorohod problems occur
include heavy traffic analysis of queueing networks (see, e.g, \cite{AtarDupuis1999,DupuisRamanan1998,DupuisRamanan2000b,Kushner2001,RamananReiman2003,RamananReiman2008,Reiman1984,Robert2003}), control theory, game theory and
mathematical economics (see, e.g., \cite{AtarDupuis2002,Kruk2000,Ramasubramanian2000,Ramasubramanian2006,SonerShreve1989}), image processing (see, e.g., \cite{Borkowski2007}) and
molecular dynamics (see, e.g., \cite{Saisho1988,Saisho1991,Saisho1994}). For further results concerning Skorohod problems, as
well as
applications of Skorohod problems, we also refer to
 \cite{AtarBudhirajaRamanan2008,BassHsu2000,BenjaminiChenRodhe2004,BurdzyToby1995,DupuisRamanan2000a,DupuisWilliams1994,ElKharroubiBenTaharYaacoubi2002,Frankowska1985,KangWilliams2007,KrukLehoczkyRamananShreve2007,Marin-RubioReal2004,Ramanan2006}
 and \cite{Taksar1992}.

An important novelty of this article is that we conduct a thorough
study of
the Skorohod problem, and the subsequent applications to stochastic
differential equations reflected at the boundary, in the setting of
time-dependent domains. To our knowledge, the Skorohod problem is
indeed less
developed in time-dependent domains. In particular, a first treatment
of the
Skorohod problem in time-dependent domains was given by Costantini, Gobet
and El Karoui \cite{CostantiniGobetKaroui2006}, who proved existence and
uniqueness of solutions to the Skorohod problem with normal reflection in
smooth time-dependent domains. Moreover, existence and uniqueness for
deterministic problems of Skorohod type in time-dependent intervals have
recently also been established by Burdzy, Chen and Sylvester
\cite{BurdzyChenSylvester2003}, Burdzy, Kang and Ramanan \cite{BurdzyKangRamanan2009}. The main
contribution of this article is that we are able to generalize the results
in \cite{Costantini1992}, concerning c\`{a}dl\`{a}g solutions to the
Skorohod problem with oblique reflection, to time-dependent domains assuming
less regularity on the domains compared to \cite{CostantiniGobetKaroui2006}.
Note also that in \cite{NystromOnskog2009c} we use the results of this
article to construct a numerical method for weak approximation of stochastic
differential equations with oblique reflection in time-dependent domains.
Finally, as in \cite{BurdzyKangRamanan2009}, we note, in particular, that
reflecting Brownian motions in time-dependent domains arise in queueing
theory (see, e.g., \cite{KonstantopoulosAnantharam1995,MandelbaumMassey1995}), statistical physics,
(see, e.g., \cite{BurdzyNualart2002,SoucaliucWerner2002}), control theory (see,
e.g., %
\cite{ElKarouiKaratzas1991b,ElKarouiKaratzas1991a}) and
finance (see,
e.g., \cite{ElKarouiKapoudjianPardouxPengQuenez1997}). In particular, in
future articles we hope to be able to further explore the results and
techniques developed in this article in several applications.

To properly formulate the multidimensional Skorohod problem considered in
this article, and our results, we in the following first have to introduce
some notation. Given $d\geq1$, we let $ \langle\cdot,\cdot
\rangle$ denote the standard inner product on $%
\mathbb{R}
^{d}$ and we let $ \vert z \vert= \langle z,z \rangle
^{1/2}$ be the Euclidean norm of $z.$ Whenever $z\in
\mathbb{R}
^{d}, r>0$, we let $B_{r} ( z ) = \{ y\in
\mathbb{R}
^{d}\dvtx  \vert z-y \vert<r \} $ and $S_{r} ( z )
= \{ y\in
\mathbb{R}
^{d}\dvtx  \vert z-y \vert=r \} $. Moreover,\vspace*{1pt} given $D\subset
\mathbb{R}
^{d+1}$, $E\subset
\mathbb{R}
^{d}$, we let $\bar{D}$, $\bar{E}$ be the closure of $D$ and $E$,
respectively, and we let $d ( y,E ) $ denote the Euclidean distance
from $y\in
\mathbb{R}
^{d}$ to $E$. Given $d\geq1$, $T>0$ and an open, connected set
$D^{\prime
}\subset
\mathbb{R}
^{d+1}$, we will refer to
%
\begin{equation}
D=D^{\prime}\cap([0,T]\times
\mathbb{R}
^{d}), \label{timedep}
\end{equation}
as a time-dependent domain. Given $D$ and $t\in [ 0,T ] $, we
define the time sections of $D$ as $D_{t}= \{ z\dvtx  ( t,z ) \in
D \} $, and we assume that%
%
\begin{equation}
D_{t}\neq\emptyset\mbox{ and that }D_{t}\mbox{ is bounded and connected
for every }t\in [ 0,T ] . \label{timedep+}
\end{equation}
We let $\partial D$ and $\partial D_{t}$, for $t\in [ 0,T ] $,
denote the boundaries of $D$ and $D_{t}$, respectively. A convex cone of
vectors in $%
\mathbb{R}
^{d}$ is a subset $\Gamma\subset
\mathbb{R}
^{d}$ such that $\alpha u+\beta v\in\Gamma$ for all $\alpha,\beta
\in
\mathbb{R}
_{+}$ and all $u,v\in\Gamma$. We let $\Gamma=\Gamma_{t}(z)=\Gamma(t,z)$
be a function defined on $%
\mathbb{R}
^{d+1}$ such that $\Gamma_{t} ( z ) $ is a closed convex cone of
vectors in $%
\mathbb{R}
^{d}$ for every $z\in\partial D_{t}$, $t\in [ 0,T ] $. To give an
example of a closed convex cone, we consider the set $C=C_{\Omega
}=\{\lambda\gamma\dvtx  \lambda>0, \gamma\in\Omega\}$, where $\Omega$ is
a closed, connected subset of $S_{1}(0)$ satisfying $\gamma_{1}\cdot
\gamma
_{2}>-1$ for all $ \gamma_{1},\gamma_{2}\in\Omega$. Given $C$, we
define $%
C^{\ast}=\{\alpha u+\beta v\dvtx  \alpha,\beta\in
\mathbb{R}
_{+}, u,v\in C\}$. Then $C^{\ast}$ is an example of a closed convex cone
and we note that $C^{\ast}=C_{\Omega^{\ast}}^{\ast}$, where $\Omega
^{\ast}$ can be viewed as the ``convex hull'' of $\Omega$ on $S_{1}(0)$.
Given $\Gamma=\Gamma_{t}(z)$, we let $\Gamma_{t}^{1}(z):=\Gamma
_{t}(z)\cap S_{1}(0)$. Given $T>0$, we let $\mathcal{D} ( [ 0,T%
] ,%
\mathbb{R}
^{d} ) $ denote the set of c\`{a}dl\`{a}g functions $w=w_{t}\dvtx  [ 0,T%
] \rightarrow
\mathbb{R}
^{d}$, that is, functions which are right continuous with left limits.
For $w\in
\mathcal{D} ( [ 0,T ] ,%
\mathbb{R}
^{d} ) $ we introduce the norm
%
\begin{equation}
\Vert w \Vert_{t_{1},t_{2}}=\sup_{t_{1}\leq r\leq s\leq
t_{2}} \vert w_{s}-w_{r} \vert\label{vvv3}
\end{equation}
for $0\leq t_{1}\leq t_{2}\leq T$ and, given $\delta>0$, we let
%
\begin{equation}
\mathcal{D}^{\delta} ( [ 0,T ] ,%
\mathbb{R}
^{d} ) = \Bigl\{ w\in\mathcal{D} ( [ 0,T ] ,%
\mathbb{R}
^{d} ) \dvtx \sup_{t} \vert w_{t}-w_{t^{-}} \vert<\delta \Bigr\}
\label{vvv1}
\end{equation}
denote the set of c\`{a}dl\`{a}g functions with jumps bounded by
$\delta$.
We denote the set of functions $\lambda=\lambda_{t}\dvtx  [ 0,T ]
\rightarrow
\mathbb{R}
^{d}$ with bounded variation by $\mathcal{BV} ( [ 0,T ] ,%
\mathbb{R}
^{d} ) $ and we let $ \vert\lambda \vert$ denote the total
variation of $\lambda\in\mathcal{BV} ( [ 0,T ] ,%
\mathbb{R}
^{d} ) $.

In this article we consider the Skorohod problem in the following form.

\begin{definition}
\label{skorohodprob} Let $d\geq1$ and $T>0$. Let $D\subset
\mathbb{R}
^{d+1}$ be a time-dependent domain satisfying \eqref{timedep+} and let
$%
\Gamma=\Gamma_{t}(z)$ be, for every $z\in\partial D_{t}$, $t\in
\lbrack
0,T]$, a closed convex cone of vectors in $%
\mathbb{R}
^{d}$. Given $w\in\mathcal{D} ( [ 0,T ] ,%
\mathbb{R}
^{d} ) $, with $w_{0}\in\overline{D_{0}}$, we say that the pair $%
( x,\lambda ) $ is a solution to the Skorohod problem for $ (
D,\Gamma,w ) $, on $ [ 0,T ] $, if $ ( x,\lambda )
\in\mathcal{D} ( [ 0,T ] ,%
\mathbb{R}
^{d} ) \times\mathcal{BV} ( [ 0,T ] ,%
\mathbb{R}
^{d} ) $ and if $ ( w,x,\lambda ) $ satisfies, for all $t\in%
[ 0,T ] $,
%
\begin{eqnarray} \label{SP1}
x_{t} &=&w_{t}+\lambda_{t},\qquad x_{t}\in\overline{D_{t}}, \\\label{SP2}
\lambda_{t}
&=&\int_{0}^{t^{+}}\gamma_{s}\,d \vert\lambda \vert
_{s},\qquad \gamma_{s}\in\Gamma_{s}^{1} ( x_{s} ) \, d \vert
\lambda \vert\mbox{-a.e on }\bigcup_{s\in [ 0,t ] }\partial
D_{s}
\end{eqnarray}
and
%
\begin{equation}
d \vert\lambda \vert\bigl ( \{ t\in [ 0,T ]
\dvtx  ( t,x_{t} ) \in D \} \bigr) =0. \label{SP4}
\end{equation}
\end{definition}

The main results of this article will be proved for time-dependent
domains $%
D\subset
\mathbb{R}
^{d+1}$ satisfying \eqref{timedep+}. However, several additional
restrictions will be imposed on $D$, on the cones of reflection $\Gamma
$ as
well as on the interaction between $D$ and $\Gamma$. In the following we
will outline these assumptions in order to be able to properly state our
existence result concerning the Skorohod problem with oblique reflection.
However, while these assumptions are introduced quite briefly here, the
intuition behind the assumptions, as well as the implications of the
assumptions, is explained in more detail in Section~\ref{geometry} below.

\subsection*{Geometry of the time-slice $D_{t}$}

We let $N_{t} (
z ) $ denote the cone of inward normal vectors at $z\in\partial D_{t}$%
, $t\in\lbrack0,T]$; see \eqref{conenorm} below for a definition. In
particular, we assume that $N_{t} ( z ) \neq\emptyset$ whenever $%
z\in\partial D_{t}$, $t\in\lbrack0,T]$. Note that we allow for the
possibility of several inward normal vectors at the same boundary point.
Given $N_{t} ( z ) $, we let $N_{t}^{1}(z):=N_{t}(z)\cap S_{1}(0)$.
Then the spatial domain $D_{t}$ is said to verify the uniform exterior
sphere condition if there exists a radius $r_{0}>0$ such that%
%
\begin{equation}
B_{r_{0}}(z-r_{0}n)\subseteq([0,T]\times
\mathbb{R}
^{d}\setminus D_t)\cap(%
\mathbb{R}
^{d+1}\setminus D), \label{extsphere-}
\end{equation}
whenever $z\in\partial{D_{t}}$, $n\in N_{t}^{1} ( z ) $. Note
that $B_{r_{0}}(z-r_{0}n)$ is the open Euclidean ball with center $z-r_{0}n$
and radius $r_{0}$. We say that a time-dependent domain $D$ satisfies a
uniform exterior sphere condition in time if the uniform exterior sphere
condition in \eqref{extsphere-} holds, with the same radius $r_{0}$,
for all
spatial domains $D_{t}$, $t\in [ 0,T ] $.

\subsection*{Temporal variation of the domain}

Following \cite{CostantiniGobetKaroui2006}, we let
%
\begin{equation}
l ( r ) =\mathop{\mathop{\sup}_{s,t\in[ 0,T]}}_{ \vert
s-t \vert\leq r}\sup_{z\in\overline{D_{s}}}d ( z,D_{t} )
\label{defl}
\end{equation}
be the modulus of continuity of the variation of $D$ in time. In particular,
in several of our estimates related to the Skorohod problem we will assume
that
%
\begin{equation}
\lim_{r\rightarrow0^{+}}l ( r ) =0. \label{limitzero}
\end{equation}

\subsection*{Cones of reflection}

Following \cite{Costantini1992}, we
assume that
%
\begin{eqnarray} \label{cone1}
&&\gamma_{1}\cdot\gamma_{2}>-1\mbox{ holds whenever } \gamma
_{1},\gamma
_{2}\in\Gamma_{t}^{1} ( z ) \mbox{ and}
\nonumber
\\[-8pt]
\\[-8pt]
&&\qquad \mbox{for all }z\in\partial
D_{t} , t\in [ 0,T ] .\nonumber
\end{eqnarray}
The assumption in \eqref{cone1} eliminates the possibility of $\Gamma$
containing vectors in opposite directions. We also assume that the set%
%
\begin{equation}
G^{\Gamma}= \{ ( t,z,\gamma ) \dvtx \gamma\in\Gamma_{t} (
z ) , z\in\partial D_{t},t\in [ 0,T ] \} \mbox{ is
closed.} \label{cone2}
\end{equation}
The interpretation of the condition in \eqref{cone2}
is discussed in Section~\ref{geometry}. In addition, we need the following
assumption concerning the variation
of the cones $\Gamma_{t}(z)$. Let
%
\begin{equation}
h(E,F)=\max(\sup\{d(z,E)\dvtx z\in F\},\sup\{d(z,F)\dvtx z\in E\}) \label{haus}
\end{equation}
denote the Hausdorff distance between the sets $E,F\subset
\mathbb{R}
^{d}$. Moreover, let $ \{ (s_{n},z_{n}) \} $ be a sequence of
points in $%
\mathbb{R}
^{d+1}$, $s_{n}\in\lbrack0,T]$, $z_{n}\in\partial D_{s_{n}}$, such
that $%
\lim_{n\rightarrow\infty}s_{n}=s\in\lbrack0,T]$, $\lim
_{n\rightarrow
\infty}z_{n}=z\in\partial D_{s}$. We assume, for any such sequence of
points $ \{ (s_{n},z_{n}) \} $, that
%
\begin{equation}
\lim_{n\rightarrow\infty}h(\Gamma_{s_{n}} ( z_{n} ) ,\Gamma
_{s} ( z ) )=0. \label{Gammanpre}
\end{equation}

\subsection*{Interaction between the geometry and the cones of
reflection}

For $z\in\partial D_{s}$, $s\in [ 0,T ] $, and $\rho
,\eta>0$ we define%
%
\begin{equation}
a_{s,z} ( \rho,\eta ) =\max_{u\in S_{1} ( 0 )
}\min_{s\leq t\leq s+\eta}\min_{y\in\partial D_{t}\cap\overline
{B_{\rho
} ( z ) }}\min_{\gamma\in\Gamma_{t}^{1} ( y )
} \langle\gamma,u \rangle \label{adef}
\end{equation}
and
%
\begin{equation} \hspace*{32pt}
c_{s,z} ( \rho,\eta ) =\max_{s\leq t\leq s+\eta}\max_{y\in
\partial D_{t}\cap\overline{B_{\rho} ( z ) }}\max_{\hat{z}\in
\overline{D_{t}}\cap\overline{B_{\rho} ( z ) }, \hat{z}\neq
y}\max_{\gamma\in\Gamma_{t}^{1} ( y ) } \biggl( \frac{ \langle
\gamma,y-\hat{z} \rangle}{ \vert y-\hat{z} \vert}\vee
0 \biggr) . \label{cdef}
\end{equation}
For technical reasons we also introduce the quantity%
%
\begin{equation}
e_{s,z} ( \rho,\eta ) =\frac{c_{s,z} ( \rho,\eta ) }{%
( a_{s,z} ( \rho,\eta ) ) ^{2}\vee a_{s,z} ( \rho
,\eta ) /2}. \label{edef}
\end{equation}
In the proof of certain a priori estimates for the Skorohod problem,
established in the bulk of the article, we will consider time-dependent
domains satisfying \eqref{timedep+} and the uniform exterior sphere
condition in time, with radius $r_{0}$. In addition, we will assume that
there exist $0<\rho_{0}<r_{0}$ and $\eta_{0}>0$, such that%
%
\begin{eqnarray}\label{crita}
\inf_{s\in [ 0,T ] }\inf_{z\in\partial D_{s}}a_{s,z} ( \rho
_{0},\eta_{0} ) &=&
a>0, \\\label{crite}
\sup_{s\in [ 0,T ] }\sup_{z\in\partial D_{s}}e_{s,z} ( \rho
_{0},\eta_{0} ) &=&
e<1.
\end{eqnarray}
Interpretations of \eqref{adef}, \eqref{cdef}, \eqref{crita} and %
\eqref{crite} are given in Section~\ref{geometry}.

\subsection*{Existence of good projections}

Let $0<\delta_{0}<r_{0}$,
$h_{0}>1$ and let $\Gamma=\Gamma_{t}(z)=\Gamma(t,z)$ be given for
all $%
z\in\partial D_{t}$, $t\in\lbrack0,T]$. We say that $([0,T]\times
\mathbb{R}
^{d})\setminus\overline{D}$ has the $(\delta_{0},h_{0})$-property of good
projections along $\Gamma$ if there exists, for any $y\in
\mathbb{R}
^{d}\setminus\overline{D}_{t}$, $t\in\lbrack0,T]$, such that
%
\begin{equation}
d ( y,D_{t} ) <\delta_{0}, \label{deltanolldef}
\end{equation}
at least one projection of $y$ onto $\partial D_{t}$ along $\Gamma_{t}$,
denoted $\pi_{\partial D_{t}}^{\Gamma_{t}} ( y ) $, which
satisfies%
%
\begin{equation}
\vert y-\pi_{\partial D_{t}}^{\Gamma_{t}} ( y ) \vert
\leq h_{0}d ( y,D_{t} ) . \label{rhonolldef}
\end{equation}

Concerning the existence and continuity of solutions to the Skorohod
problem, as defined in Definition \ref{skorohodprob}, we prove the following
two theorems.

\begin{theorem}
\label{Theorem 1} Let $T>0$ and let $D\subset
\mathbb{R}
^{d+1}$ be a time-dependent domain satisfying \eqref{timedep+}, %
\eqref{limitzero} and a uniform exterior sphere condition in time with
radius $r_{0}$ in the sense of \eqref{extsphere-}. Let $\Gamma=\Gamma
_{t} ( z ) $ be a closed convex cone of vectors in $%
\mathbb{R}
^{d}$ for every $z\in\partial D_{t}$, $t\in [ 0,T ] $, and assume
that $\Gamma$ satisfies \eqref{cone1}, \eqref{cone2} and~\eqref{Gammanpre}.
Assume that \eqref{crita} and \eqref{crite} hold for some $0<\rho
_{0}<r_{0}$%
, $\eta_{0}>0$, $a$~and~$e$. Finally, assume that $([0,T]\times
\mathbb{R}
^{d})\setminus\overline{D}$ has the $(\delta_{0},h_{0})$-property of good
projections along $\Gamma$, for some $0<\delta_{0}<\rho_{0}$, $h_{0}>1$,
as defined in \eqref{deltanolldef} and\vspace*{1pt} \eqref{rhonolldef}. Then,
given $w\in
\mathcal{D}^{( {\delta_{0}}/{4}\wedge {\rho_{0}}/({4h_{0}}))}
( %
[ 0,T ] ,%
\mathbb{R}
^{d} ) $, with $w_{0}\in\overline{D_{0}}$, there exists a solution $%
( x,\lambda ) $ to the Skorohod problem for $ ( D,\Gamma
,w ) $, in the sense of Definition \ref{skorohodprob}, with $x\in
\mathcal{D}^{\rho_{0}} ( [0,T],%
\mathbb{R}
) $.\vadjust{\goodbreak}
\end{theorem}

\begin{theorem}
\label{Theorem 2}Assume that the assumptions stated in Theorem \ref{Theorem
1} are satisfied and let $\rho_{0}$ be as in the statement of Theorem
\ref%
{Theorem 1}. Let $w\dvtx  [ 0,T ] \rightarrow
\mathbb{R}
^{d}$ be a continuous function and let $ ( x,\lambda ) $ be any
solution to the Skorohod problem for $ ( D,\Gamma,w ) $ in the
sense of Definition \ref{skorohodprob}. If $x\in\mathcal{D}^{\rho
_{0}} ( [ 0,T ] ,%
\mathbb{R}
^{d} ) $, then $x$ is continuous.
\end{theorem}

In the following remarks we have gathered comments concerning the importance
of the assumptions imposed in Theorems \ref{Theorem 1} and   \ref
{Theorem 2}, as well as comments concerning situations when these
assumptions are fulfilled.

\begin{remark}
\label{remarklex1} Our proofs of Theorems \ref{Theorem 1} and
\ref%
{Theorem 2} rely, as outlined below, on certain a priori estimates
proved in
Section \ref{sectcomp}. These estimates are proved assuming that
$D\subset
\mathbb{R}
^{d+1}$ is a time-dependent domain satisfying \eqref{timedep+}, %
\eqref{limitzero} and a uniform exterior sphere condition in time with
radius $r_{0}$ in the sense of \eqref{extsphere-}. Furthermore, to derive
these estimates, we also assume that \eqref{crita} and \eqref{crite}
hold for
some $0<\rho_{0}<r_{0}$, $\eta_{0}>0$, $a$ and $e$ and that
$([0,T]\times
\mathbb{R}
^{d})\setminus\overline{D}$ has the $(\delta_{0},h_{0})$-property of good
projections along $\Gamma$, for some $0<\delta_{0}<\rho_{0}$, $h_{0}>1$,
as defined in \eqref{deltanolldef} and \eqref{rhonolldef}. In
particular, we
do not have to assume that $\Gamma=\Gamma_{t} ( z ) $ satisfies %
\eqref{cone1}, \eqref{cone2} and \eqref{Gammanpre} in order to
derive the
results in Section \ref{sectcomp}.
\end{remark}

\begin{remark}
\label{remarklex1+} In Section \ref{sectskorgen} we proceed toward the
final proof of Theorem \ref{Theorem 1}. In particular, we use the a priori
estimates of Section \ref{sectcomp} to derive general results
concerning the
convergence of solutions to Skorohod problems in time-dependent
domains. We
note that our assumptions on $D$ do not exclude the possibility of
holes in $%
D$ and $D_{t}$, for some $t\in\lbrack0,T]$. Nevertheless, the assumptions
on $D$ ensure that the number of holes in $D_{t}$ stays the same for
all $%
t\in [ 0,T ] $ and that these holes cannot shrink too much as time
changes. This observation, Lemma \ref{llhatequivalence} below and its proof
allow us to conclude the validity of the conclusion in Remark \ref{remsim},
which, in turn, is used to complete the proofs in Section \ref{sectskorgen}.
Simple examples show that the conclusion in Remark \ref{remsim} would not
hold if we, for instance, allowed the number of holes in $D_{t}$ to change
as a function of $t$ and if we, in particular, allowed the holes to vanish.
\end{remark}

\begin{remark}
\label{remarklex1++} The assumption that $\Gamma=\Gamma_{t} ( z )
$ satisfies \eqref{cone1}, \eqref{cone2} and \eqref{Gammanpre} is
used to
complete the proofs in Section \ref{sectskorgen}. In particular,
focusing on
Theorem \ref{SPconvergence+}, which is the convergence result actually used
in the proof of Theorem~\ref{Theorem 1}, we note that we need to
assume %
\eqref{Gammanpre} in order to be able derive \eqref{bbbb2again}. We
then use %
\eqref{cone2} to complete the argument in the proof of Theorem \ref%
{SPconvergence+}. Note also the difference between \eqref{cone2} and %
\eqref{Gammanpre}. Assumption \eqref{cone2} simply states that if $%
(s_{n},z_{n},\gamma_{n})$ is a sequence such that $\gamma_{n}\in
\Gamma
_{s_{n}} ( z_{n} ) , z_{n}\in\partial D_{s_{n}},s_{n}\in [
0,T ] $,\vspace*{1pt} and if $(s_{n},z_{n},\gamma_{n})\rightarrow(t,z,\gamma)$ in
$%
\mathbb{R}
\times
\mathbb{R}
^{d}\times
\mathbb{R}
^{d}$, for some $(s,z,\gamma)\in
\mathbb{R}
\times
\mathbb{R}
^{d}\times
\mathbb{R}
^{d}$, then $\gamma\in\Gamma_{s} ( z ) , z\in\partial
D_{s},s\in [ 0,T ] $. Assumption \eqref{Gammanpre}, on the other
hand, is a statement concerning the convergence, in the Hausdorff distance
sense, of the cones $\{\Gamma_{s_{n}}(z_{n})\}$. Finally, to comment on
assumption \eqref{cone1}, which was also imposed in \cite{Costantini1992},
we note that \eqref{cone1} is only used in the proofs of Theorems \ref%
{SPconvergence} and   \ref{SPconvergence+} and, in particular,
in the
verification of \eqref{ddd5} and \eqref{ddd3}. Assumption \eqref{cone1}
eliminates the possibility of $\Gamma$ containing vectors in opposite
directions and we have not been able to complete our argument without this
assumption. However, there are articles dealing with Skorohod type lemmas
and reflected Brownian motion; see \cite{BurdzyToby1995}, in particular,
where this assumption is not required. As noted in \cite{BurdzyToby1995},
the inclusion of vectors in opposite directions can be viewed as a critical
case and \cite{BurdzyToby1995} considers a related problem in a particular
setting in the plane. In our general case we leave this question as a
subject for future research.
\end{remark}

\begin{remark}
\label{remarklex2} For examples of cases when the geometric assumptions
imposed in Theorems \ref{Theorem 1} and   \ref{Theorem 2} are fulfilled,
we refer to   \hyperref[appx]{Appendix}. However, we here briefly discuss
Theorems \ref%
{Theorem 1} and   \ref{Theorem 2} in the context of convex
domains. In
particular, let $T>0$ and let $D\subset
\mathbb{R}
^{d+1}$ be a time-dependent domain satisfying \eqref{timedep+} and %
\eqref{limitzero}. Assume, in addition, that $D_{t}$ is convex
whenever $%
t\in\lbrack0,T]$. Let $\Gamma=\Gamma_{t} ( z ) $ be as in the
statement of Theorem \ref{Theorem 1}. Assume that
%
\begin{eqnarray}\label{critacon}
\lim_{\eta\rightarrow0}\lim_{\rho\rightarrow0}\inf_{s\in [ 0,T%
] }\inf_{z\in\partial D_{s}}a_{s,z} ( \rho,\eta ) &=&a>0,
\\\label{critecon}
\lim_{\eta\rightarrow0}\lim_{\rho\rightarrow0}\sup_{s\in [ 0,T%
] }\sup_{z\in\partial D_{s}}e_{s,z} ( \rho,\eta ) &=&e<1.
\end{eqnarray}
If $D_{t}$ is convex whenever $t\in\lbrack0,T]$, then there exists, for
every $0<\delta_{0}$ given, $h_{0}>1$ such that $([0,T]\times
\mathbb{R}
^{d})\setminus\overline{D}$ has the $(\delta_{0},h_{0})$-property of good
projections along $\Gamma$ as defined in \eqref{deltanolldef} and %
\eqref{rhonolldef}. In this case the conclusion of Theorem~\ref
{Theorem 1}
is, as can be seen from the proofs below, that given $w\in\mathcal{D}
( %
[ 0,T ] ,%
\mathbb{R}
^{d} ) $, with $w_{0}\in\overline{D_{0}}$, there exists a solution $%
( x,\lambda ) $ to the Skorohod problem for $ ( D,\Gamma
,w ) $, in the sense of Definition \ref{skorohodprob}, with $x\in
\mathcal{D} ( [0,T],%
\mathbb{R}
) $. Moreover, if $w$ is a continuous function, then $x$ is
continuous. In particular, if the time-slices $\{D_{t}\}$ are convex, then
the restrictions, in Theorems \ref{Theorem 1} and   \ref{Theorem
2}, on
the jump-sizes in terms of $\delta_{0}$, $\rho_{0}$ can be removed.
Moreover, this is consistent with the results in \cite{Costantini1992} valid
in time-independent domains; see Theorem~4.1 and Proposition~2.3 in
\cite{Costantini1992}.
\end{remark}

We next formulate a subsequent application of Theorems \ref{Theorem 1} and
  \ref{Theorem 2} to the problem of constructing weak solutions to
stochastic differential equations\vspace*{1pt} in $\overline{D}$ with reflection
along $%
\Gamma_{t}$ on $\partial D_{t}$ for all $t\in [ 0,T ] $. Given $%
T>0$, we let $\mathcal{C} ( [ 0,T ] ,%
\mathbb{R}
^{d} ) $ denote the class of continuous functions from $ [ 0,T%
] $ to $%
\mathbb{R}
^{d}$. In the following, we let $m$ be a positive integer and we let
$b\dvtx
\mathbb{R}
_{+}\times
\mathbb{R}
^{d}\rightarrow
\mathbb{R}
^{d}$ and $\sigma\dvtx %
\mathbb{R}
_{+}\times
\mathbb{R}
^{d}\rightarrow
\mathbb{R}
^{d\times m}$ be given functions which are bounded and continuous.

\begin{definition}
\label{weaksol} Let $d\geq1$ and $T>0$. Let $D\subset
\mathbb{R}
^{d+1}$ be a time-dependent domain satisfying \eqref{timedep+}, let
$\Gamma
=\Gamma_{t} ( z ) $ be a closed convex cone of vectors in $%
\mathbb{R}
^{d}$ for every $z\in\partial D_{t}$, $t\in\lbrack0,T]$, and let
$\hat{z}%
\in\overline{D_{0}}$. A weak solution to the stochastic differential
equation in $\overline{D}$ with coefficients $b$ and $\sigma$, reflection\vspace*{1pt}
along $\Gamma_{t}$ on $\partial D_{t}$, $t\in [ 0,T ] $, and with
initial condition $\hat{z}$ at $t=0$, is a stochastic process $ ( X^{0,
\hat{z}},\Lambda^{0,\hat{z}} ) $ with paths in $\mathcal{C} ( %
[ 0,T ] ,%
\mathbb{R}
^{d} ) \times\mathcal{BV} ( [ 0,T ] ,%
\mathbb{R}
^{d} ) $, which is defined on a filtered probability space $ (
\Omega,\mathcal{F,} \{ \mathcal{F}_{t} \} ,P ) $ and
satisfies, $P$-almost surely, whenever $t\in [ 0,T ] $,
%
\begin{eqnarray}\label{RSDEth1}
X_{t}^{0,\hat{z}}
&=&\hat{z}+\int_{0}^{t}b ( s,X_{s}^{0,\hat{z}} )
\,ds+\int_{0}^{t}\sigma ( s,X_{s}^{0,\hat{z}} ) \,dW_{s}+\Lambda
_{t}^{0,\hat{z}}, \\\label{RSDEth1+}
\Lambda_{t}^{0,\hat{z}}
&=&\int_{0}^{t}\gamma_{s}\,d|\Lambda
^{0,\hat{z}%
}|_{s},\qquad \gamma_{s}\in\Gamma_{s} ( X_{s}^{0,\hat{z}} ) \cap
S_{1}(0),\  d|\Lambda^{0,\hat{z}}|\mbox{-a.e.}, \\ \label{RSDEth2}
X_{t}^{0,\hat{z}}
&\in&\overline{D_{t}},\qquad d \vert\Lambda^{0,\hat{z%
}} \vert ( \{ t\in [ 0,T ] \dvtx X_{t}^{0,\hat{z}}\in
D_{t} \} ) =0.
\end{eqnarray}
Here $W$ is a $m$-dimensional Wiener process on $ ( \Omega,\mathcal
{F,}%
\{ \mathcal{F}_{t} \} ,P ) $ and $ ( X^{0,\hat{z}%
},\Lambda^{0,\hat{z}} ) $ is $ \{ \mathcal{F}_{t} \} $%
-adapted.
\end{definition}

Concerning weak solutions to stochastic differential equations in
$\overline{%
D}$ with oblique reflection along $\partial D$, we prove the following
theorem.

\begin{theorem}
\label{Theorem 3}Let $T>0$, $D\subset
\mathbb{R}
^{d+1}$ and $\Gamma=\Gamma_{t} ( z ) $ be as in the statement of
Theorem \ref{Theorem 1}. Let $b\dvtx %
\mathbb{R}
_{+}\times
\mathbb{R}
^{d}\rightarrow
\mathbb{R}
^{d}$ and $\sigma\dvtx %
\mathbb{R}
_{+}\times
\mathbb{R}
^{d}\rightarrow
\mathbb{R}
^{d\times m}$ be given, bounded and continuous functions on $\overline{D}$
and let $\hat z\in\overline{D_{0}}$. Then there exists a weak
solution, in
the sense of Definition \ref{weaksol}, to the stochastic differential
equation in $\overline{D}$ with coefficients $b$ and $\sigma$, reflection
along $\Gamma_{t}$ on $\partial D_{t}$, $t\in [ 0,T ] $, and with
initial condition~$\hat z$ at $t=0$.
\end{theorem}

We note that Theorem \ref{Theorem 3} generalizes the corresponding results
in \cite{Costantini1992,CostantiniGobetKaroui2006} and \cite{Saisho1987}. Furthermore, we note that there has recently been considerable
activity in the study of reflected diffusions in time-dependent intervals.
In particular, in this context we mention \cite{BurdzyChenSylvester2003,BurdzyChenSylvester2004AP,BurdzyChenSylvester2004JMAA} and
\cite{BurdzyKangRamanan2009} and we refer the interested reader to these
articles for more information as well as for references to other related
articles.

The rest of the article is organized as follows. In Section \ref
{contri} we
first briefly outline two general and important themes present in the proofs
of the results in this article. The first theme concerns a priori estimates
and compactness for solutions to Skorohod problems and the second theme
concerns convergence results for sequences of solutions to Skorohod
problems. Second, we discuss the proofs of Theorems~\ref{Theorem 1},
\ref{Theorem 2} and   \ref{Theorem 3} and we try to point out
the new
difficulties occurring due to the time-dependent character of the domain.
This section is included for further reference and, in particular, to convey
some of the ideas to the reader. In Section~\ref{Prel} we introduce
additional notation, outline the restrictions imposed on $D$ and~$\Gamma$
and collect a few notions and facts from the Skorohod topology. There is
also an appendix attached to Section \ref{Prel}, \hyperref[appx]{Appendix}. In
\hyperref[appx]{Appendix} we state sufficient conditions for the $(\delta
_{0},h_{0})$-property of good projections along  $\Gamma$ and we give
examples of time-dependent domains satisfying the assumptions stated in
Theorems \ref{Theorem 1},   \ref{Theorem 2} and   \ref
{Theorem 3}%
. Section \ref{sectcomp} is devoted to estimates for solutions to the
Skorohod problem, with oblique reflection, which have bounded jumps and also
to the corresponding estimates for certain approximations of the Skorohod
problem. In Section \ref{sectskorgen} we first prove Theorem \ref%
{SPconvergence}, containing a general result concerning convergence of
solutions to Skorohod problems in time-dependent domains. Furthermore, we
establish the somewhat similar result for certain approximations of the
Skorohod problem. The latter estimates are then used in the proof of
Theorem %
\ref{Theorem 1}. The final proofs of Theorems \ref{Theorem 1},
\ref%
{Theorem 2} and   \ref{Theorem 3} are given in Section \ref{sectskor}.
The article ends with the  \hyperref[appx]{Appendix}, discussed above.

\section{A brief outline of proofs and our contribution}\label{contri}

Concerning proofs, we note that the arguments in this article follow two
general and important themes which we here, to start with, briefly
outline.

\subsection*{A priori estimates and compactness}

To explain the a
priori estimates, we let $T>0$, $D\subset
\mathbb{R}
^{d+1}$ and $\Gamma=\Gamma_{t} ( z ) $ be as in the statement of
Theorem \ref{Theorem 1}, and we let $w\in\mathcal{D} ( [ 0,T ]
,%
\mathbb{R}
^{d} ) $ with $w_{0}\in\overline{D_{0}}$. Assume that $ (
x,\lambda ) $ is a solution to the Skorohod problem for $ (
D,\Gamma,w ) $ such that $x\in\mathcal{D}^{\rho_{0}} ( [ 0,T%
] ,%
\mathbb{R}
^{d} ) $. Under these assumptions, we prove (see Theorem \ref%
{compacttheorem} below) that there exist positive constants $L_{1} (
w,T ) $, $L_{2} ( w,T ) $, $L_{3} ( w,T ) $ and $%
L_{4} ( w,T ) $ such that
%
\begin{eqnarray}\label{apriori}
\Vert x \Vert_{t_{1},t_{2}} &\leq&L_{1} ( w,T )
\Vert w \Vert_{t_{1},t_{2}}+L_{2} ( w,T ) l (
t_{2}-t_{1} ) ,
\nonumber
\\[-8pt]
\\[-8pt]
|\lambda |_{t_{2}}-|\lambda |_{t_{1}} &\leq &L_{3} ( w,T )
 \Vert w \Vert _{t_{1},t_{2}}+L_{4} ( w,T ) l (
t_{2}-t_{1} ) ,\nonumber
\end{eqnarray}
whenever $0\leq t_{1}\leq t_{2}\leq T$. Furthermore, we prove that if $%
\mathcal{W}\subset\mathcal{D} ( [ 0,T ] ,%
\mathbb{R}
^{d} ) $ is relatively compact in the Skorohod topology and $w_{0}\in
\overline{D_{0}}$, whenever $w\in\mathcal{W}$, then there exist positive
constants $L_{1}^{T}$, $L_{2}^{T}$, $L_{3}^{T}$ and $L_{4}^{T}$, such
that%
%
\begin{equation}
\sup_{w\in\mathcal{W}}L_{i} ( w,T ) \leq L_{i}^{T}<\infty \qquad
\mbox{for }i=1,2,3,4. \label{apriori2}
\end{equation}

\subsection*{Convergence results for sequences of solutions to Skorohod
problems}

The a priori estimates and compactness result in \eqref{apriori}
and \eqref{apriori2} are useful for proving convergence of solutions to
Skorohod problems. To explain this further, let $\{D^{n}\}
_{n=1}^{\infty}$
be a sequence of time-dependent domains $D^{n}\subset
\mathbb{R}
^{d+1}$ satisfying \eqref{timedep+} and let $\{\Gamma^{n}\}
_{n=1}^{\infty
}=\{\Gamma_{t}^{n} ( z ) \}_{n=1}^{\infty}$ be a sequence of
closed convex cones of vectors in $%
\mathbb{R}
^{d}$. Assume that $\{D^{n}\}_{n=1}^{\infty}$ and $\{\Gamma
^{n}\}_{n=1}^{\infty}$ satisfy the conditions stated in Theorem \ref%
{Theorem 1} with constants that are ``uniform with respect to $n$'' in a sense
made precise in Section \ref{sectskorgen}. Let $ \{ w^{n} \} $,
with $w_{0}^{n}\in\overline{D_{0}^{n}}$, be a sequence in $\mathcal
{D}%
( [ 0,T ] ,%
\mathbb{R}
^{d} ) $ which is relatively compact in the Skorohod topology and which
converges to $w\in\mathcal{D} ( [ 0,T ] ,%
\mathbb{R}
^{d} ) $ with $w_{0}\in\overline{D_{0}}$. Furthermore, let $D\subset
\mathbb{R}
^{d+1}$ be a time-dependent domain satisfying \eqref{timedep+}, let
$\Gamma
=\Gamma_{t}(z)$ be, for every $z\in\partial D_{t}$, $t\in\lbrack
0,T]$, a
closed convex cone of vectors in $%
\mathbb{R}
^{d}$ satisfying \eqref{cone1} and \eqref{cone2}. Assume that the
sequences $%
\{D^{n}\}_{n=1}^{\infty}$ and $\{\Gamma^{n}\}_{n=1}^{\infty}$
converge to
$D$ and $\Gamma$, respectively, in a sense specified in Theorem~\ref%
{SPconvergence}. If there exists, for all $n\geq1$, a solution $ (
x^{n},\lambda^{n} ) $ to the Skorohod problem for $ ( D^{n},\Gamma
^{n},w^{n} ) $ such that $x_{t}^{n}\in\overline{D_{t}^{n}}$, for all $
t\in [ 0,T ] $, and $x^{n}\in\mathcal{D}^{\rho_{0}} ( [
0,T ] ,%
\mathbb{R}
^{d} ) $, then it follows, using \eqref{apriori} and \eqref{apriori2},
that $ \{ ( w^{n},x^{n},\lambda^{n}, \vert\lambda
^{n} \vert ) \} $ is relatively compact in $\mathcal{D}%
( [ 0,T ] ,%
\mathbb{R}
^{d} ) \times\mathcal{D} ( [ 0,T ] ,%
\mathbb{R}
^{d} ) \times\mathcal{D} ( [ 0,T ] ,%
\mathbb{R}
^{d} ) \times\mathcal{D} ( [ 0,T ] ,%
\mathbb{R}
_{+} ) $. Hence, we are able to conclude that $ \{ (
x^{n},\lambda^{n} ) \} $ converges to some $(x,\lambda)\in
\mathcal{D} ( [ 0,T ] ,%
\mathbb{R}
^{d} ) \times\mathcal{D} ( [ 0,T ] ,%
\mathbb{R}
^{d} ) $ with $x\in\overline{D}$ and we can, in addition, prove that $
(x,\lambda)$ is indeed a solution to the Skorohod problem for $ (
D,\Gamma,w ) $. This result, found in Theorem \ref{SPconvergence}
below, constitutes a general convergence result for sequences of solutions
to Skorohod problems based on the a priori estimates and compactness result
in \eqref{apriori} and~\eqref{apriori2}.

Although Theorem \ref{Theorem 1} does not follow directly from the results
outlined above, we claim that \eqref{apriori}, \eqref{apriori2} and
Theorem %
\ref{SPconvergence}, stated below, are of independent interest and may be
useful in other applications involving the Skorohod problem. To start an
outline of the actual proofs of Theorems \ref{Theorem 1},   \ref%
{Theorem 2} and   \ref{Theorem 3}, we note that to prove Theorem
\ref%
{Theorem 1} we use arguments similar to those outlined above, but in this
case we have to construct, given $(D,\Gamma,w)$, an approximating
sequence $%
\{(D^{n},\Gamma^{n},w^{n})\}$ such that a solution $(x^{n},\lambda^{n})$
to the Skorohod problem for $(D^{n},\Gamma^{n},w^{n})$ can be found
explicitly.

\subsection*{Proof of Theorems \protect\ref{Theorem 1}, \protect
\ref{Theorem 2}
and \protect\ref{Theorem 3}}

To discuss the construction of $%
\{(D^{n},\break\Gamma^{n}, w^{n})\}$ and $\{(x^{n},\lambda^{n})\}$ used in the
proof of Theorem \ref{Theorem 1}, we consider $w\in\mathcal{D} ( [
0,T ] ,%
\mathbb{R}
^{d} ) $, with $w_{0}\in\overline{D_{0}}$ and with jumps bounded by
some constant, and we now let $\{\tau_{k}\}_{k=0}^{N}$ define a
partition $%
\Delta$ of the interval $[0,T]$, that is, $0=\tau_{0}<\tau
_{1}<\cdots<\tau
_{N-1}<\tau_{N}=T$. Given $\Delta$, we let
%
\begin{equation}
\Delta^{\ast}:=\max_{k\in\{0,\ldots,N-1\}}\tau_{k+1}-\tau_{k},
\label{korlim+}
\end{equation}
and, given $\Delta$ and $w$, we define
%
\begin{equation}
w_{t}^{\Delta}=w_{\tau_{k-1}} \qquad \mbox{whenever }t\in [ \tau
_{k-1},\tau_{k} ) , k\in \{ 1,\ldots,N \} ,
\label{approxw}
\end{equation}
and $w_{T}^{\Delta}=w_{T}$. Then $w^{\Delta}\in\mathcal{D} ( [
0,T ] ,%
\mathbb{R}
^{d} ) $ is a step function approximation of $w$. Furthermore, assume
that $\Delta$ and $w^{\Delta}$ are such that
%
\begin{equation}
\Vert w^{\Delta}\Vert_{\tau_{k-1},\tau_{k}}+l(\Delta^{\ast})<\delta_{0},
\label{korlim}
\end{equation}
whenever $k\in\{1,\ldots,N\}$. Recall that $\delta_0$ is the constant
appearing in the notion of the $(\delta_{0},h_{0})$-property of good
projections. We next define
%
\begin{eqnarray}\label{approxw+}
D_{t}^{\Delta}&=&D_{\tau_{k-1}},\nonumber
\\[-8pt]
\\[-8pt]
  \Gamma_{t}^{\Delta}&=&\Gamma_{\tau
_{k-1}}%
  \qquad\mbox{whenever }t\in [ \tau_{k-1},\tau_{k} ) , k\in
\{ 1,\ldots,N \} ,\nonumber
\end{eqnarray}
and $D_{T}^{\Delta}=D_{T}, \Gamma_{T}^{\Delta}=\Gamma_{T}$. Given $%
w^{\Delta}$, $D^{\Delta}$ and $\Gamma^{\Delta}$ as above, we define a
pair of processes $ ( x^{\Delta},\lambda^{\Delta} ) $ as
follows. We let%
%
\begin{equation}
x_{t}^{\Delta}=w_{0},\qquad \lambda_{t}^{\Delta}=0 \qquad\mbox{for }t\in [
0,\tau_{1} ). \label{approxw++}
\end{equation}
If $x_{\tau_{k-1}}^{\Delta}\in\overline{D_{\tau_{k-1}}^{\Delta}}$ for
some $k\in \{ 1,\ldots,N \} $, then, by the triangle inequality and %
\eqref{korlim},
%
\begin{equation}
d ( x_{\tau_{k-1}}^{\Delta}+w_{\tau_{k}}^{\Delta}-w_{\tau
_{k-1}}^{\Delta},\overline{D_{\tau_{k}}^{\Delta}} ) \leq
\Vert w^{n}\Vert _{\tau_{k-1},\tau_{k}}+l(\Delta^{\ast})<\delta_{0}.
\label{jjlex}
\end{equation}
Hence, by the $ ( \delta_{0},h_{0} ) $-property of good
projections, it follows that if $x_{\tau_{k-1}}^{\Delta}+w_{\tau
_{k}}^{\Delta}-w_{\tau_{k-1}}^{\Delta}\notin\overline{D_{\tau
_{k}}^{\Delta}%
}$, then there exists a point
%
\begin{equation}
\pi_{\partial D_{\tau_{k}}^{\Delta}}^{\Gamma_{\tau_{k}}^{\Delta
}} ( x_{\tau_{k-1}}^{\Delta}+w_{\tau_{k}}^{\Delta}-w_{\tau
_{k-1}}^{\Delta} )\in\partial D_{\tau_{k}}^{\Delta}, \label{jjjlex}
\end{equation}
which is the projection of $x_{\tau_{k-1}}^{\Delta}+w_{\tau
_{k}}^{\Delta}-w_{\tau_{k-1}}^{\Delta}$ onto $\partial D_{\tau
_{k}}^{\Delta} $ along $\Gamma_{\tau_{k}}^\Delta$.\vspace*{1pt} Furthermore, if $
x_{\tau_{k-1}}^{\Delta}+w_{\tau_{k}}^{\Delta}-w_{\tau
_{k-1}}^{\Delta}\in%
\overline{D_{\tau_{k}}^{\Delta}}$, then we let
%
\begin{equation}
\pi_{\partial D_{\tau_{k}}^{\Delta}}^{\Gamma_{\tau_{k}}^{\Delta
}} ( x_{\tau_{k-1}}^{\Delta}+w_{\tau_{k}}^{\Delta}-w_{\tau
_{k-1}}^{\Delta} )=x_{\tau_{k-1}}^{\Delta}+w_{\tau
_{k}}^{\Delta}-w_{\tau_{k-1}}^{\Delta}. \label{jjjjlex}
\end{equation}
Based on this argument, we define, whenever $t\in [ \tau_{k},\tau
_{k+1} ) $, $k\in \{ 1,\ldots,N-1 \} $,
%
\begin{eqnarray}\label{approxw+++}
x_{t}^{\Delta} &=&\pi_{\partial D_{\tau_{k}}^{\Delta}}^{\Gamma
_{\tau
_{k}}^{\Delta}} ( x_{\tau_{k-1}}^{\Delta}+w_{\tau_{k}}^{\Delta
}-w_{\tau_{k-1}}^{\Delta} ) ,\nonumber
\\[-8pt]
\\[-8pt]
\lambda_{t}^{\Delta} &=&\lambda_{\tau_{k-1}}^{\Delta
}+\bigl(x_{t}^{\Delta
}-(x_{\tau_{k-1}}^{\Delta}+w_{\tau_{k}}^{\Delta}-w_{\tau
_{k-1}}^{\Delta
})\bigr).\nonumber
\end{eqnarray}
Finally, we define $x_{T}^{\Delta}$ and $\lambda_{T}^{\Delta}$ as in
\eqref{approxw+++} by putting $k=N$ in \eqref{approxw+++}. By construction,
the pair $(x^{\Delta},\lambda^{\Delta})$ is a solution to the Skorohod
problem for $ ( D^{\Delta},\Gamma^{\Delta},w^{\Delta} )$.
Moreover, using the assumption on the size of the jumps of $w$ stated in
Theorem \ref{Theorem 1}, we will be able to make the construction so
that we
can conclude that $x^{\Delta}\in\mathcal{D}^{\rho_{0}} ( [ 0,T%
] ,%
\mathbb{R}
^{d} ) $. As the next step we then apply Theorem \ref%
{compacttheoremapprox} stated below, showing the existence of positive
constants $\hat{L}_{1} ( w,T ) $, $\hat{L}_{2} ( w,T ) $, $%
\hat{L}_{3} ( w,T ) $ and $\hat{L}_{4} ( w,T ) $ such that
%
\begin{eqnarray}\label{lll1}\qquad
\Vert x^{\Delta} \Vert_{t_{1},t_{2}} &\leq&\hat{L}_{1} (
w,T ) \Vert w \Vert_{t_{1},t_{2}}+\hat{L}_{2} (
w,T ) \bigl(l ( t_{2}-t_{1} ) +l(\Delta^{\ast})\bigr),\nonumber
\\[-8pt]
\\[-8pt]\qquad
|\lambda^{\Delta}|_{t_{2}}-|\lambda^{\Delta}|_{t_{1}} &\leq&\hat
{L}%
_{3} ( w,T ) \Vert w \Vert_{t_{1},t_{2}}+\hat{L}%
_{4} ( w,T ) \bigl(l ( t_{2}-t_{1} ) +l(\Delta^{\ast})\bigr),\nonumber
\end{eqnarray}
whenever $0\leq t_{1}\leq t_{2}\leq T$. Provided with the estimates in %
\eqref{lll1}, we are then able to prove Theorem \ref{Theorem 1} by
means of
compactness arguments similar to those outlined above. Indeed, we construct
an appropriate sequence of partitions $\{\Delta_{n}\}_{n=1}^{\infty}$,
based on $w$, such that $x^{\Delta_{n}}\in\mathcal{D}^{\rho_{0}} ( %
[ 0,T ] ,%
\mathbb{R}
^{d} ) $ for $n$ larger than some $n_{0}$ and such that $(x^{\Delta
_{n}},\lambda^{\Delta_{n}})$ is a solution to the Skorohod problem
for $%
( D^{\Delta_{n}},\Gamma^{\Delta_{n}},w^{\Delta_{n}} ) $ for $%
n\geq n_{0}$. Then, using \eqref{lll1}, we conclude that $ \{ (
w^{\Delta_{n}},x^{\Delta_{n}},\lambda^{\Delta_{n}},|\lambda
^{\Delta
_{n}}| ) \} $ is a relatively compact sequence in the Skorohod
topology and that $ \{ ( x^{\Delta_{n}},\lambda^{\Delta
_{n}} ) \} $ converges in the sense of the Skorohod topology to a
pair of functions $ ( x,\lambda ) $. Note that an important
difference here, compared to the situation outlined above, is that $%
D^{\Delta}$ and $\Gamma^{\Delta}$ as defined in \eqref{approxw+} are
discontinuous in time. To be able to handle this situation, we employ some
additional arguments, similar to the ones in the proof of Theorem \ref%
{SPconvergence}, in order to prove that $(x,\lambda)$ is a solution to the
Skorohod problem for $ ( D,\Gamma,w ) $ on $ [ 0,T ] $.
This completes the proof of Theorem \ref{Theorem 1}. Concerning
Theorem \ref%
{Theorem 2}, we see that this theorem follows immediately from the continuity
of $w$ and \eqref{limitzero} using the estimates in \eqref{apriori}. To
prove Theorem \ref{Theorem 3}, we argue somewhat similarly as in the
proof of
Theorem \ref{Theorem 1} and we refer to the bulk of the article for
details.

To conclude, we note that the proof of Theorem \ref{Theorem 1} is more
involved compared to the proof of the corresponding result for
time-independent domains established in \cite{Costantini1992} and that new
difficulties occur, naturally, due to the fact that we are considering
time-dependent domains. In the time-independent case a solution
$(x,\lambda
) $ to the Skorohod problem for $(D,\Gamma,w)$ is constructed as the limit
of a sequence $\{(x^{\Delta_{n}},\lambda^{\Delta_{n}})\}$, where $%
(x^{\Delta_{n}},\lambda^{\Delta_{n}})$ is a solution to a Skorohod
problem based on $w^{\Delta_{n}}$. In this case $(x^{\Delta
_{n}},\lambda
^{\Delta_{n}})$ is a solution to a Skorohod problem for $(D,\Gamma
,w^{\Delta_{n}})$, while in our case $(x^{\Delta_{n}},\lambda
^{\Delta
_{n}})$ is a solution to a Skorohod problem for $(D^{\Delta
_{n}},\Gamma
^{\Delta_{n}},w^{\Delta_{n}}) $. Hence, in the time-dependent case
we, at
each step, also have to discretize and approximate $D$ and $\Gamma$
due to
the time-dependent character of the domain. In particular, the fact
that $%
D^{\Delta_{n}}$ and $\Gamma^{\Delta_{n}}$, as defined in \eqref
{approxw+}%
, are discontinuous in time induces several new difficulties which we have
to overcome in order to complete the proof of Theorem \ref{Theorem
1}.

\section{Preliminaries}\label{Prel}

In this section we introduce notation, collect a number of preliminary
results concerning the geometry of time-dependent domains and recall a few
notions and facts from the Skorohod topology.

\subsection{Notation}

Points in Euclidean $ ( d+1 ) $-space $%
\mathbb{R}
^{d+1}$ are denoted by $ ( t,z ) = ( t,z_{1},\ldots
,z_{d} ) $. Given a differentiable function $f=f ( t,z ) $
defined on $%
\mathbb{R}
\times
\mathbb{R}
^{d}$, we let $\partial_{z_{i}}f ( t,z ) $ denote the partial
derivative of $f$ at $(t,z)$ with respect to $z_{i}$ and we let $\nabla_{z}f$
denote the gradient $ ( \partial_{z_{1}}f,\ldots,\partial_{z_{d}}f )
$. Higher order derivatives of $f$ with respect to the space variables will
often be denoted by $\partial_{z_{i}z_{j}}f ( t,z ) $, $\partial
_{z_{i}z_{j}z_{k}}f ( t,z ) $ and so on. Furthermore, given a
multi-index $\beta= ( \beta_{1},\ldots,\beta_{d} ) $, $\beta
_{i}\in
\mathbb{Z}
_{+}$, we define $ \vert\beta \vert=\beta_{1}+\cdots+\beta_{d}$
and we let $\partial_{z}^{\beta}f ( t,z ) $ denote the associated
partial derivative of $f ( t,z ) $ with respect to the space
variables. Time derivatives of $f$ will be denoted by $\partial
_{t}^{j}f ( t,z ) $ where $j\in
\mathbb{Z}
_{+}$. As in the \hyperref[intro]{Introduction}, we let $ \langle\cdot,\cdot
\rangle$ denote the standard inner product on $%
\mathbb{R}
^{d}$ and we let $ \vert z \vert= \langle z,z \rangle
^{1/2}$ be the Euclidean norm of $z$. Whenever $z\in
\mathbb{R}
^{d}, r>0$, we let $B_{r} ( z ) = \{ y\in
\mathbb{R}
^{d}\dvtx  \vert z-y \vert<r \} $ and $S_{r} ( z )
= \{ y\in
\mathbb{R}
^{d}\dvtx  \vert z-y \vert=r \} $. In addition, $dz$ denotes
the Lebesgue $d$-measure on $%
\mathbb{R}
^{d}$. Moreover, given $E\subset
\mathbb{R}
^{d}$, we let $\bar{E}$ and $\partial E$ be the closure and boundary
of~$E$,
respectively, and we let $d ( z,E ) $ denote the Euclidean distance
from $z\in
\mathbb{R}
^{d}$ to $E$. Given $ ( t,z ) , ( s,y ) \in
\mathbb{R}
^{d+1}$, we let $d_{p} ( ( t,z ) , ( s,y ) )
=\max \{ \vert z-y \vert, \vert t-s \vert
^{1/2} \} $ denote the parabolic distance between $ ( t,z ) $
and $ ( s,y ) $ and for $F\subset
\mathbb{R}
^{d+1}$, we let $d_{p} ( ( t,z ) ,F ) $ denote the
parabolic distance from $ ( t,z ) \in
\mathbb{R}
^{d+1}$ to $F$. Moreover, for $ ( t,z ) \in
\mathbb{R}
^{d+1}$ and $r>0$, we introduce the parabolic cylinder $C_{r} (
t,z ) = \{ ( s,y ) \in
\mathbb{R}
^{d+1}\dvtx  \vert y-z \vert<r, \vert t-s \vert<r^{2} \}
$. Given two real numbers $a$ and $b$, we let $a\vee b=\max \{
a,b \} $ and $a\wedge b=\min \{ a,b \} $. Finally, given a
Borel set $E\subset
\mathbb{R}
^{d+1}$, we let $\chi_{E}$ denote the characteristic function
associated to
$E$.

Given a time-dependent domain $D^{\prime}$, a function $f$ defined on $
D^{\prime}$ and a constant $\alpha\in ( 0,1 ] $, we adopt the
definition on page 46 in \cite{Lieberman1996} and introduce
%
\begin{eqnarray}
\vert f \vert_{1+\alpha,D^{\prime}} &=&\sum_{ \vert\beta
\vert\leq1}\sup_{D^{\prime}} \vert\partial_{z}^{\beta
}f \vert+\sup_{ ( t,z ) \in D^{\prime
}}\sup_{ ( s,z ) \in D^{\prime}\setminus \{ (
t,z ) \} }\frac{ \vert f ( t,z ) -f (
s,z ) \vert}{ \vert t-s \vert^{ ( \alpha+1 )
/2}} \nonumber
\\[-8pt]
\\[-8pt]
&&{}+\sum_{ \vert\beta \vert=1}\sup_{ ( t,z ) \in
D^{\prime}}\sup_{ ( s,y ) \in D^{\prime}\setminus \{ (
t,z ) \} }\frac{ \vert\partial_{z}^{\beta}f (
t,z ) -\partial_{z}^{\beta}f ( s,y ) \vert}{ [
d_{p} ( ( t,z ) , ( s,y ) ) ] ^{\alpha}}.\nonumber
\end{eqnarray}
The third term on the right-hand side of $ \vert f \vert_{1+\alpha
,D^{\prime}}$ is superfluous for our purposes, but we include it here for
agreement with the theory of partial differential equations in
time-dependent domains (see \cite{NystromOnskog2009c}). Using the norm
$%
\vert f \vert_{1+\alpha,D^{\prime}}$, we let $\mathcal{H}%
_{1+\alpha} ( D^{\prime} ) $ denote the Banach space of functions
$f$ on $D^{\prime}$ with finite $ \vert f \vert_{1+\alpha
,D^{\prime}}$-norm.

\subsection{Geometry of time-dependent domains}\label{geometry}

We here outline the geometric restrictions which we impose on the
time-dependent domains and cones of reflections. Concerning $D$, we first
prove the following auxiliary lemma.

\begin{lemma}
\label{llhatequivalence} Let $T>0$ and let $D\subset
\mathbb{R}
^{d+1}$ be a time-dependent domain satisfying \eqref{timedep+} and %
\eqref{limitzero} and assume that $D$ satisfies a uniform exterior sphere
condition in time with radius $r_{0}$ in the sense of \eqref
{extsphere-}. Let%
%
\begin{equation}
\hat{l} ( r ) :=\mathop{\mathop{\sup}_{ s,t\in[ 0,T]}}_{ \vert
s-t \vert\leq r}\sup_{z\in\partial D_{s}}d ( z,\partial
D_{t} ) . \label{lhat}
\end{equation}
Then, $l ( r ) =\hat{l} ( r ) $ for all $r>0$ such that $%
l ( r ) <r_{0}$.
\end{lemma}

\begin{pf} In the following we let $\epsilon>0$ be arbitrary
and we consider $r$ small enough to ensure that $l(r)<r_{0}$. With
$\epsilon
$ and $r$ fixed, we let $s,t\in\lbrack0,T]$, $ \vert s-t \vert
\leq r$, be such that
%
\begin{equation}
\tilde{l}_{1}\leq\hat{l} ( r ) \leq\tilde{l}_{1}+\epsilon\qquad\mbox{where }\tilde{l}_{1}=\sup_{z\in\partial D_{s}}d ( z,\partial
D_{t} ) . \label{tildel}
\end{equation}
Naturally,
%
\begin{equation}
\tilde{l}_{1}=\max \Bigl\{\sup_{z\in\partial D_{s}\cap(%
\mathbb{R}
^{d}\setminus D_{t})}d ( z,\partial D_{t} ) ,\sup_{z\in\partial
D_{s}\cap D_{t}}d ( z,\partial D_{t} ) \Bigr\}. \label{lhat1}
\end{equation}
Assume $z\in\partial D_{s}\cap(%
\mathbb{R}
^{d}\setminus D_{t})$. Then we immediately obtain
%
\begin{equation}\qquad
\sup_{z\in\partial D_{s}\cap(%
\mathbb{R}
^{d}\setminus D_{t})}d ( z,\partial D_{t} ) =\sup_{z\in\partial
D_{s}\cap(%
\mathbb{R}
^{d}\setminus D_{t})}d ( z,D_{t} ) \leq\sup_{z\in\overline{D_{s}}%
}d ( z,D_{t} ) \leq l ( r ) . \label{lhat2}
\end{equation}
Assume, on the contrary, that $z\in\partial D_{s}\cap D_{t}$. In this case,
as $l ( r ) <r_{0}$ and $D_{s}$ satisfies the uniform exterior
sphere condition with radius $r_{0}$, we can conclude that there exists at
least one point $y_{z}\in\partial D_{t}\cap \{ z+n_{\lambda}\in
\mathbb{R}
^{d}:n_{\lambda}\in N_{s} ( z ) \cap S_{\lambda} ( 0 )
,0<\lambda<r_{0} \} $ and obviously $d ( z,\partial D_{t} )
\leq \vert z-y_{z} \vert$. Furthermore, again applying the
uniform exterior sphere condition, we see that $z$ minimizes the distance
from $y_{z}\in\partial D_{t}$ to $D_{s}$. Hence,%
%
\begin{eqnarray}\label{lhat3}
\sup_{z\in\partial D_{s}\cap D_{t}}d ( z,\partial D_{t} ) &\leq
&\sup_{z\in\partial D_{s}\cap D_{t}} \vert z-y_{z} \vert\leq
\sup_{z\in\partial D_{s}\cap D_{t}}d ( y_{z},D_{s} ) \nonumber
\\[-8pt]
\\[-8pt]
&\leq&\sup_{y_{z}\in\partial D_{t}}d ( y_{z},D_{s} ) \leq
\sup_{y\in\overline{D_{t}}}d ( y,D_{s} ) \leq l ( r ) .\nonumber
\end{eqnarray}
Combining \eqref{lhat1}--\eqref{lhat3}, we conclude that $\tilde
{l}_{1}\leq
l ( r ) $. Using \eqref{tildel}, it is clear that $\hat{l} (
r ) \leq l ( r ) +\epsilon$ and then, as $\epsilon$ is
arbitrary, $\hat{l} ( r ) \leq l ( r ) $. We next consider
the opposite inequality. With $\epsilon$ and $r$ fixed, we let $s,t\in
\lbrack0,T]$, $ \vert s-t \vert\leq r$, be such that
%
\begin{equation}
\tilde{l}_{2}\leq{l} ( r ) \leq\tilde{l}_{2}+\epsilon\qquad
\mbox{where }\tilde{l}_{2}=\sup_{z\in\overline{D_{s}}}d ( z,D_{t} ) .
\label{tildel+}
\end{equation}
In this case we have
%
\begin{equation}\qquad
\tilde{l}_{2}=\max \Bigl\{\sup_{z\in\overline{D_{s}}\cap(%
\mathbb{R}
^{d}\setminus D_{t})}d ( z,D_{t} ) ,0 \Bigr\}=\max \Bigl\{%
\sup_{z\in\overline{D_{s}}\cap(%
\mathbb{R}
^{d}\setminus D_{t})}d ( z,\partial D_{t} ) ,0 \Bigr\}
\label{lhat1+}
\end{equation}
and in the following we can assume, without loss of generality, that
$\tilde{%
l}_{2}>0$. Then, by the uniform exterior sphere condition, and the fact that
$l ( r ) <r_{0}$, we see\vspace*{1pt} that every point $z\in\overline{D_{s}}%
\cap(%
\mathbb{R}
^{d}\setminus D_{t})$ can be written as $z=y_{z}+n_{\lambda}$ for some
$%
y_{z}\in\partial D_{t}$ and some $n_{\lambda}\in N_{t} ( y_{z} )
\cap S_{\lambda} ( 0 ) $, $0<\lambda<r_{0}$. Furthermore, there
exists a point $\tilde{z}=y_{z}+n_{\tilde{\lambda}}\in\partial
D_{s}$, with
$0<\lambda<\tilde{\lambda}<r_{0}$. Once again applying the uniform exterior
sphere condition, we see that $y_{z}$ minimizes the distance from
$\tilde{z}$
to $\partial D_{t}$ and we obtain
%
\begin{eqnarray}\label{lhat5}\qquad
\sup_{z\in\overline{D_{s}}\cap(%
\mathbb{R}
^{d}\setminus D_{t})}d ( z,\partial D_{t} ) &\leq&\sup_{z\in
\overline{D_{s}}\cap(%
\mathbb{R}
^{d}\setminus D_{t})} \vert z-y_{z} \vert \nonumber
\\[-8pt]
\\[-8pt]\qquad
&\leq&\sup_{z\in\overline{D_{s}}\cap(%
\mathbb{R}
^{d}\setminus D_{t})} \vert\tilde{z}-y_{z} \vert\leq\sup_{\tilde{%
z}\in\partial D_{s}}d ( \tilde{z},\partial D_{t} ) \leq\hat{l}%
( r ) .\nonumber
\end{eqnarray}
Combining \eqref{lhat1+}--\eqref{lhat5}, we conclude that $\tilde
{l}_{2}\leq
\hat{l} ( r ) $. Using \eqref{tildel+}, it is clear that $l (
r ) \leq\hat{l} ( r ) +\epsilon$ and then, as $\epsilon$ is
arbitrary, $l ( r ) \leq\hat{l} ( r ) $. This completes
the proof of the lemma.
\end{pf}

\begin{remark}
\label{remsim} Note also that the prerequisites of Lemma \ref%
{llhatequivalence} ensure that the number of holes in $D_{t}$ stays the
same for all $t\in [ 0,T ] $ and, in particular, that these holes
cannot shrink too much as time changes. Furthermore, Lemma \ref%
{llhatequivalence} and its proof allow us to conclude that
%
\begin{equation}
h ( \overline{D_{s}},\overline{D_{t}} ) =h ( D_{s},D_{t} )
=h ( \partial D_{s},\partial D_{t} ) \label{conseqsimply}
\end{equation}
whenever $s,t\in [ 0,T ]$, $|s-t| \leq r$, $l ( r )
<r_{0} $.
\end{remark}

Concerning $\Gamma$, we let $\Gamma=\Gamma_{t}(z)=\Gamma(t,z)$ be a
function defined on $%
\mathbb{R}
^{d+1}$ such that $\Gamma_{t} ( z ) $ is a closed convex cone of
vectors in $%
\mathbb{R}
^{d}$ for every $z\in\partial D_{t}$, $t\in\lbrack0,T]$ and we assume
that $\Gamma$ satisfies \eqref{cone1} and \eqref{cone2}. To
understand the
condition in~\eqref{cone2}, that is, the assumption that the graph
$G^{\Gamma
}$ is closed, we observe that one motivation for using a cone of reflection,
rather than a single-valued direction of reflection, is to be able to deal
with discontinuities in the direction of reflection. Such discontinuities
arise, for instance, in the normal direction for a convex polygon. At a
point of discontinuity of the direction of reflection one can use the cone
generated by all the limit vectors (if they exist) of the direction of
reflection. For a cone of reflection, the assumption that the graph $%
G^{\Gamma}$ is closed provides a form of continuity of the cone. In fact,
for a cone of the form $\Gamma_{t}(z)=\{\lambda\gamma_{t}(z),\lambda
\geq
0\}$, for some $%
\mathbb{R}
^{d}$-valued function $\gamma_{t}(z)$, the assumption that $G^{\Gamma
}$ is
closed is equivalent to the assumption that the function $\gamma
_{t}(z)$ is
continuous as a function of $(t,z)$.

The cone $N_{t} ( z ) $ of inward normal vectors at $z\in\partial
D_{t}$, $t\in\lbrack0,T]$, is defined as being equal to the set consisting
of the union of the set $ \{ 0 \}$ and the set
%
\begin{equation}\hspace*{30pt}
\{ v\in
\mathbb{R}
^{d}\dvtx v\neq0,\exists\rho>0\mbox{ such that }B_{\rho} ( z-\rho
v/ \vert v \vert ) \subset([0,T]\times
\mathbb{R}
^{d})\setminus D \} . \label{conenorm}
\end{equation}
Note that this definition does not rule out the possibility of several unit
inward normal vectors at the same boundary point. Given $N_{t} (
z ) $, we let $N_{t}^{1}(z):=N_{t}(z)\cap S_{1}(0)$, so that $%
N_{t}^{1}(z)$ contains the set of vectors in $N_{t}(z)$ with unit length.
Moreover, based on $N_{t} ( z ) $, we introduce the set
%
\begin{equation}
G^{N}= \{ ( t,z,n ) \dvtx n\in N_{t} ( z ) , z\in
\partial D_{t},t\in [ 0,T ] \} . \label{Gnormal}
\end{equation}
The spatial domain $D_{t}$ is said to verify the uniform exterior sphere
condition if there exists a radius $r_{0}>0$ such that \eqref{extsphere-}
holds. It is easy to see that \eqref{extsphere-} is equivalent to the
statement that
%
\begin{equation}
\langle n,y-z \rangle+\frac{1}{2r_{0}} \vert y-z \vert
^{2}\geq0  \label{extsphere}
\end{equation}
for all $y\in\overline{D_{t}}$, $n\in N_{t}^{1} ( z ) $ and $z\in
\partial D_{t}$. Moreover, as deduced from Remark 2.1 in \cite{Costantini1992},
the uniform exterior sphere condition in time asserts that
$N_{t} ( z ) $ is a closed convex cone for all $z\in\partial D_{t}$%
, $t\in\lbrack0,T]$ and that $G^{N}$ is closed.

For $z\in\partial D_{s}$, $s\in [ 0,T ] $, and $\rho,\eta>0$,
recall the definition of the quantity $a_{s,z} ( \rho,\eta )$
introduced in \eqref{adef},
\[
a_{s,z} ( \rho,\eta ) =\max_{u\in S_{1} ( 0 )
}\min_{s\leq t\leq s+\eta}\min_{y\in\partial D_{t}\cap\overline
{B_{\rho
} ( z ) }}\min_{\gamma\in\Gamma_{t}^{1} ( y )
} \langle\gamma,u \rangle.
\]
The vector $u$ that maximizes the minimum of $ \langle\gamma
,u \rangle$ over all vectors $\gamma\in\Gamma_{t}^{1} (
y ) $ in a time--space neighborhood of a point $ ( s,z ) $, $%
z\in\partial D_{s}$, $s\in [ 0,T ] $, can be regarded as the best
approximation of the $\Gamma_{t}^{1} ( y ) $-vectors in that
neighborhood. With this interpretation $a_{s,z} ( \rho,\eta ) $
represents the cosine of the largest angle between the best approximation
and a $\Gamma_{t}^{1} ( y ) $-vector in the neighborhood. Hence,
in a sense, $a_{s,z} ( \rho,\eta ) $ quantifies the variation of $%
\Gamma$ in a space--time neighborhood of $(s,z)$. For $z\in\partial
D_{s}$%
, $s\in [ 0,T ] $ and $\rho,\eta>0$, recall the definition of
the quantity $c_{s,z} ( \rho,\eta )$ introduced in \eqref{cdef},
\[
c_{s,z} ( \rho,\eta ) =\max_{s\leq t\leq s+\eta}\max_{y\in
\partial D_{t}\cap\overline{B_{\rho} ( z ) }}\max_{\hat{z}\in
\overline{D_{t}}\cap\overline{B_{\rho} ( z ) }, \hat{z}\neq
y}\max_{\gamma\in\Gamma_{t}^{1} ( y ) } \biggl( \frac{ \langle
\gamma,y-\hat{z} \rangle}{ \vert y-\hat{z} \vert}\vee
0 \biggr).
\]
This quantity is close to one if the vectors $\gamma\in\Gamma
_{t}^{1} ( y ) $, in a time--space neighborhood, deviate much from
the normal vectors and/or the domain is very concave. Hence, in a
sense, $%
c_{s,z} ( \rho,\eta ) $ quantifies the skewness of $\Gamma$ and
the concavity of $D$. Note that (\ref{cdef}) implies%
%
\begin{equation}
\langle\gamma,\hat{z}-y \rangle+c_{s,z} ( \rho,\eta
) \vert y-\hat{z} \vert\geq0  \label{cdef2}
\end{equation}
for all $y\in\partial D_{t}\cap\overline{B_{\rho} ( z ) }$, $%
\hat{z}\in\overline{D_{t}}\cap\overline{B_{\rho} ( z ) }$, $ %
\hat{z}\neq y$ and $\gamma\in\Gamma_{t}^{1} ( y ) $ with $z\in
\partial D_{t}$, $t\in [ s,s+\eta ] \subset [ 0,T ] $.
This condition exhibits some similarity with the uniform exterior sphere
property (\ref{extsphere}). Finally, recall the definition of the
quantity $%
e_{s,z} ( \rho,\eta )$ introduced in \eqref{edef},
\[
e_{s,z} ( \rho,\eta ) =\frac{c_{s,z} ( \rho,\eta ) }{%
( a_{s,z} ( \rho,\eta ) ) ^{2}\vee a_{s,z} ( \rho
,\eta ) /2}.
\]
Furthermore, as stated in the \hyperref[intro]{Introduction}, in the subsequent section we
prove estimates related to the Skorohod problem in time-dependent domains
satisfying \eqref{timedep+} and the uniform exterior sphere condition in
time, with radius $r_{0}$. Moreover, to derive these estimates, we also
assume that there exist $0<\rho_{0}<r_{0}$ and $\eta_{0}>0$, such
that the
assumptions in \eqref{crita} and \eqref{crite} hold, that is,
\begin{eqnarray*}
\inf_{s\in [ 0,T ] }\inf_{z\in\partial D_{s}}a_{s,z} ( \rho
_{0},\eta_{0} ) &=&a>0, \\
\sup_{s\in [ 0,T ] }\sup_{z\in\partial D_{s}}e_{s,z} ( \rho
_{0},\eta_{0} ) &=&e<1.
\end{eqnarray*}

\begin{remark}
The function $a_{s,z} ( \rho,\eta ) $ is a straightforward
generalization of the function
%
\begin{equation}
\alpha_{z} ( \rho ) =\max_{u\in S_{1} ( 0 ) }\min_{y\in
\partial\Omega\cap\overline{B_{\rho} ( z ) }}\min_{n\in
N^{1} ( y ) } \langle n,u \rangle, \label{Sai1}
\end{equation}
introduced by Tanaka \cite{Tanaka1979} in his treatment of the Skorohod
problem. Here $\Omega\subset
\mathbb{R}
^{d}$ is a bounded spatial domain and $N^{1}(y)$ is the set of unit inward
normals at $y\in\partial\Omega$. In \cite{Costantini1992,LionsSznitman1984,Saisho1987} and
\cite{Tanaka1979} the
condition%
%
\begin{equation}
\lim_{\rho\rightarrow0}\inf_{z\in\partial\Omega}\alpha_{z} (
\rho
) =\alpha>0 \label{Sai2}
\end{equation}
is used to rule out the case of tangential normal directions (see \cite
{Costantini1992} for equivalent characterizations of domains satisfying this
criterion).
\end{remark}

\begin{remark}
The functions $c_{s,z} ( \rho,\eta ) $ and $e_{s,z} ( \rho
,\eta ) $, introduced in \eqref{cdef} and \eqref{edef}, are
straightforward generalizations of the functions $\widetilde{c}$ and $%
\widetilde{e}$, respectively, introduced in \cite{Costantini1992}. Moreover,
the related functions $c$ and $e$, also introduced in \cite{Costantini1992},
are useful only in the context of convex domains. Hence, as we here consider
general (possibly nonconvex) domains, only generalizations of the functions~$\widetilde{c}$ and $\widetilde{e}$ are useful. For notational simplicity,
we have removed the tilde in our definition of the generalized versions
of $%
\widetilde{c}$ and $\widetilde{e}$.
\end{remark}

Given $T>0$, let $D\subset
\mathbb{R}
^{d+1}$ be a time-dependent domain satisfying \eqref{timedep+} and a uniform
exterior sphere condition in time with radius $r_{0}$ in the sense of %
\eqref{extsphere-}. Given a point $(t,z)\in([0,T]\times
\mathbb{R}
^{d})\setminus\overline{D}$ in a neighborhood of $D$, in this article
we heavily use the projection of $(t,z)$ onto $\partial D$ along the
vectors in
the cone $\Gamma$. While such projections can be defined in several ways,
in this article we here only consider projections in space along
vectors $%
\gamma\in\Gamma_{t} ( y ) $, $y\in\partial D_{t}$, onto $%
\partial D_{t}$. With this restriction, the analysis of Section~\ref{sectcomp} in
\cite{Costantini1992} can be used to derive sufficient conditions for the
existence of a projection of a point $z\in
\mathbb{R}
^{d}\setminus\overline{D_{t}}$, onto $\partial D_{t}$, along $\Gamma_{t}$.
In other words, we can determine whether or not there exist, for a given
point $z\in
\mathbb{R}
^{d}\setminus\overline{D_{t}}$, a point $y\in\partial D_{t}$ and a
vector $\gamma\in\Gamma_{t} ( y ) $ such that $y-z\parallel
\gamma$. In particular, it can be understood when, for $0<\delta
_{0}<r_{0}$%
, $h_{0}>1$ and $\Gamma=\Gamma_{t}(z)=\Gamma(t,z)$ given,
$([0,T]\times
\mathbb{R}
^{d})\setminus\overline{D}$ has the $(\delta_{0},h_{0})$-property of good
projections along $\Gamma$ in the sense defined in the introduction;
see %
\eqref{deltanolldef} and \eqref{rhonolldef}. We refer to the
\hyperref[appx]{Appendix}, for more on this, as well as for a discussion of examples
of time-dependent domains satisfying the restrictions imposed in
Theorems \ref%
{Theorem 1},   \ref{Theorem 2} and   \ref{Theorem 3}.

\subsection{C\`{a}dl\`{a}g functions and the Skorohod topology}\label
{skorohodtop}

Let $T>0$ and let $x\in\mathcal{D}( [ 0,T ] ,%
\mathbb{R}
^{d})$. Given a bounded set $I\subset [ 0,T ] $, we let%
%
\begin{equation}
\widehat{w} ( x,I ) =\sup_{u,r\in I} \vert
x_{u}-x_{r} \vert.
\end{equation}
Then, using Lemma 1 on page 122 in \cite{Billingsley1999}, we see that there
exists, for $\epsilon>0$ given, a sequence of points $t_{0},\ldots
,t_{\nu}$,
such that%
%
\begin{equation}
0=t_{0}<t_{1}<\cdots<t_{\nu}=T,\qquad \widehat{w} ( x, [
t_{i-1},t_{i} ) ) <\epsilon.
\end{equation}
In particular, there can only be finitely many points $t\in [ 0,T ]
$ at which the jump $ \vert x_{t}-x_{t^{-}} \vert$ exceeds a given
positive number. To proceed, in the following we use the notation and
exposition of Chapter 3 in \cite{EthierKurtz1986}. We let $q (
x,y ) = \vert x-y \vert\wedge1$ whenever $x,y\in
\mathbb{R}
^{d}$ and we let $d_{\mathcal{D}}([0,T],x,y)$ be the metric on the
space $%
\mathcal{D}( [ 0,T ] ,%
\mathbb{R}
^{d})$ introduced, for the interval $[0,T]$, as in display (5.2) in
\cite{EthierKurtz1986}. Then, by Theorem 5.6 in \cite{EthierKurtz1986}, we see
that $(\mathcal{D}( [ 0,T ] ,%
\mathbb{R}
^{d}),d_{\mathcal{D}}([0,T],\cdot,\cdot))$ is a complete metric
space and
the topology on $\mathcal{D}( [ 0,T ] ,%
\mathbb{R}
^{d})$, induced by the metric $d_{\mathcal{D}}([0,T],\cdot,\cdot)$, is
known as the Skorohod topology on $\mathcal{D}( [ 0,T ] ,%
\mathbb{R}
^{d})$. Recall that if $x,y\in\mathcal{D}( [ 0,T ] ,%
\mathbb{R}
^{d})$, then $d_{\mathcal{D}}([0,T],x,y)=0$ implies that $x_{t}=y_{t}$ for
every $t$. Furthermore, if $\{x^{n}\}$ is a sequence in $\mathcal{D}( [
0,T ] ,%
\mathbb{R}
^{d})$ and $x\in\mathcal{D}( [ 0,T ] ,%
\mathbb{R}
^{d})$, then the statement that $d_{\mathcal
{D}}([0,T],x^{n},x)\rightarrow0$
as $n\rightarrow\infty$ is equivalent to the statement that there
exists $%
\{\lambda_{n}\}\subset\Lambda$ (see \cite{EthierKurtz1986} for the
definition of the space $\Lambda$) such that (5.6) in \cite{EthierKurtz1986}
holds and such that
%
\begin{equation}
\lim_{n\rightarrow\infty}\sup_{0\leq t\leq T}\bigl|x_{t}^{n}-x_{\lambda
_{n}(t)}\bigr|=0.
\end{equation}
For a proof of this result we refer to Proposition 5.3 in \cite{EthierKurtz1986}. Furthermore, to understand the relatively compact
sets in
$\mathcal{D}( [ 0,T ] ,%
\mathbb{R}
^{d})$, we introduce and use a modulus of continuity. In particular,
for $%
x\in\mathcal{D}( [ 0,T ] ,%
\mathbb{R}
^{d})$ and $\delta>0$ we define the quantity
%
\begin{equation}
w^{\prime} ( x,\delta,T ) =\inf_{ \{ t_{i} \}
}\max_{i}\sup_{u,r\in [ t_{i},t_{i+1} ) } \vert
x_{u}-x_{r} \vert, \label{b11}
\end{equation}
where the infimum is taken with respect to all partitions of the form $%
0=t_{0}<t_{1}<\cdots<t_{n-1}<T\leq t_{n}$, with $\min_{i} \vert
t_{i}-t_{i-1} \vert>\delta$. Furthermore, given $\mathcal{W}\subset
\mathcal{D} ( [ 0,T ] ,%
\mathbb{R}
^{d} ) $, we let
%
\begin{equation}
\mu ( \mathcal{W},\delta,T ) =\sup_{w\in\mathcal{W}}w^{\prime
} ( w,\delta,T ) . \label{d11}
\end{equation}
Using this notation, we first quote Theorem 6.3 in \cite{EthierKurtz1986}
which states that $\mathcal{W}\subset\mathcal{D} ( [ 0,T ] ,%
\mathbb{R}
^{d} ) $ is relatively compact in the Skorohod topology if and only if
for every rational $t\in\lbrack0,T]$ there exists a relatively
compact set
$A_{t}\subset
\mathbb{R}
^{d}$ such that $w_{t}\in A_{t}$ for all $w\in\mathcal{W}$ and such that
%
\begin{equation}
\lim_{\delta\rightarrow0}\mu ( \mathcal{W},\delta,T ) =0.
\label{e11}
\end{equation}
Finally, we also note the following. Given $\delta^{\prime}>\delta$,
let%
%
\begin{equation}
\tilde{w}^{\prime} ( x,\delta,\delta^{\prime},T )
=\inf_{ \{ t_{i} \} }\max_{i}\sup_{u,r\in [ t_{i},t_{i+1}%
) } \vert x_{u}-x_{r} \vert, \label{c11}
\end{equation}
where the infimum is taken with respect to all partitions as above but with
the additional restriction that $\max_{i} \vert t_{i}-t_{i-1} \vert
<\delta^{\prime}$. Furthermore, given $\mathcal{W}\subset\mathcal
{D}%
( [ 0,T ] ,%
\mathbb{R}
^{d} ) $, we let
%
\begin{equation}
\tilde{\mu} ( \mathcal{W},\delta,\delta^{\prime},T ) =\sup_{w\in
\mathcal{W}}\tilde{w}^{\prime} ( w,\delta,\delta^{\prime},T ) .
\label{d11e}
\end{equation}
Then
%
\begin{equation}\qquad
\tilde{w}^{\prime} ( x,\delta,\delta^{\prime},T ) ={w}^{\prime
} ( x,\delta,T ) \quad \mbox{and}\quad \tilde{\mu} ( \mathcal{W},\delta
,\delta^{\prime},T ) ={\mu} ( \mathcal{W},\delta,T ) .
\label{d11+}
\end{equation}

\section{Estimates for solutions and approximations to Skorohod problems}
\label{sectcomp}

In this section we first prove certain estimates for solutions to the
Skorohod problem for $(D,\Gamma,w)$, assuming that $D$ satisfies the
assumptions stated in Theorem \ref{Theorem 1} and that $w\in\mathcal
{D}%
( [ 0,T ] ,%
\mathbb{R}
^{d} ) $ with $w_{0}\in\overline{D_{0}}$. In particular, we prove that
the modulus of continuity of c\`{a}dl\`{a}g solutions to the Skorohod
problem for $(D,\Gamma,w)$, with bounded jumps, can be estimated from above
by the modulus of continuity of $w$ and the modulus of continuity $l$. This
result is derived in two steps. In the first step we prove (see Lemma
\ref%
{compactlemma} below) a local compactness result which is valid in a spatial
neighborhood of a given boundary point and on a constructed time interval.
In the second step we then prove that corresponding global estimates (see
Theorem~\ref{compacttheorem} below) can be derived based on the local
compactness result. In particular, Theorem \ref{compacttheorem} is the main
result we establish in this context. In Section \ref{SPest} we derive these
estimates for solutions to the Skorohod problem and in Section \ref{SPappest}
we establish the corresponding results for approximations to the Skorohod
problem. In particular, given $w\in\mathcal{D} ( [ 0,T ] ,%
\mathbb{R}
^{d} ) $ with $w_{0}\in\overline{D_{0}}$ and a partition $\{\tau
_{k}\}_{k=0}^{N}$, which\vspace*{1pt} we denote by $\Delta$, of the interval
$[0,T]$, we
define $w^{\Delta}$, $D^{\Delta}$, $\Gamma^{\Delta}$, $x^{\Delta}$
and $%
\lambda^{\Delta}$ as in~\eqref{approxw}, \eqref{approxw+}, %
\eqref{approxw++} and \eqref{approxw+++}. Then, by construction, the
pair $%
(x^{\Delta},\lambda^{\Delta})$ is a solution to the Skorohod problem
for $%
( D^{\Delta},\Gamma^{\Delta},w^{\Delta} ) $. In Lemma~\ref{compactlemmaapprox} and Theorem~\ref{compacttheoremapprox} we prove
estimates for solutions to the Skorohod problem for $ ( D^{\Delta
},\Gamma^{\Delta},w^{\Delta} ) $, which are similar to the ones
established in Lemma \ref{compactlemma} and Theorem~\ref
{compacttheorem} for
the Skorohod problem for $ ( D,\Gamma,w ) $. We note that the
reason for this twofold approach is that since we are considering
time-dependent domains, the condition in \eqref{limitzero} will in general
not hold for $D^{\Delta}$.

\subsection{Estimates for solutions to Skorohod problems}\label{SPest}

Given $a>0$ and $e\in(0,1)$, we define the positive functions $%
K_{1},K_{2},K_{3}$ and $K_{4}$ as follows:
%
\begin{eqnarray}\label{functions}
K_{1} ( a,e ) &=&\frac{a+2a^{2}e+2+ae}{a ( 1-e ) },\qquad
K_{2} ( a,e ) =\frac{2a^{2}e+2+ae}{a ( 1-e ) }, \nonumber
\\[-8pt]
\\[-8pt]
K_{3} ( a,e ) &=&\frac{1+K_{1} ( a,e ) }{a},\qquad
K_{4} ( a,e ) =\frac{1+K_{2} ( a,e ) }{a}.\nonumber
\end{eqnarray}
In this section we first prove the following two general results for
solutions to the Skorohod problem.

\begin{lemma}
\label{compactlemma} Let $T>0$ and let $D\subset
\mathbb{R}
^{d+1}$ be a time-dependent domain satisfying \eqref{timedep+}, %
\eqref{limitzero} and a uniform exterior sphere condition in time with
radius $r_{0}$ in the sense of \eqref{extsphere-}. Let $\Gamma=\Gamma
_{t} ( z ) $ be a closed convex cone of vectors in $%
\mathbb{R}
^{d}$ for every $z\in\partial D_{t}$, $t\in [ 0,T ] $. Assume
that \eqref{crita} and~\eqref{crite} hold for some $0<\rho
_{0}<r_{0}$, $%
\eta_{0}>0$, $a$ and $e$. Finally, assume that $([0,T]\times
\mathbb{R}
^{d})\setminus\overline{D}$ has the $(\delta_{0},h_{0})$-property of good
projections along $\Gamma$, for some $0<\delta_{0}<r_{0}$, $h_{0}>1$ as
defined in \eqref{deltanolldef} and~\eqref{rhonolldef}. Let $w\in
\mathcal{D}%
( [ 0,T ] ,%
\mathbb{R}
^{d} ) $ with $w_{0}\in\overline{D_{0}}$ and let $ ( x,\lambda
) $ be a solution to the Skorohod problem for $ ( D,\Gamma
,w ) $. Consider a fixed but arbitrary $s\in [ 0,T ] $, such
that $x_{s}\in\partial D_{s}$, and note that it follows from \eqref{crita}
and \eqref{crite} that there exist $0<\rho<r_{0}$ and $\eta>0$ such
that%
%
\begin{equation}
a_{s,x_{s}} ( \rho,\eta ) >0,\qquad e_{s,x_{s}} ( \rho,\eta
) <1.
\end{equation}
Then, for $0\leq s\leq t_{1}\leq t_{2}<\tau_{\rho,\eta}$,
%
\begin{eqnarray}\label{lemma1r1}
\Vert x \Vert_{t_{1},t_{2}} &\leq&K_{1} ( a,e )
\Vert w \Vert_{t_{1},t_{2}}+K_{2} ( a,e ) l (
t_{2}-t_{1} ) , \\\label{lemma1r2}
\vert\lambda \vert_{t_{2}}- \vert\lambda \vert
_{t_{1}} &\leq&K_{3} ( a,e ) \Vert w \Vert
_{t_{1},t_{2}}+K_{4} ( a,e ) l ( t_{2}-t_{1} ) ,
\end{eqnarray}
where $a=a_{s,x_{s}}(\rho,\eta)$, $e=e_{s,x_{s}}(\rho,\eta)$. Here
$\tau
_{\rho,\eta}$ is defined as follows. If there exists some $t$ such
that $%
s\leq t<(s+\eta)\wedge T$ and $ \vert x_{t}-x_{s} \vert+l (
t-s ) \geq\rho$, then
%
\begin{equation}
\tau_{\rho,\eta}=\inf\{t\dvtx  s\leq t<(s+\eta)\wedge T, \vert
x_{t}-x_{s} \vert+l ( t-s ) \geq\rho\}, \label{stoptid1}
\end{equation}
whereas if $ \vert x_{t}-x_{s} \vert+l ( t-s ) <\rho$
for all $s\leq t<(s+\eta)\wedge T$, then
%
\begin{equation}
\tau_{\rho,\eta}=(s+\eta)\wedge T. \label{stoptid2}
\end{equation}
\end{lemma}

\begin{theorem}
\label{compacttheorem} Let $T>0$, $D\subset
\mathbb{R}
^{d+1}$, $r_{0}$, $\Gamma=\Gamma_{t} ( z ) $, $0<\rho_{0}<r_{0}$%
, $\eta_{0}>0$, $a$, $e$, $\delta_{0}$ and $h_{0}$ be as in the statement
of Lemma \ref{compactlemma}. Let $w\in\mathcal{D} ( [ 0,T ] ,%
\mathbb{R}
^{d} ) $ with $w_{0}\in\overline{D_{0}}$ and let $ ( x,\lambda
) $ be a solution to the Skorohod problem for $ ( D,\Gamma
,w ) $. Moreover, assume in addition that $x\in\mathcal{D}^{\rho
_{0}} ( [ 0,T ] ,%
\mathbb{R}
^{d} ) $. Then there exist positive constants $L_{1} ( w,T ) $%
, $L_{2} ( w,T ) $, $L_{3} ( w,T ) $ and $L_{4} (
w,T ) $ such that
%
\begin{eqnarray}
\Vert x \Vert_{t_{1},t_{2}} &\leq&L_{1} ( w,T )
\Vert w \Vert_{t_{1},t_{2}}+L_{2} ( w,T ) l (
t_{2}-t_{1} ) , \\
|\lambda|_{t_{2}}-|\lambda|_{t_{1}} &\leq&L_{3} ( w,T )
\Vert w \Vert_{t_{1},t_{2}}+L_{4} ( w,T ) l (
t_{2}-t_{1} ) ,
\end{eqnarray}
whenever $0\leq t_{1}\leq t_{2}\leq T$. Furthermore, if $\mathcal
{W}\subset
\mathcal{D} ( [ 0,T ] ,%
\mathbb{R}
^{d} ) $ is relatively compact in the Skorohod topology and satisfies $%
w_{0}\in\overline{D_{0}}$, whenever $w\in\mathcal{W}$, then there
exist positive constants $L_{1}^{T}$, $L_{2}^{T}$, $L_{3}^{T}$ and $%
L_{4}^{T} $, such that%
%
\begin{equation}
\sup_{w\in\mathcal{W}}L_{i} ( w,T ) \leq L_{i}^{T}<\infty\qquad
\mbox{for }i=1,2,3,4.
\end{equation}
\end{theorem}

\begin{remark}
Versions of Lemma \ref{compactlemma} and Theorem \ref
{compacttheorem}, valid
only in the setting of time-independent domains, are proved in Lemma~2.1,
Theorems~2.2 and~2.4 in \cite{Costantini1992}. Our contribution is that we
are able to establish similar results when $D\subset
\mathbb{R}
^{d+1}$ is a time-dependent domain. Furthermore, concerning related results
in the setting of time-dependent domains, we note that if $D$ is an
$\mathcal{%
H}_{2}$-domain and if $\Gamma_{t}(z)=\{\lambda\gamma_{t}(z),\lambda
\geq
0\}$, for some $S_{1} ( 0 ) $-valued continuous function $\gamma
_{t} ( z ) $ such that%
%
\begin{equation}
\inf_{z\in\partial D_{t}, t\in\lbrack0,T]}\langle\gamma
_{t}(z),n_{t}(z)\rangle>\frac{\sqrt{3}}{2},
\end{equation}
then a version of Theorem \ref{compacttheorem} is proved in Theorem
C.3 in
\cite{CostantiniGobetKaroui2006}. Note also that if $D$ is an
$\mathcal{H}%
_{2}$-domain, then there exists a unique unit inward normal,
$n_{t}(z)$, at $%
z\in\partial D_{t}$, $t\in\lbrack0,T]$.
\end{remark}

\begin{remark}
\label{conenotres} Unlike in the statements of Theorems \ref{Theorem
1}, \ref%
{Theorem 2} and \ref{Theorem 3}, we need not assume that $\Gamma$
satisfies %
\eqref{cone1}, \eqref{cone2} and \eqref{Gammanpre} in the
prerequisites of
Lemma~\ref{compactlemma} and Theorem \ref{compacttheorem}. This
remark also
applies to Lemma \ref{compactlemmaapprox} and Theorem \ref%
{compacttheoremapprox} stated below.
\end{remark}

\begin{pf*}{Proof of Lemma \ref{compactlemma}}
To simplify the
notation, we in the following let $a=a_{s,x_{s}} ( \rho,\eta ) $, $%
c=c_{s,x_{s}} ( \rho ) $ and $e=e_{s,x_{s}} ( \rho,\eta
) $. Moreover, we let $u$ be a unit vector such that
%
\begin{equation}
\langle\gamma_{r},u \rangle\geq a
\end{equation}
for all $\gamma_{r}\in\Gamma_{r}^{1} ( y ) $, $y\in\partial
D_{r}\cap\overline{B_{\rho} ( x_{s} ) }$ and $r\in [ s,\tau
_{\rho,\eta} ] \subset [ s,(s+\eta)\wedge T ] $. The
existence of such a vector follows from the definition of $a_{s,x_{s}} (
\rho,\eta ) $. Using properties \eqref{SP1}--\eqref{SP2} in
Definition %
\ref{skorohodprob}, we see that
%
\begin{equation}
\langle x_{t_{2}}-x_{t_{1}},u \rangle= \langle
w_{t_{2}}-w_{t_{1}},u \rangle+\int_{t_{1}^{+}}^{t_{2}^{+}}
\underbrace{ \langle\gamma_{r},u \rangle}_{
\geq a}\,d \vert
\lambda \vert_{r}  \label{hahaha}
\end{equation}
for any $0\leq s\leq t_{1}\leq t_{2}<\tau_{\rho,\eta}$. Based on %
\eqref{hahaha}, we deduce that
%
\begin{equation}
\vert\lambda \vert_{t_{2}}- \vert\lambda \vert
_{t_{1}}\leq\frac{1}{a} ( \vert w_{t_{2}}-w_{t_{1}} \vert
+ \vert x_{t_{2}}-x_{t_{1}} \vert ) . \label{lemma11}
\end{equation}
Furthermore, again using properties \eqref{SP1}--\eqref{SP2} in
Definition %
\ref{skorohodprob}, we also see that
%
\begin{eqnarray} \label{haha1}\qquad
\vert x_{t_{2}}-x_{t_{1}} \vert^{2}
&=& \biggl(
w_{t_{2}}-w_{t_{1}}+\int_{t_{1}^{+}}^{t_{2}^{+}}\gamma_{r}\,d \vert
\lambda \vert_{r} \biggr) ^{2} \nonumber
\\[-8pt]
\\[-8pt]\qquad
&=&
\vert w_{t_{2}}-w_{t_{1}} \vert^{2}+ \biggl(
\int_{t_{1}^{+}}^{t_{2}^{+}}\gamma_{r}\,d \vert\lambda \vert
_{r} \biggr) ^{2}+2\int_{t_{1}}^{t_{2}^{+}} \langle
w_{t_{2}}-w_{t_{1}},\gamma_{r} \rangle\, d \vert\lambda \vert
_{r}, \nonumber
\end{eqnarray}
whenever $0\leq s\leq t_{1}\leq t_{2}<\tau_{\rho,\eta}$. Note that the
integrand in the last term in this display can be rewritten as
%
\begin{eqnarray} \label{haha2}\hspace*{34pt}
\langle w_{t_{2}}-w_{t_{1}},\gamma_{r} \rangle&=& \langle
w_{t_{2}}-w_{r},\gamma_{r} \rangle+ \langle
w_{r}-w_{t_{1}},\gamma_{r} \rangle\nonumber
\\[-8pt]
\\[-8pt] \hspace*{34pt}
&=& \langle w_{t_{2}}-w_{r},\gamma_{r} \rangle+ \langle
x_{r}-x_{t_{1}},\gamma_{r} \rangle- \biggl\langle \biggl(
\int_{t_{1}^{+}}^{r^{+}}\gamma_{u}\,d \vert\lambda \vert
_{u} \biggr) ,\gamma_{r} \biggr\rangle. \nonumber\hspace*{-24pt}
\end{eqnarray}
In particular, combining \eqref{haha1} and \eqref{haha2}, we see that
%
\begin{eqnarray} \label{lemma15}
\vert x_{t_{2}}-x_{t_{1}} \vert^{2} &=& \vert
w_{t_{2}}-w_{t_{1}} \vert
^{2}+2\int_{t_{1}^{+}}^{t_{2}^{+}} \langle w_{t_{2}}-w_{r},\gamma
_{r} \rangle\, d \vert\lambda \vert_{r}
\nonumber\\
&&{}+2\int_{t_{1}^{+}}^{t_{2}^{+}} \langle x_{r}-x_{t_{1}},\gamma
_{r} \rangle\, d \vert\lambda \vert_{r}+ \biggl(
\int_{t_{1}^{+}}^{t_{2}^{+}}\gamma_{r}\,d \vert\lambda \vert
_{r} \biggr) ^{2} \\
&&{}-2\int_{t_{1}^{+}}^{t_{2}^{+}} \biggl\langle \biggl(
\int_{t_{1}{}^{+}}^{r^{+}}\gamma_{u} ( x_{u} )\, d \vert\lambda
\vert_{u} \biggr) ,\gamma_{r} \biggr\rangle\, d \vert\lambda
\vert_{r}.\nonumber
\end{eqnarray}
We now intend to derive bounds from above for all integrals in (\ref
{lemma15}%
). To do this, we first note that the first integral on the right-hand side
of (\ref{lemma15}) is bounded from above by%
%
\begin{equation}
2\int_{t_{1}^{+}}^{t_{2}^{+}} \langle w_{t_{2}}-w_{r},\gamma
_{r} \rangle\, d \vert\lambda \vert_{r}\leq
2\int_{t_{1}^{+}}^{t_{2}^{+}} \vert w_{t_{2}}-w_{r} \vert\,
d \vert\lambda \vert_{r}. \label{lemma19}
\end{equation}
To find an upper bound of the second integral in (\ref{lemma15}), we must
take into account that, due to the fact that our domain is
time-dependent, $%
x_{t_{1}}$ might not belong to $\overline{D_{r}}$. Recall that we are
assuming that $N_{r}(y)\neq\emptyset$ for all $y\in\partial D_{r}$,
$r\in
\lbrack0,T]$, and, given $r\in\lbrack0,T]$, $y\in
\mathbb{R}
^{d}\setminus\overline{D_{r}}$, in the following we denote a
projection of $%
y$ onto $\partial D_{r}$ along $N_{r}$ by $\pi_{\partial
D_{r}}^{N_{r}} ( y ) $. Furthermore, whenever $y\in\overline{D_{r}}
$ we let $\pi_{\partial D_{r}}^{N_{r}} ( y ) =y$.\vspace*{1pt} Using this
notation, and the definition of $\tau_{\rho,\eta}$, we see that
%
\begin{equation}
\vert\pi_{\partial D_{r}}^{N_{r}} ( x_{t_{1}} )
-x_{s} \vert\leq \vert x_{t_{1}}-x_{s} \vert+l (
r-t_{1} ) \leq\rho. \label{haha3}
\end{equation}
Equation \eqref{haha3} implies that $\pi_{\partial D_{r}}^{N_{r}} (
x_{t_{1}} ) \in\overline{B_{\rho} ( x_{s} ) }\cap
\partial D_{r}\subset\overline{B_{\rho} ( x_{s} ) }\cap\overline{%
D_{r}}$. Next, writing%
%
\begin{equation}
\langle x_{r}-x_{t_{1}},\gamma_{r} \rangle= \langle
x_{r}-\pi_{\partial D_{r}}^{N_{r}} ( x_{t_{1}} ) ,\gamma
_{r} \rangle+ \langle\pi_{\partial D_{r}}^{N_{r}} (
x_{t_{1}} ) -x_{t_{1}},\gamma_{r} \rangle, \label{haha4}
\end{equation}
and using the fact that $x_{r}\in\overline{B_{\rho} ( x_{s} ) }%
\cap\partial D_{r}$ a.e. when $d \vert\lambda \vert_{r}\neq0$,
together with a version of (\ref{cdef2}), we deduce that
%
\begin{equation}\quad
\langle x_{r}-\pi_{\partial D_{r}}^{N_{r}} ( x_{t_{1}} )
,\gamma_{r} \rangle\leq c \vert x_{r}-\pi_{\partial
D_{r}}^{N_{r}} ( x_{t_{1}} ) \vert\leq c \vert
x_{r}-x_{t_{1}} \vert+cl ( r-t_{1} ) . \label{haha5}
\end{equation}
Furthermore,
%
\begin{equation}
\langle\pi_{\partial D_{r}}^{N_{r}} ( x_{t_{1}} )
-x_{t_{1}},\gamma_{r} \rangle\leq \vert\pi_{\partial
D_{r}}^{N_{r}} ( x_{t_{1}} ) -x_{t_{1}} \vert\leq l (
r-t_{1} ) . \label{haha7}
\end{equation}
Using the estimates derived above, we conclude that the second integral
in (%
\ref{lemma15}) has the upper bound
%
\begin{eqnarray}\label{lemma20}
2\int_{t_{1}^{+}}^{t_{2}^{+}} \langle x_{r}-x_{t_{1}},\gamma
_{r} \rangle \,d \vert\lambda \vert_{r} &\leq&
2c\int_{t_{1}^{+}}^{t_{2}^{+}} \vert x_{r}-x_{t_{1}} \vert
\,d \vert\lambda \vert_{r} \nonumber
\\[-8pt]
\\[-8pt]
&&{}+2 ( c+1 ) l ( t_{2}-t_{1} ) ( \vert\lambda
\vert_{t_{2}}- \vert\lambda \vert_{t_{1}} ) .\nonumber
\end{eqnarray}
Next, we use Lemma 2.1(ii) in \cite{Saisho1987} and rewrite the third
integral in (\ref{lemma15}) as
%
\begin{eqnarray}\label{haha6}
\biggl( \int_{t_{1}^{+}}^{t_{2}^{+}}\gamma_{r}\,d \vert\lambda
\vert_{r} \biggr) ^{2}
&=&2\int_{t_{1}^{+}}^{t_{2}^{+}} \biggl\langle
\biggl( \int_{t_{1}^{+}}^{r^{+}}\gamma_{u} ( x_{u} )\, d \vert
\lambda \vert_{u} \biggr) ,\gamma_{r} \biggr\rangle\, d \vert
\lambda \vert_{r} \nonumber
\\[-8pt]
\\[-8pt]
&&{}-\sum_{t_{1}<r\leq t_{2}} \underbrace{ \vert\gamma
_{r} \vert^{2}}_{=1} ( \vert\lambda \vert_{r}- \vert
\lambda \vert_{r^{-}} ) ^{2}.\nonumber
\end{eqnarray}
Based on the last display, it is clear that and the third and fourth integral
in (\ref{lemma15}) reduce to the term%
%
\begin{equation}
-\sum_{t_{1}<r\leq t_{2}} ( \vert\lambda \vert
_{r}- \vert\lambda \vert_{r^{-}} ) ^{2}. \label{lemma21}
\end{equation}
Putting the relations (\ref{lemma15})--(\ref{lemma21}) together, we
obtain%
%
\begin{eqnarray}\label{lemma12}\qquad
\vert x_{t_{2}}-x_{t_{1}} \vert^{2}
&\leq& \vert
w_{t_{2}}-w_{t_{1}} \vert^{2}+2\int_{t_{1}^{+}}^{t_{2}^{+}} \vert
w_{t_{2}}-w_{r} \vert\, d \vert\lambda \vert
_{r}+2c\int_{t_{1}^{+}}^{t_{2}^{+}} \vert x_{r}-x_{t_{1}} \vert
\,d \vert\lambda \vert_{r} \nonumber
\\[-8pt]
\\[-8pt]\qquad
&&{}-\sum_{t_{1}<r\leq t_{2}} ( \vert\lambda \vert
_{r}- \vert\lambda \vert_{r^{-}} ) ^{2}+2 ( c+1 )
l ( t_{2}-t_{1} ) ( \vert\lambda \vert
_{t_{2}}- \vert\lambda \vert_{t_{1}} ) .\nonumber
\end{eqnarray}
If we now combine (\ref{lemma12}) and the properties \eqref
{SP1}--\eqref{SP2}
in Definition \ref{skorohodprob}, we first get%
%
\begin{eqnarray}\label{haha8}
\vert x_{t_{2}}-x_{t_{1}} \vert^{2} &\leq& \vert
w_{t_{2}}-w_{t_{1}} \vert^{2}+2\int_{t_{1}^{+}}^{t_{2}^{+}} \vert
w_{t_{2}}-w_{r} \vert\, d \vert\lambda \vert
_{r}+2c\int_{t_{1}^{+}}^{t_{2}^{+}} \vert w_{r}-w_{t_{1}} \vert
\,d \vert\lambda \vert_{r} \nonumber \\
&&{}+2c\int_{t_{1}^{+}}^{t_{2}^{+}} ( \vert\lambda \vert
_{r}- \vert\lambda \vert_{t_{1}} )\, d \vert\lambda
\vert_{r}-\sum_{t_{1}<r\leq t_{2}} ( \vert\lambda
\vert_{r}- \vert\lambda \vert_{r^{-}} ) ^{2}
\\
&&{}+2 ( c+1 ) l ( t_{2}-t_{1} ) ( \vert\lambda
\vert_{t_{2}}- \vert\lambda \vert_{t_{1}} ) ,\nonumber
\end{eqnarray}
and then, again using Lemma 2.1(ii) in \cite{Saisho1987} as well as the
fact that $0\leq c\leq1$, we conclude that
%
\begin{eqnarray} \label{haha9}\quad
\vert x_{t_{2}}-x_{t_{1}} \vert^{2} &\leq& \vert
w_{t_{2}}-w_{t_{1}} \vert^{2}+2\int_{t_{1}^{+}}^{t_{2}^{+}} \vert
w_{t_{2}}-w_{r} \vert\, d \vert\lambda \vert
_{r}+2c\int_{t_{1}^{+}}^{t_{2}^{+}} \vert w_{r}-w_{t_{1}} \vert
\,d \vert\lambda \vert_{r} \nonumber
\\[-8pt]
\\[-8pt]\quad
&&{}+c ( \vert\lambda \vert_{t_{2}}- \vert\lambda
\vert_{t_{1}} ) ^{2}+2 ( c+1 ) l (
t_{2}-t_{1} ) ( \vert\lambda \vert_{t_{2}}- \vert
\lambda \vert_{t_{1}} ) .\nonumber
\end{eqnarray}
Relation (\ref{lemma11}) and the inequality in the last display yield%
%
\begin{eqnarray}\label{lemma14}
\Vert x \Vert_{t_{1},t_{2}}^{2}
&\leq& \biggl( 1+\frac{2 (
c+1 ) }{a}+\frac{c}{a^{2}} \biggr) \Vert w \Vert
_{t_{1},t_{2}}^{2}+2 \biggl( \frac{c+1}{a}+\frac{c}{a^{2}} \biggr) \Vert
x \Vert_{t_{1},t_{2}} \Vert w \Vert_{t_{1},t_{2}}
\nonumber\\
&&{}+\frac{c}{a^{2}} \Vert x \Vert_{t_{1},t_{2}}^{2}+\frac{2 (
c+1 ) }{a}l ( t_{2}-t_{1} ) \Vert x \Vert
_{t_{1},t_{2}}\\
&&{}+\frac{2 ( c+1 ) }{a}l ( t_{2}-t_{1} )
\Vert w \Vert_{t_{1},t_{2}}. \nonumber
\end{eqnarray}
In addition, combining (\ref{lemma11}) and (\ref{lemma12}), we obtain
%
\begin{eqnarray} \label{lemma13}\quad
\Vert x \Vert_{t_{1},t_{2}}^{2} &\leq& \biggl( 1+\frac{2}{a}%
\biggr) \Vert w \Vert_{t_{1},t_{2}}^{2}+\frac{2 ( c+1 )
}{a} \Vert x \Vert_{t_{1},t_{2}} \Vert w \Vert
_{t_{1},t_{2}}+\frac{2c}{a} \Vert x \Vert_{t_{1},t_{2}}^{2} \nonumber
\\[-8pt]
\\[-8pt]\quad
&&{}+\frac{2 ( c+1 ) }{a}l ( t_{2}-t_{1} ) \Vert
x \Vert_{t_{1},t_{2}}+\frac{2 ( c+1 ) }{a}l (
t_{2}-t_{1} ) \Vert w \Vert_{t_{1},t_{2}}.\nonumber
\end{eqnarray}
The inequalities (\ref{lemma14}) and (\ref{lemma13}) can both be
written on
the form%
%
\begin{equation}\hspace*{30pt}
A \Vert x \Vert_{t_{1},t_{2}}^{2}-B \Vert x \Vert
_{t_{1},t_{2}} \Vert w \Vert_{t_{1},t_{2}}-C \Vert
w \Vert_{t_{1},t_{2}}^{2}-D \Vert x \Vert
_{t_{1},t_{2}}-D \Vert w \Vert_{t_{1},t_{2}}\leq0,
\label{lemma16}
\end{equation}
where the positive constants $A$, $B$ and $C$ are easily shown to satisfy
the condition $A+B=C$ in both cases. We claim that
%
\begin{equation}
\Vert x \Vert_{t_{1},t_{2}}\leq\frac{C}{A} \Vert
w \Vert_{t_{1},t_{2}}+\frac{D}{A}. \label{haha11}
\end{equation}
Indeed, suppose, on the contrary, that
%
\begin{equation}
0\leq \Vert w \Vert_{t_{1},t_{2}}<\frac{A}{C} \Vert
x \Vert_{t_{1},t_{2}}-\frac{D}{C}. \label{haha12}
\end{equation}
Then, by \eqref{haha12},
%
\begin{eqnarray}\label{haha13}
&&A \Vert x \Vert_{t_{1},t_{2}}^{2}-B \Vert x \Vert
_{t_{1},t_{2}} \Vert w \Vert_{t_{1},t_{2}}-C \Vert
w \Vert_{t_{1},t_{2}}^{2}-D \Vert x \Vert
_{t_{1},t_{2}}-D \Vert w \Vert_{t_{1},t_{2}} \nonumber
\\
&&\qquad>A \Vert x \Vert_{t_{1},t_{2}}^{2}+B \Vert x \Vert
_{t_{1},t_{2}} \biggl( -\frac{A}{C} \Vert x \Vert_{t_{1},t_{2}}+%
\frac{D}{C} \biggr) -C \biggl( \frac{A}{C} \Vert x \Vert
_{t_{1},t_{2}}-\frac{D}{C} \biggr) ^{2} \nonumber
\\[-8pt]
\\[-8pt]
&&\qquad\quad {}-D \Vert x \Vert_{t_{1},t_{2}}+D \biggl( -\frac{A}{C} \Vert
x \Vert_{t_{1},t_{2}}+\frac{D}{C} \biggr) \nonumber\\
&&\qquad=A \Vert x \Vert_{t_{1},t_{2}}^{2}
\underbrace{ \biggl( 1-\frac{B}{C}-\frac{A}{C} \biggr) }_{=\frac{C-B-A}{C}=0}{}+ D \Vert
x \Vert_{t_{1},t_{2}}\underbrace{ \biggl(
\frac{B}{C}+\frac{2A}{C}-\frac{A}{C}-1 \biggr) }_{=\frac{A+B-C}{C}=0}
=0.\nonumber
\end{eqnarray}
Obviously \eqref{haha13} contradicts (\ref{lemma16}) and, hence, the
claim in %
\eqref{haha11} is proved. To complete the proof of Lemma \ref{compactlemma},
we first note that (\ref{lemma14}) implies
%
\begin{eqnarray}\label{lemma18}
\Vert x \Vert_{t_{1},t_{2}} &\leq&\frac{C}{A} \Vert
w \Vert_{t_{1},t_{2}}+\frac{D}{A}\nonumber\\
&=&\frac{1+ {2 ( c+1 ) }/{a}%
+ {c}/{a^{2}}}{1- {c}/{a^{2}}} \Vert w \Vert_{t_{1},t_{2}}+%
\frac{ {2 ( c+1 ) }/{a}}{1- {c}/{a^{2}}}l (
t_{2}-t_{1} )
\\
&=&\frac{a^{2}+2ac+2a+c}{a^{2}-c} \Vert w \Vert_{t_{1},t_{2}}+%
\frac{2a ( c+1 ) }{a^{2}-c}l ( t_{2}-t_{1} ) ,\nonumber
\end{eqnarray}
and that (\ref{lemma13}) implies%
%
\begin{eqnarray} \label{lemma17}\qquad
\Vert x \Vert_{t_{1},t_{2}}
&\leq&\frac{C}{A} \Vert
w \Vert_{t_{1},t_{2}}+\frac{D}{A}=\frac{1+ {2}/{a}}{1-
{2c}/{a}}%
\Vert w \Vert_{t_{1},t_{2}}+\frac{ {2 ( c+1 ) }/{a}}{%
1- {2c}/{a}}l ( t_{2}-t_{1} ) \nonumber
\\[-8pt]
\\[-8pt]\qquad
&=&
\frac{a+2}{a-2c} \Vert w \Vert_{t_{1},t_{2}}+\frac{2 (
c+1 ) }{a-2c}l ( t_{2}-t_{1} ) .\nonumber
\end{eqnarray}
From the definition $e=\frac{c}{a^{2}\vee a/2}$ we know that if $
{a}/{2%
}\leq a^{2}$, then we can set $c=a^{2}e$ in (\ref{lemma18}) and obtain%
%
\begin{eqnarray}\label{haha14}
\Vert x \Vert_{t_{1},t_{2}}
&\leq&\frac{a^{2}+2a^{3}e+2a+a^{2}e}{%
a^{2} ( 1-e ) } \Vert w \Vert_{t_{1},t_{2}}+\frac{%
2a^{3}e+2a}{a^{2} ( 1-e ) }l ( t_{2}-t_{1} ) \nonumber
\\[-8pt]
\\[-8pt]
&=&\frac{a+2a^{2}e+2+ae}{a ( 1-e ) } \Vert w \Vert
_{t_{1},t_{2}}+\frac{2a^{2}e+2}{a ( 1-e ) }l (
t_{2}-t_{1} ) ,\nonumber
\end{eqnarray}
whereas if $ {a}/{2}\geq a^{2}$, then we can set $2c=ae$ in (\ref%
{lemma17}) and obtain%
%
\begin{equation}
\Vert x \Vert_{t_{1},t_{2}}\leq\frac{a+2}{a ( 1-e ) }%
\Vert w \Vert_{t_{1},t_{2}}+\frac{ae+2}{a ( 1-e ) }%
l ( t_{2}-t_{1} ) . \label{haha15}
\end{equation}
Hence, in either case, we arrive at
%
\begin{equation}
\Vert x \Vert_{t_{1},t_{2}}\leq\frac{a+2a^{2}e+2+ae}{a (
1-e ) } \Vert w \Vert_{t_{1},t_{2}}+\frac{2a^{2}e+2+ae}{%
a ( 1-e ) }l ( t_{2}-t_{1} ) , \label{haha16}
\end{equation}
and the proof of estimate (\ref{lemma1r1}) is complete. Finally, we note
that estimate (\ref{lemma1r2}) now follows directly from (\ref{lemma1r1})
and (\ref{lemma11}). This completes the proof of Lemma~\ref{compactlemma}.
\end{pf*}

%

\begin{pf*}{Proof of Theorem \ref{compacttheorem}}
We first note that
the assumptions stated in Theorem \ref{compacttheorem} ensure that there
exist some $a>0$ and $0<e<1$ such that $a_{s,x_{s}} ( \rho_{0},\eta
_{0} ) \geq a$ and $e_{s,x_{s}} ( \rho_{0},\eta_{0} ) \leq e$
for all $x_{s}\in\partial D_{s}$, $s\in [ 0,T ] $. Next we
recursively define two sets of time-points $\{\hat{T}_{i}\}$ and $\{
T_{i}\}$%
. In particular, we let $\hat{T}_{0}=T_{0}=0$ and define, for $i\geq
0$, $%
T_{i+1}=T$ if $x_{t}\in D_{t}$ for all $t\in\lbrack0,T]$ and
%
\begin{equation}
T_{i+1}=\inf \{ t\dvtx  \hat{T}_{i}\leq t\leq T, x_{t}\in\partial
D_{t} \} , \label{Tdef}
\end{equation}
otherwise. Similarly, for $i\geq0$, we let $\hat
{T}_{i+1}=(T_{i+1}+\eta
_{0})\wedge T$, if $ \vert x_{t}-x_{T_{i+1}} \vert+l (
t-T_{i+1} ) <\rho_{0}$ for all $t$ such that $T_{i+1}\leq t\leq
(T_{i+1}+\eta_{0})\wedge T$, and
%
\begin{equation}\hspace*{30pt}
\hat{T}_{i+1}=\inf \{ T_{i+1}\leq t<(T_{i+1}+\eta_{0})\wedge
T\dvtx  \vert x_{t}-x_{T_{i+1}} \vert+l ( t-T_{i+1} ) \geq
\rho_{0} \} , \label{taudef}
\end{equation}
otherwise. Using \eqref{limitzero} and the fact that $x$ is a right
continuous function,\vspace*{1pt} it follows that $T_{i+1}<\hat{T}_{i+1}$ for all
$i\geq
0 $. Moreover, using (\ref{Tdef})--(\ref{taudef}), we can apply Lemma~\ref%
{compactlemma} to any pair of time points $(t_{1},t_{2})$ such that $%
T_{i}\leq t_{1}\leq t_{2}<\hat{T}_{i}$ and obtain%
%
\begin{eqnarray}\label{hahaha1}
\Vert x \Vert_{t_{1},t_{2}} &\leq&
K_{1} ( a,e )
\Vert w \Vert_{t_{1},t_{2}}+K_{2} ( a,e ) l (
t_{2}-t_{1} ) ,\nonumber
\\[-8pt]
\\[-8pt]
|\lambda|_{t_{2}}-|\lambda|_{t_{1}} &\leq&
K_{3} ( a,e )
\Vert w \Vert_{t_{1},t_{2}}+K_{4} ( a,e ) l (
t_{2}-t_{1} ) ,\nonumber
\end{eqnarray}
whenever $T_{i}\leq t_{1}\leq t_{2}<\hat{T}_{i}$ where $K_{1}$,
$K_{2}$, $%
K_{3}$ and $K_{4}$ are defined as in Lemma~\ref{compactlemma} based on $a$
and $e$ introduced above. Next, we want to find a similar estimate
whenever $%
\hat{T}_{i}\leq t_{1}\leq t_{2}<T_{i+1}$. If $\hat{T}_{i}=T_{i+1}$,
we are
done and, hence, we assume that $\hat{T}_{i}<T_{i+1}$. In that case
$x_{t}\in
D_{t}$ for all $\hat{T}_{i}\leq t<T_{i+1}$ and, as a consequence, the
changes in $x$ and $w$ coincide on this time interval. Finally, considering
the case $T_{i}\leq t_{1}<\hat{T}_{i}\leq t_{2}<T_{i+1}$, we have%
%
\begin{eqnarray} \label{hahaha2}
|x_{t_{2}}-x_{t_{1}}| &\leq&|w_{t_{2}}-w_{\hat{T}_{i}}|+|x_{\hat
{T}_{i}}-x_{%
\hat{T}_{i}^{-}}|+|x_{\hat{T}_{i}^{-}}-x_{t_{1}}|, \nonumber
\\[-8pt]
\\[-8pt]
|\lambda|_{t_{2}}-|\lambda|_{t_{1}} &\leq&(|\lambda|_{\hat{T}%
_{i}}-|\lambda|_{\hat{T}_{i}^{-}})+(|\lambda|_{\hat
{T}_{i}^{-}}-|\lambda
|_{t_{1}}).\nonumber
\end{eqnarray}
The terms $|x_{\hat{T}_{i}^{-}}-x_{t_{1}}|$ and $|\lambda|_{\hat{T}%
_{i}^{-}}-|\lambda|_{t_{1}}$ in \eqref{hahaha2} can be handled\vspace*{1pt} using %
\eqref{hahaha1}. Regarding the terms $|x_{\hat{T}_{i}}-x_{\hat{T}%
_{i}^{-}}| $ and $|\lambda|_{\hat{T}_{i}}-|\lambda|_{\hat{T}_{i}^{-}}$
in \eqref{hahaha2}, we can, since $ \vert x_{\hat{T}_{i}}-x_{\hat{T}%
_{i}^{-}} \vert\leq\rho_{0}$, use (\ref{cdef2}) and the definition
of the Skorohod problem to first conclude that%
%
\begin{eqnarray}
\vert w_{\hat{T}_{i}}-w_{\hat{T}_{i}^{-}} \vert^{2}
&=& \vert
x_{\hat{T}_{i}}-x_{\hat{T}_{i}^{-}} \vert^{2}+ \bigl\vert\gamma_{\hat{%
T}_{i}} ( \vert\lambda \vert_{\hat{T}_{i}}- \vert
\lambda \vert_{\hat{T}_{i}^{-}} ) \bigr\vert^{2} \nonumber\\
&&{}-2 ( x_{\hat{T}_{i}}-x_{\hat{T}_{i}^{-}} ) \cdot\gamma_{\hat{T}%
_{i}} ( \vert\lambda \vert_{\hat{T}_{i}}- \vert\lambda
\vert_{\hat{T}_{i}^{-}} ) \nonumber
\\[-8pt]
\\[-8pt]
&\geq& \vert x_{\hat{T}_{i}}-x_{\hat{T}_{i}^{-}} \vert^{2}+ (
\vert\lambda \vert_{\hat{T}_{i}}- \vert\lambda \vert
_{\hat{T}_{i}^{-}} ) ^{2} \nonumber\\
&&{}
-2c_{x_{\hat{T}_{i}},\hat{T}_{i}} ( \rho_{0},\eta_{0} )
\vert x_{\hat{T}_{i}}-x_{\hat{T}_{i}^{-}} \vert ( \vert
\lambda \vert_{\hat{T}_{i}}- \vert\lambda \vert_{\hat{T}%
_{i}^{-}} ) .\nonumber
\end{eqnarray}
Then
%
\begin{eqnarray}\label{hahaha4}
\vert w_{\hat{T}_{i}}-w_{\hat{T}_{i}^{-}} \vert^{2} &\geq& \bigl(
1-c_{x_{\hat{T}_{i}},\hat{T}_{i}} ( \rho_{0},\eta_{0} ) \bigr)
\vert x_{\hat{T}_{i}}-x_{\hat{T}_{i}^{-}} \vert^{2}\nonumber
\\[-8pt]
\\[-8pt]
&&{}+ \bigl( 1-c_{x_{\hat{T}_{i}},\hat{T}_{i}} ( \rho_{0},\eta_{0} )
\bigr) ( \vert\lambda \vert_{\hat{T}_{i}}- \vert
\lambda \vert_{\hat{T}_{i}^{-}} ) ^{2},\nonumber
\end{eqnarray}
and, as $ ( a_{s,y} ( \rho_{0},\eta_{0} ) ) ^{2}\vee
a_{s,y} ( \rho_{0},\eta_{0} ) /2\leq1$, for all $y\in\partial
D_{s}$, $s\in [ 0,T ] $, we obtain%
%
\begin{eqnarray}\label{hahaha5}
c_{x_{\hat{T}_{i}},\hat{T}_{i}} ( \rho_{0},\eta_{0} ) &\leq&\frac{%
c_{x_{\hat{T}_{i}},\hat{T}_{i}} ( \rho_{0},\eta_{0} ) }{ (
a_{x_{\hat{T}_{i}},\hat{T}_{i}} ( \rho_{0},\eta_{0} ) )
^{2}\vee a_{x_{\hat{T}_{i}},\hat{T}_{i}} ( \rho_{0},\eta_{0} )
/2}\nonumber
\\[-8pt]
\\[-8pt]
&=&e_{x_{\hat{T}_{i}},\hat{T}_{i}} ( \rho_{0},\eta_{0} ) \leq e.\nonumber
\end{eqnarray}
Combining the estimates in \eqref{hahaha4} and \eqref{hahaha5}, we
arrive at
%
\begin{eqnarray}\label{theorem14}\qquad
\vert x_{\hat{T}_{i}}-x_{\hat{T}_{i}^{-}} \vert&\leq&\frac{1}{%
\sqrt{1-c_{x_{\hat{T}_{i}},\hat{T}_{i}} ( \rho_{0} ) }} \vert
w_{\hat{T}_{i}}-w_{\hat{T}_{i}^{-}} \vert\leq\frac{1}{\sqrt{1-e}}%
\vert w_{\hat{T}_{i}}-w_{\hat{T}_{i}^{-}} \vert, \nonumber
\\[-8pt]
\\[-8pt]\qquad
\vert\lambda \vert_{\hat{T}_{i}}- \vert\lambda \vert
_{\hat{T}_{i}^{-}} &\leq&\frac{1}{\sqrt{1-e}} \vert w_{\hat
{T}_{i}}-w_{%
\hat{T}_{i}^{-}} \vert.\nonumber
\end{eqnarray}
Introducing the notation
%
\begin{eqnarray}\label{blablabla}
K_{1} &=&K_{1} ( a,e ) +1+\frac{1}{\sqrt{1-e}},\qquad
K_{2}=K_{2} ( a,e ) , \nonumber
\\[-8pt]
\\[-8pt]
K_{3} &=&K_{3} ( a,e ) +\frac{1}{\sqrt{1-e}},\qquad
K_{4}=K_{4} ( a,e ) ,\nonumber
\end{eqnarray}
we can use the deductions in \eqref{hahaha1}--\eqref{theorem14} to conclude
that
%
\begin{eqnarray}\label{hahaha1again}
\Vert x \Vert_{t_{1},t_{2}} &\leq&K_{1} \Vert w \Vert
_{t_{1},t_{2}}+K_{2}l ( t_{2}-t_{1} ) ,\nonumber
\\[-8pt]
\\[-8pt]
|\lambda|_{t_{2}}-|\lambda|_{t_{1}} &\leq&K_{3} \Vert w \Vert
_{t_{1},t_{2}}+K_{4}l ( t_{2}-t_{1} ) ,\nonumber
\end{eqnarray}
whenever $T_{i}\leq t_{1}\leq t_{2}<T_{i+1}$. We now intend to make use of
the estimates in~\eqref{hahaha1again} to complete the proof of Theorem
\ref%
{compacttheorem}. Note that above we have constructed a set of
time-points $%
\{T_{i}\}_{i=0}^{M+1}$, where $M$ is so far undetermined, and
%
\begin{equation}
0=T_{0}<T_{1}<\cdots<T_{M}<T=T_{M+1}. \label{blabla2}
\end{equation}
If $M\geq1$, let
%
\begin{equation}\hspace*{30pt}
0\leq t_{1}\leq u\leq r\leq t_{2}\leq T,\qquad T_{h-1}\leq u<T_{h},\
T_{v}\leq r<T_{v+1},\  h-1\leq v. \label{blabla3}
\end{equation}
Then, using \eqref{hahaha1again}, we have
%
\begin{eqnarray}\label{theorem12}\qquad
\vert x_{r}-x_{u} \vert&\leq& ( M+1 ) \bigl(
K_{1} \Vert w \Vert_{t_{1},t_{2}}+K_{2}l ( t_{2}-t_{1} )
\bigr) +\sum_{i=h}^{v} \vert x_{T_{i}}-x_{T_{i}^{-}} \vert,
\nonumber
\\[-8pt]
\\[-8pt]\qquad
|\lambda|_{r}-|\lambda|_{u} &\leq& ( M+1 )\bigl (
K_{3} \Vert w \Vert_{t_{1},t_{2}}+K_{4}l ( t_{2}-t_{1} )
\bigr) +\sum_{i=h}^{v} \vert\lambda \vert_{T_{i}}- \vert
\lambda \vert_{T_{i}^{-}}.\nonumber
\end{eqnarray}
Moreover, arguing exactly as in the deduction of \eqref{theorem14}, we
obtain
%
\begin{eqnarray} \label{theorem14+}
\vert x_{T_{i}}-x_{T_{i}^{-}} \vert&\leq&\frac{1}{\sqrt{1-e}}%
\vert w_{T_{i}}-w_{T_{i}^{-}} \vert, \nonumber
\\[-8pt]
\\[-8pt]
\vert\lambda \vert_{T_{i}}- \vert\lambda \vert
_{T_{i}^{-}} &\leq&\frac{1}{\sqrt{1-e}} \vert
w_{T_{i}}-w_{T_{i}^{-}} \vert,\nonumber
\end{eqnarray}
whenever $1\leq i\leq M+1$. Hence, to complete the proof of Theorem
\ref%
{compacttheorem}, we have to estimate $M$. To do this, we consider
$\mathcal{W}%
\subset\mathcal{D} ( [ 0,T ] ,%
\mathbb{R}
^{d} ) $, which is assumed to be relatively compact in the Skorohod
topology and for which $w_{0}\in\overline{D_{0}}$ whenever $w\in
\mathcal{W}
$. We shall prove that the $M$ introduced above is bounded for every such
set $\mathcal{W}$. To do this, we use the notation introduced in
Section \ref%
{skorohodtop} concerning the Skorohod topology. In the following let
$\delta
^{\prime}$ be a fixed number such that
%
\begin{equation}\quad
\delta^{\prime}=\min\{\eta_{0},\hat{\delta}{}^{\prime}\}  \qquad\mbox{where }%
\hat{\delta}{}^{\prime}\mbox{ is such that }l(\hat{\delta}{}^{\prime
})\leq
\rho_{0}/ \bigl( 2 ( K_{2}+1 ) \bigr) . \label{ahah1}
\end{equation}
Note that the existence of $\hat{\delta}{}^{\prime}$ follows
immediately from %
\eqref{limitzero}. Using the definition of $\delta^{\prime}$ and the fact
that $\mathcal{W}\subset\mathcal{D} ( [ 0,T ] ,%
\mathbb{R}
^{d} ) $ is relatively compact in the Skorohod topology, we see, by (%
\ref{e11}) and (\ref{d11+}), that
%
\begin{equation}
\lim_{\delta\rightarrow0}\tilde{\mu} ( \mathcal{W},\delta,\delta
^{\prime},T ) =0. \label{gaga}
\end{equation}
In particular, using \eqref{gaga}, we can find a $0<\delta<\delta
^{\prime
} $ such that for every $w\in\mathcal{W}$ there exists a partition $
\{
t_{j} \} _{j=0}^{M}$, in general depending on $w$, such that
%
\begin{equation}
\delta<|t_{j+1}-t_{j}|<\delta^{\prime}\qquad \mbox{for }j\in \{
0,\ldots,M-1 \}  \label{h11-}
\end{equation}
and
%
\begin{equation}
\max_{0\leq j\leq M-1}\sup_{u,r\in [ t_{j},t_{j+1} ) } \vert
w_{u}-w_{r} \vert<\frac{\rho_{0}}{2K_{1}}. \label{h11}
\end{equation}
We claim that none of the intervals $\{[t_{j},t_{j+1})\}$ in this partition
can contain more than one point from the sequence $\{T_{i}\}$. To prove this,
we suppose, on the contrary, that there exist $i$ and $j$ such that $%
t_{j}\leq T_{i}<T_{i+1}<t_{j+1}$. Then, by construction,
%
\begin{equation}
t_{j}\leq T_{i}<\hat{T}_{i}\leq T_{i+1}<t_{j+1}. \label{aassump}
\end{equation}
We intend to estimate $ \vert x_{\hat{T}_{i}}-x_{T_{i}} \vert+l(%
\hat{T}_{i}-T_{i})$. We first note that if $ \vert
x_{t}-x_{T_{i}} \vert+l(t-T_{i})<\rho_{0}$ for all $t$ such that $%
T_{i}\leq t\leq(T_{i}+\eta_{0})\wedge T$, then $\hat
{T}_{i}=(T_{i}+\eta
_{0})\wedge T$. However, using \eqref{ahah1} and \eqref{h11-}, it is clear
that neither $\hat{T}_{i}=(T_{i}+\eta_{0})$ nor $\hat{T}_{i}=T$ can occur.
Hence, we can assume that $\hat{T}_{i}$ is given by \eqref{taudef}
and, as a
consequence, that
%
\begin{equation}
\vert x_{\hat{T}_{i}}-x_{T_{i}} \vert+l(\hat{T}_{i}-T_{i})\geq
\rho_{0}. \label{aassump+}
\end{equation}
But on the other hand, using \eqref{hahaha1again}, we first see that
%
\begin{eqnarray}
\vert x_{\hat{T}_{i}}-x_{T_{i}} \vert+l(\hat{T}_{i}-T_{i})&\leq&
\Vert x \Vert_{{T_{i}},{\hat{T}_{i}}}+l(\hat{T}_{i}-T_{i})\nonumber
\\[-8pt]
\\[-8pt]&\leq&
K_{1} \Vert w \Vert_{{T_{i}},{\hat{T}_{i}}}+ ( K_{2}+1 )
l(\hat{T}_{i}-T_{i}),\nonumber
\end{eqnarray}
and then, using \eqref{ahah1}, \eqref{h11-} and \eqref{h11}, we deduce
%
\begin{equation}
\vert x_{\hat{T}_{i}}-x_{T_{i}} \vert+l(\hat{T}_{i}-T_{i})<K_{1}%
\frac{\rho_{0}}{2K_{1}}+ ( K_{2}+1 ) l(\delta^{\prime})<\rho
_{0},
\end{equation}
which contradicts the assumption $t_{j}\leq T_{i}<T_{i+1}<t_{j+1}$. Hence,
none of the intervals $\{[t_{j},t_{j+1})\}$ in the partition can contain
more than one point from the sequence $\{T_{i}\}$ and, in particular, we
conclude that
%
\begin{equation}
M\leq\frac{T}{\delta}+1. \label{theorem15}
\end{equation}
Combining \eqref{theorem12}, \eqref{theorem14+} and (\ref
{theorem15}), we
see that%
%
\begin{eqnarray}
\Vert x \Vert_{t_{1},t_{2}}
&\leq&
 ( M+1 ) \bigl(
K_{1} \Vert w \Vert_{t_{1},t_{2}}+K_{2}l ( t_{2}-t_{1} )
\bigr) +M|x_{T_{i}}-x_{T_{i}^{-}}|
\nonumber\\
&\leq&
\biggl( \frac{T}{\delta}+2 \biggr) \bigl( K_{1} \Vert w \Vert
_{t_{1},t_{2}}+K_{2}l ( t_{2}-t_{1} ) \bigr) + \biggl( \frac{T}{%
\delta}+1 \biggr) \frac{1}{\sqrt{1-e}} \Vert w \Vert_{t_{1},t_{2}}
\nonumber
\\[-8pt]
\\[-8pt]
&\leq&
\biggl( K_{1} \biggl( \frac{T}{\delta}+2 \biggr) +\frac{1}{\sqrt{1-e}}%
\biggl( \frac{T}{\delta}+1 \biggr) \biggr) \Vert w \Vert
_{t_{1},t_{2}} \nonumber\\
&&{}+K_{2} \biggl( \frac{T}{\delta}+2 \biggr) l ( t_{2}-t_{1} ) ,\nonumber
\end{eqnarray}
and, similarly, that
%
\begin{eqnarray}
\vert\lambda \vert_{t_{2}}- \vert\lambda \vert
_{t_{1}}&\leq& \biggl( K_{3} \biggl( \frac{T}{\delta}+2 \biggr) +\frac{1}{\sqrt{%
1-e}} \biggl( \frac{T}{\delta}+1 \biggr) \biggr) \Vert w \Vert
_{t_{1},t_{2}} \nonumber
\\[-8pt]
\\[-8pt]
&&{}+
K_{4} \biggl( \frac{T}{\delta}+2 \biggr) l ( t_{2}-t_{1} ).\nonumber
\end{eqnarray}
The deductions in the last two displays complete the proof of Theorem
\ref%
{compacttheorem}.
\end{pf*}

\subsection{Estimates for approximations to Skorohod problems}\label{SPappest}

Let $T>0$, $D\subset
\mathbb{R}
^{d+1}$ and $\Gamma=\Gamma_{t} ( z ) $ satisfy the assumptions
stated in Theorem \ref{compacttheorem}. In this section we derive estimates
for approximations to the Skorohod problem for $(D,\Gamma,w)$. In
particular, given $w\in\mathcal{D} ( [ 0,T ] ,%
\mathbb{R}
^{d} ) $ with $w_{0}\in\overline{D_{0}}$ and a partition $\{\tau
_{k}\}_{k=0}^{N}$ of the interval $[0,T]$, which we denote by $\Delta
$, we
define $w^{\Delta}$ as in \eqref{approxw}. Recall that $\Delta^{\ast}$
was defined in \eqref{korlim+}. Furthermore, in the following, we will
assume that \eqref{korlim} holds whenever $k\in\{1,\ldots,N\}$.
Based on the
assumption in \eqref{korlim}, we define $D^{\Delta}$, $\Gamma
^{\Delta}$, $%
x^{\Delta}$ and $\lambda^{\Delta}$ as in \eqref{approxw+}, %
\eqref{approxw++} and \eqref{approxw+++}. Then, by construction, the
pair $%
(x^{\Delta},\lambda^{\Delta})$ is a solution to the Skorohod problem
for $%
( D^{\Delta},\Gamma^{\Delta},w^{\Delta} ) $. In this section
we prove the following results.

\begin{lemma}
\label{compactlemmaapprox} Let $T>0$, $D\subset
\mathbb{R}
^{d+1}$, $r_{0}$, $\Gamma=\Gamma_{t} ( z ) $, $0<\rho_{0}<r_{0}$%
, $\eta_{0}>0$, $a$, $e$, $\delta_{0}$ and $h_{0}$ be as in the statement
of Theorem \ref{compacttheorem}. Given $a>0$ and $e\in(0,1)$, let the
functions $K_{1}$, $K_{2}$, $K_{3}$ and $K_{4}$ be defined as in %
\eqref{functions} and let $w\in\mathcal{D} ( [ 0,T ] ,%
\mathbb{R}
^{d} ) $ with $w_{0}\in\overline{D_{0}}$. Let $\Delta=\{\tau
_{k}\}_{k=0}^{N}$ be a partition of the interval $[0,T]$, let
$w^{\Delta}$
be defined as in \eqref{approxw} and assume that \eqref{korlim}
holds. Given
$\Delta$ and $w^{\Delta}$, let $D^{\Delta}$, $\Gamma^{\Delta}$, $%
x^{\Delta}$ and $\lambda^{\Delta}$ be defined as in \eqref
{approxw+}, %
\eqref{approxw++} and \eqref{approxw+++}. Consider a fixed but
arbitrary $s\in
\lbrack0,T]$, such that $x_{s}^{\Delta}\in\partial D_{s}^{\Delta}$.
Then, for $0\leq s\leq t_{1}\leq t_{2}<\tau_{\rho_{0},\eta
_{0}}^{\Delta}$%
, $(w^{\Delta},x^{\Delta},\lambda^{\Delta})$ satisfies the estimates
%
\begin{eqnarray} \label{lemma1r1+}\quad
\Vert x^{\Delta} \Vert_{t_{1},t_{2}} &\leq&K_{1} (
a,e ) \Vert w \Vert_{t_{1},t_{2}}+K_{2} ( a,e )
\bigl(l ( t_{2}-t_{1} ) +l(\Delta^{\ast})\bigr), \\\label{lemma1r2+}
\quad
\vert\lambda^{\Delta} \vert_{t_{2}}- \vert\lambda
^{\Delta} \vert_{t_{1}} &\leq&K_{3} ( a,e ) \Vert
w \Vert_{t_{1},t_{2}}+K_{4} ( a,e ) \bigl(l (
t_{2}-t_{1} ) +l(\Delta^{\ast})\bigr).
\end{eqnarray}
Here $\tau_{\rho_{0},\eta_{0}}^{\Delta}$ is defined as follows. If there
exists some $t$ such that $s\leq t<(s+\eta_{0})\wedge T$ and $ \vert
x_{t}^{\Delta}-x_{s}^{\Delta} \vert+l ( t-s ) +l(\Delta
^{\ast})\geq\rho_{0}$, then
%
\begin{equation}\hspace*{32pt}
\tau_{\rho_{0},\eta_{0}}^{\Delta}=\inf\{t\dvtx  s\leq t<(s+\eta
_{0})\wedge
T, \vert x_{t}^{\Delta}-x_{s}^{\Delta} \vert+l ( t-s )
+l(\Delta^{\ast})\geq\rho_{0}\},
\end{equation}
whereas if $ \vert x_{t}^{\Delta}-x_{s}^{\Delta} \vert+l (
t-s ) +l(\Delta^{\ast})<\rho_{0}$ for all $s\leq t<(s+\eta
_{0})\wedge T$, then
%
\begin{equation}
\tau_{\rho_{0},\eta_{0}}^{\Delta}=(s+\eta_{0})\wedge T.
\end{equation}
\end{lemma}

\begin{theorem}
\label{compacttheoremapprox} Let $T>0$, $D\subset
\mathbb{R}
^{d+1}$, $r_{0}$, $\Gamma=\Gamma_{t} ( z ) $, $0<\rho_{0}<r_{0}$%
, $\eta_{0}>0$, $a$, $e$, $\delta_{0}$ and $h_{0}$ be as in the statement
of Theorem \ref{compacttheorem}. Given $a>0$ and $e\in(0,1)$, let the
functions $K_{1}$, $K_{2}$, $K_{3}$ and $K_{4}$ be defined as in %
\eqref{functions} and let $w\in\mathcal{D} ( [ 0,T ] ,%
\mathbb{R}
^{d} ) $ with $w_{0}\in\overline{D_{0}}$. Let $\Delta=\{\tau
_{k}\}_{k=0}^{N}$ be a partition of the interval $[0,T]$, let
$w^{\Delta}$
be defined as in \eqref{approxw} and assume that \eqref{korlim}
holds. Let $%
\Delta$ be such that $l(\Delta^{\ast})\leq\rho_{0}/(4(K_{2}(a,e)+1))$
and let
%
\begin{eqnarray}\label{ahah1+}
\delta^{\prime}&=&\min\{ \eta_{0},\hat{\delta}{}^{\prime}\} \qquad
\mbox{where }%
\hat{\delta}{}^{\prime}\mbox{ is such that}\nonumber
\\[-8pt]
\\[-8pt]
l( \hat{\delta}{}^{\prime})
+l(\Delta^{\ast})&\leq&\rho_{0}/\bigl(2\bigl(K_{2}(a,e)+1\bigr)\bigr).\nonumber
\end{eqnarray}
Given $\Delta$ and $w^{\Delta}$, let $D^{\Delta}$, $\Gamma^{\Delta
}$, $%
x^{\Delta}$ and $\lambda^{\Delta}$ be defined as in \eqref
{approxw+}, %
\eqref{approxw++} and \eqref{approxw+++}. Moreover, assume that
$x^{\Delta}\in
\mathcal{D}^{\rho_{0}} ( [ 0,T ] ,%
\mathbb{R}
^{d} ) $. Then there exist positive constants $\hat{L}_{1} (
w,T ) $ , $\hat{L}_{2} ( w,T ) $, $\hat{L}_{3} (
w,T ) $ and $\hat{L}_{4} ( w,T ) $, independent of $\Delta$,
such that
%
\begin{eqnarray}\quad
\Vert x^{\Delta} \Vert_{t_{1},t_{2}} &\leq&\hat{L}_{1} (
w,T ) \Vert w \Vert_{t_{1},t_{2}}+\hat{L}_{2} (
w,T ) \bigl(l ( t_{2}-t_{1} ) +l(\Delta^{\ast})\bigr),\nonumber
\\[-8pt]
\\[-8pt]\quad
|\lambda^{\Delta}|_{t_{2}}-|\lambda^{\Delta}|_{t_{1}} &\leq&\hat
{L}%
_{3} ( w,T ) \Vert w \Vert_{t_{1},t_{2}}+\hat{L}%
_{4} ( w,T ) \bigl(l ( t_{2}-t_{1} ) +l(\Delta^{\ast})\bigr),\nonumber
\end{eqnarray}
whenever $0\leq t_{1}\leq t_{2}\leq T$.
\end{theorem}

\begin{pf*}{Proof of Lemma \ref{compactlemmaapprox}}
Naturally, the
proof of this lemma is similar to the proof of Lemma \ref
{compactlemma} and,
thus, we only describe the main differences compared to the proof of
Lemma %
\ref{compactlemma}. First we note, by the assumptions on $D$ and the
construction of $D^{\Delta}$ based on $D$, that there exists a unit
vector $%
u$ such that
%
\begin{equation}
\langle\gamma_{r}^{\Delta},u \rangle\geq a
\end{equation}
for all $\gamma_{r}^{\Delta}\in\Gamma_{r}^{\Delta,1} ( y ) $, $%
y\in\partial D_{r}^{\Delta}\cap\overline{B_{\rho_{0}} (
x_{s}^{\Delta} ) }$ and $r\in\lbrack s,\tau_{\rho_{0},\eta
_{0}}^{\Delta}]\subset [ s,(s+\eta_{0})\wedge T ] $. We also
note that if $t_{1}\in\lbrack\tau_{j},\tau_{j+1})$ and $t_{2}\in
\lbrack
\tau_{k},\tau_{k+1})$, for some $j,k\in\{0,\ldots,N-1\}$, then $%
|x_{t_{2}}^{\Delta}-x_{t_{1}}^{\Delta}|=0$ if $j=k$ and otherwise $%
|x_{t_{2}}^{\Delta}-x_{t_{1}}^{\Delta}|=|x_{\tau_{k}}^{\Delta
}-x_{\tau
_{j}}^{\Delta}|$. Now, using the fact that $(x^{\Delta},\lambda
^{\Delta
}) $ solves the Skorohod problem for $ ( D^{\Delta},\Gamma^{\Delta
},w^{\Delta} ) $, we conclude, in analogy with \eqref{hahaha}, that
%
\begin{equation}
\langle x_{t_{2}}^{\Delta}-x_{t_{1}}^{\Delta},u \rangle
= \langle w_{t_{2}}^{\Delta}-w_{t_{1}}^{\Delta},u \rangle
+\int_{t_{1}^{+}}^{t_{2}^{+}}
\underbrace{ \langle
\gamma_{r}^{\Delta},u \rangle}_{\geq a}\,
d \vert\lambda^{\Delta
} \vert_{r}  \label{haha+}
\end{equation}
for any $0\leq s\leq t_{1}\leq t_{2}<\tau_{\rho_{0},\eta
_{0}}^{\Delta}$,
where $\gamma_{r}^{\Delta}\in\Gamma_{r}^{\Delta,1} ( y ) $ for
some $y\in\partial D_{r}^{\Delta}$. Based on \eqref{haha+}, we obtain
%
\begin{equation}
\vert\lambda^{\Delta} \vert_{t_{2}}- \vert\lambda
^{\Delta} \vert_{t_{1}}\leq\frac{1}{a} ( \vert
w_{t_{2}}^{\Delta}-w_{t_{1}}^{\Delta} \vert+ \vert
x_{t_{2}}^{\Delta}-x_{t_{1}}^{\Delta} \vert ) .
\label{lemma11+}
\end{equation}
Furthermore, arguing as in the proof of Lemma \ref{compactlemma}, we derive
%
\begin{eqnarray} \label{lemma15+}
\vert x_{t_{2}}^{\Delta}-x_{t_{1}}^{\Delta} \vert^{2}
&=& \vert w_{t_{2}}^{\Delta}-w_{t_{1}}^{\Delta} \vert
^{2}+2\int_{t_{1}^{+}}^{t_{2}^{+}} \langle w_{t_{2}}^{\Delta
}-w_{r}^{\Delta},\gamma_{r}^{\Delta} \rangle\, d \vert\lambda
^{\Delta} \vert_{r} \nonumber\\
&&{}+2\int_{t_{1}^{+}}^{t_{2}^{+}} \langle x_{r}^{\Delta
}-x_{t_{1}}^{\Delta},\gamma_{r}^{\Delta} \rangle\, d \vert\lambda
^{\Delta} \vert_{r}+ \biggl( \int_{t_{1}^{+}}^{t_{2}^{+}}\gamma
_{r}^{\Delta}\,d \vert\lambda^{\Delta} \vert_{r} \biggr) ^{2}
\\
&&{}-2\int_{t_{1}^{+}}^{t_{2}^{+}} \biggl\langle \biggl(
\int_{t_{1}{}^{+}}^{r^{+}}\gamma_{u}^{\Delta}\,d \vert\lambda
^{\Delta
} \vert_{u} \biggr) ,\gamma_{r}^{\Delta} \biggr\rangle\, d \vert
\lambda^{\Delta} \vert_{r}.\nonumber
\end{eqnarray}
As in the proof of Lemma \ref{compactlemma}, we have to find upper
bounds of
all integrals in~(\ref{lemma15+}) and, naturally, particular attention has
to be paid to the second integral,\vspace*{1pt} as $x_{t_{1}}^{\Delta}$ might not belong
to $\overline{D_{r}^{\Delta}}$. For $y\in\partial D_{t}^{\Delta}$,
$t\in
\lbrack0,T]$, we let $N_{t}^{\Delta}(y)$ denote the set of inward normals
at $y\in\partial D_{t}^{\Delta}$ and given $r\in\lbrack0,T]$, $y\in
\mathbb{R}
^{d}\setminus\overline{D_{r}^{\Delta}}$ in the following we denote a
projection of $y$ onto $\partial D_{r}^{\Delta}$ along $N_{r}^{\Delta
}$ by
$\pi_{\partial D_{r}^{\Delta}}^{N_{r}^{\Delta}} ( y ) $.
Furthermore,\vspace*{-4pt} if $y\in\overline{D_{r}^{\Delta}}$, then we let $\pi
_{\partial D_{r}^{\Delta}}^{N_{r}^{\Delta}} ( y ) =y$. Using this
notation, and the definition of $\tau_{\rho_{0},\eta_{0}}^{\Delta
}$, we
see that
%
\begin{equation}
\vert\pi_{\partial D_{r}^{\Delta}}^{N_{r}^{\Delta}} (
x_{t_{1}}^{\Delta} ) -x_{s}^{\Delta} \vert\leq \vert
x_{t_{1}}^{\Delta}-x_{s}^{\Delta} \vert+l ( r-t_{1} )
+l(\Delta^{\ast})\leq\rho_{0}. \label{haha3+}
\end{equation}
Equation \eqref{haha3+} implies that $\pi_{\partial D_{r}^{\Delta
}}^{N_{r}^{\Delta
}} ( x_{t_{1}}^{\Delta} ) \in\overline{B_{\rho_{0}} (
x_{s}^{\Delta} ) }\cap\partial D_{r}^{\Delta}\subset\overline{%
B_{\rho_{0}} ( x_{s}^{\Delta} ) }\cap\overline{D_{r}^{\Delta}}$%
. Arguing as in~\eqref{haha4}--\eqref{haha7}, we then deduce that
%
\begin{eqnarray}
\langle x_{r}^{\Delta}-\pi_{\partial D_{r}^{\Delta}}^{N_{r}^{\Delta
}} ( x_{t_{1}}^{\Delta} ) ,\gamma_{r}^{\Delta} \rangle
&\leq&c \vert x_{r}^{\Delta}-\pi_{\partial D_{r}^{\Delta
}}^{N_{r}^{\Delta}} ( x_{t_{1}}^{\Delta} ) \vert \nonumber
\\[-8pt]
\\[-8pt]
&\leq&c \vert x_{r}^{\Delta}-x_{t_{1}}^{\Delta} \vert+cl (
r-t_{1} ) +cl(\Delta^{\ast}),\nonumber
\end{eqnarray}
and that
%
\begin{equation}
\langle\pi_{\partial D_{r}^{\Delta}}^{N_{r}^{\Delta}} (
x_{t_{1}}^{\Delta} ) -x_{t_{1}}^{\Delta},\gamma_{r}^{\Delta
} \rangle\leq \vert\pi_{\partial D_{r}^{\Delta
}}^{N_{r}^{\Delta}} ( x_{t_{1}}^{\Delta} ) -x_{t_{1}}^{\Delta
} \vert\leq l ( r-t_{1} ) +l(\Delta^{\ast}). \label{haha6+}
\end{equation}
Using the estimates derived above, we conclude that the second integral
in (%
\ref{lemma15+}) has the upper bound
%
\begin{equation}\hspace*{30pt}
2c\int_{t_{1}^{+}}^{t_{2}^{+}} \vert x_{r}^{\Delta}-x_{t_{1}}^{\Delta
} \vert \,d \vert\lambda^{\Delta} \vert_{r}+2 (
c+1 )\bigl (l ( t_{2}-t_{1} ) +l(\Delta^{\ast})\bigr) ( \vert
\lambda^{\Delta} \vert_{t_{2}}- \vert\lambda^{\Delta
} \vert_{t_{1}} ) . \label{lemma20+}
\end{equation}
Equipped with \eqref{lemma20+}, the proof of Lemma \ref{compactlemmaapprox}
can now be completed following the lines of the proof of Lemma \ref%
{compactlemma}.
\end{pf*}

%

\begin{pf*}{Proof of Theorem \ref{compacttheoremapprox}}
Proceeding
as in the proof of Theorem \ref{compacttheorem}, we first recursively define
two sets of time-points $\{\hat{T}{}_{i}^{\Delta}\}$ and $\{
T_{i}^{\Delta}\}$
in order to use Lemma~\ref{compactlemmaapprox}. In particular, we let
$\hat{T%
}{}_{0}^{\Delta}=T_{0}^{\Delta}=0$ and define, for $i\geq0$, $%
T_{i+1}^{\Delta}=T$ if $x_{t}^{\Delta}\in D_{t}^{\Delta}$ for all
$t\in
\lbrack0,T]$ and
%
\begin{equation}
T_{i+1}^{\Delta}=\inf \{ t\dvtx  \hat{T}{}_{i}^{\Delta}\leq t\leq
T, x_{t}^{\Delta}\in\partial D_{t}^{\Delta} \} ,
\end{equation}
otherwise. Similarly, for $i\geq0$ we let $\hat{T}{}_{i+1}^{\Delta
}=(T_{i+1}^{\Delta}+\eta_{0})\wedge T$, if $|x_{t}^{\Delta
}-x_{T_{i+1}}^{\Delta}|+l ( t-T_{i+1}^{\Delta} ) +l(\Delta^{\ast
})<\rho_{0}$ for all $t$ such that $T_{i+1}^{\Delta}\leq t\leq
(T_{i+1}^{\Delta}+\eta_{0})\wedge T$. Moreover, if the latter is not the
case, we then define $\hat{T}{}_{i+1}^{\Delta}$ to equal
%
\begin{equation}\label{taudef+}\hspace*{31pt}
\inf \{ T_{i+1}^{\Delta}\leq t<(T_{i+1}^{\Delta}+\eta_{0})\wedge
T\dvtx |x_{t}^{\Delta}-x_{T_{i+1}}^{\Delta}|+l ( t-T_{i+1}^{\Delta} )
+l(\Delta^{\ast})\geq\rho_{0} \}.
\end{equation}
We can then repeat the argument in \eqref{hahaha1}--\eqref
{hahaha1again} to
conclude that
%
\begin{eqnarray}\label{hahaha1again+}
\Vert x^{\Delta} \Vert_{t_{1},t_{2}} &\leq&K_{1} \Vert
w \Vert_{t_{1},t_{2}}+K_{2}\bigl(l ( t_{2}-t_{1} ) +l(\Delta
^{\ast})\bigr), \nonumber
\\[-8pt]
\\[-8pt]
|\lambda^{\Delta}|_{t_{2}}-|\lambda^{\Delta}|_{t_{1}} &\leq
&K_{3} \Vert w \Vert_{t_{1},t_{2}}+K_{4}\bigl(l (
t_{2}-t_{1} ) +l(\Delta^{\ast})\bigr),\nonumber
\end{eqnarray}
whenever $T_{i}^{\Delta}\leq t_{1}\leq t_{2}<T_{i+1}^{\Delta}$.
Furthermore, we note that the $M$ in $\{T_{i}\}_{i=0}^{M+1}$ is so far
undetermined and, as in \eqref{theorem12}--\eqref{theorem14+}, we derive
%
\begin{eqnarray}\label{theorem12+}
\Vert x^{\Delta} \Vert_{t_{1},t_{2}} &\leq& ( M+1 )
\bigl( K_{1} \Vert w \Vert_{t_{1},t_{2}}+K_{2}\bigl(l (
t_{2}-t_{1} ) +l(\Delta^{\ast})\bigr) \bigr)
\nonumber\\
&&{}+\sqrt{\frac{1}{1-e}}\sum_{i=1}^{M} \vert w_{T_{i}^{\Delta
}}-w_{T_{i}^{\Delta,-}} \vert,
\nonumber
\\[-8pt]
\\[-8pt]
|\lambda^{\Delta}|_{t_{2}}-|\lambda^{\Delta}|_{t_{1}} &\leq& (
M+1 )
\bigl( K_{3} \Vert w \Vert_{t_{1},t_{2}}+K_{4}\bigl(l (
t_{2}-t_{1} ) +l(\Delta^{\ast})\bigr) \bigr) \nonumber\\
&&{}+\sqrt{\frac{1}{1-e}}\sum_{i=1}^{M} \vert w_{T_{i}^{\Delta
}}-w_{T_{i}^{\Delta,-}} \vert,\nonumber
\end{eqnarray}
whenever $0\leq t_{1}\leq t_{2}\leq T$. To complete the proof of
Theorem \ref%
{compacttheoremapprox}, we can now proceed as in the proof of Theorem
\ref%
{compacttheorem} and conclude that $M\leq T/\delta+1$, where $\delta$ is
given in the proof of Theorem \ref{compacttheorem}. This completes the proof
of Theorem \ref{compacttheoremapprox}.
\end{pf*}

\section{Convergence and approximation of Skorohod problems} \label{sectskorgen}

In the first subsection of this section we prove the general convergence
result for sequences of Skorohod problems (see Theorem \ref{SPconvergence}
stated below) referred to in Section \ref{contri}. Then, in the second
subsection we explicitly construct, given $(D,\Gamma,w)$, an approximating
sequence $\{(D^{n},\Gamma^{n},w^{n})\}$ and, for each $n$, an explicit
solution $(x^{n},\lambda^{n})$ to the Skorohod problem for $ (
D^{n},\Gamma^{n},w^{n} ) $. We then prove that the constructed
sequence $\{(x^{n},\lambda^{n})\}$ of solutions converges to a
solution to
the Skorohod problem for $(D,\Gamma,w)$.

\subsection{Convergence of a sequence of solutions to Skorohod problems}

Let\vspace*{1pt} $T>0$ and let $D\subset
\mathbb{R}
^{d+1}$ be a time-dependent domain satisfying \eqref{timedep+}. Let
$\Gamma
=\Gamma_{t} ( z ) $ be a closed convex cone of vectors in $%
\mathbb{R}
^{d}$ for every $z\in\partial D_{t}$, $t\in\lbrack0,T]$ and assume
that $%
\Gamma$ satisfies \eqref{cone1} and \eqref{cone2}. Let $\{D^{n}\}
_{n=1}^{%
\infty}$ be a sequence of time-dependent domains $D^{n}\subset
\mathbb{R}
^{d+1}$ and let $\{\Gamma^{n}\}_{n=1}^{\infty}=\{\Gamma_{t}^{n} (
z ) \}_{n=1}^{\infty}$ be a sequence of\vspace*{1pt} closed convex cones of vectors
in $%
\mathbb{R}
^{d}$. Let $w\in\mathcal{D} ( [ 0,T ] ,%
\mathbb{R}
^{d} ) $ with $w_{0}\in\overline{D_{0}}$ and let $ \{
w^{n} \} $ with $w_{0}^{n}\in\overline{D_{0}^{n}}$ be a sequence of c%
\`{a}dl\`{a}g functions converging to $w$ in the Skorohod topology. Assume
that there exists a solution $ ( x^{n},\lambda^{n} ) $ to the
Skorohod problem for $ ( D^{n},\Gamma^{n},w^{n} ) $. Then in
Theorem \ref{SPconvergence} we prove, by making appropriate
assumptions on $%
D $, $\Gamma$, $\{D^{n}\}_{n=1}^{\infty}$ and $\{\Gamma
^{n}\}_{n=1}^{\infty}$, that if $D^{n}\rightarrow D$ and $\Gamma
^{n}\rightarrow\Gamma$ in the sense\vspace*{1pt} defined in Theorem \ref
{SPconvergence}%
, then $ ( x^{n},\lambda^{n} ) $ converges to $ ( x,\lambda
) $ and $(x,\lambda)$ is a solution to the Skorohod problem for $%
( D,\Gamma,w ) $. However, to state Theorem \ref{SPconvergence},
we need to introduce some additional notions and notation. In
particular, in
the following we let $a_{s,z}^{n}$ and $e_{s,z}^{n}$ be defined as in %
\eqref{adef} and \eqref{edef} but with respect to $(D^{n},\Gamma
^{n})$. We
assume that $D^{n}$, for $n\geq1$, satisfies the uniform exterior sphere
condition in time with radius $r_{0}$, independent of $n$, and that there
exist $0<\rho_{0}<r_{0}$ and $\eta_{0}>0$ such that, for all $n\geq1$,
%
\begin{eqnarray}\label{crita+}
\inf_{s\in [ 0,T ] }\inf_{z\in\partial
D_{s}^{n}}a_{s,z}^{n} ( \rho_{0},\eta_{0} ) &=&a_{n}>0,
\\
\label{crite+}
\sup_{s\in [ 0,T ] }\sup_{z\in\partial
D_{s}^{n}}e_{s,z}^{n} ( \rho_{0},\eta_{0} ) &=&e_{n}<1.
\end{eqnarray}
Furthermore, we let
%
\begin{equation}
l_{n} ( r ) =\mathop{\mathop{\sup}_{ s,t\in[ 0,T]}}_{ \vert
s-t \vert\leq r}\sup_{z\in\overline{D_{s}^{n}}}d (
z,D_{t}^{n} ) , \label{defl+}
\end{equation}
and we assume that
%
\begin{equation}
\lim_{r\rightarrow0^{+}}\sup_{n\geq1}l_{n} ( r ) =0.
\label{limitzero+}
\end{equation}
Note also that if we define $\hat{l}_{n} ( r ) $ as in (\ref{lhat})
but with $D$ replaced by $D^{n}$, then Lemma~\ref{llhatequivalence}
and %
\eqref{limitzero+} imply that
%
\begin{equation}
\lim_{r\rightarrow0^{+}}\sup_{n\geq1}\hat{l}_{n} ( r ) =0.
\label{limitzero+aa}
\end{equation}
Moreover, we assume that there exists $\hat{R}>0$ such that $%
D_{t}^{n}\subset B(0,\hat{R})$ and $D_{t}\subset B(0,\hat{R})$, for
all $%
n\geq1$ and $t\in\lbrack0,T]$, and we let
%
\begin{equation}
R=2\sup_{t\in [ 0,T ] }\sup_{n}\max\{\operatorname
{diam}(D_{t}^{n}),\operatorname{diam}
( D_{t} ) \}, \label{defR}
\end{equation}
where diam$(D_{t}^{n})$, diam$(D_{t})$ are the Euclidean diameters of the
spatial regions $D_{t}^{n}$ and $D_{t}$, respectively. Recalling that the
set $G^{\Gamma}$ was introduced in \eqref{cone2}, we here also introduce,
for $n\geq1$,
%
\begin{equation}
G^{\Gamma^{n}}= \{ ( s,z,\gamma ) \dvtx \gamma\in\Gamma
_{t}^{n} ( z ) , z\in\partial D_{s}^{n}, s\in [
0,T ] \} , \label{cone2+}
\end{equation}
and we let, whenever $t\in\lbrack0,T]$,
%
\begin{eqnarray} \label{Gtdef}
G_{t} &=&G^{\Gamma}\cap \bigl( [ 0,t ] \times{B_{R} (
0 ) }\times S_{1} ( 0 ) \bigr) , \nonumber
\\[-8pt]
\\[-8pt]
G_{t}^{n} &=&G^{\Gamma^{n}}\cap \bigl( [ 0,t ] \times{%
B_{R} ( 0 ) }\times S_{1} ( 0 ) \bigr) .\nonumber
\end{eqnarray}
In the following we need to measure the distance between the sets $G_{T}$
and $G_{T}^{n}$ and, hence, we introduce an appropriate Hausdorff
distance for
subsets of $[0,T]\times{B_{R} ( 0 ) }\times S_{1} ( 0 ) $%
. In particular, we let, given $(s,z,\gamma)\in\lbrack0,T]\times{%
B_{R} ( 0 ) }\times S_{1} ( 0 ) $ and $(\hat{s},\hat{z},%
\hat{\gamma})\in\lbrack0,T]\times{B_{R} ( 0 ) }\times
S_{1} ( 0 ) $,
%
\begin{equation}
E((s,z,\gamma),(\hat{s},\hat{z},\hat{\gamma}))=|s-\hat
{s}|+|z-\hat{z}%
|+|\gamma-\hat{\gamma}| \label{Gdefhaus-}
\end{equation}
denote the (Euclidean) distance between $(s,z,\gamma)$ and $(\hat
{s},\hat{z}%
,\hat{\gamma})$. Furthermore, based on $E$, we define, given $%
F_{1},F_{2}\subseteq\lbrack0,T]\times{B_{R} ( 0 ) }\times
S_{1} ( 0 ) $ and $(s,z,\gamma)\in\lbrack0,T]\times{B_{R} (
0 ) }\times S_{1} ( 0 ) $, the distances $E((s,z,\gamma
),F_{1})$, $E((s,z,\gamma),F_{2})$ and\break $E(F_{1},F_{2})$ in the natural way.
Furthermore, for $F_{1}$ and $F_{2}$ as above, we introduce a Hausdorff
distance between $F_{1}$ and $F_{2}$ as
%
\begin{eqnarray}\label{Gdefhaus}
H(F_{1},F_{2}) &=&\max\{A,B\}, \nonumber \\
A &=&\sup\{E((s,z,\gamma),F_{2})\dvtx  (s,z,\gamma)\in F_{1}\},
\\
B &=&\sup\{E((\hat{s},\hat{z},\hat{\gamma}),F_{1})\dvtx  (\hat{s},\hat
{z},\hat{%
\gamma})\in F_{2}\}.\nonumber
\end{eqnarray}
In the following we say that $G_{T}^{n}$ converges to $G_{T}$ if
%
\begin{equation}
H(G_{T}^{n},G_{T})\rightarrow0 \qquad\mbox{as }n\rightarrow\infty.
\label{Gdefhaus+}
\end{equation}
Imposing the assumptions on $D$, $\Gamma$ stated above and assuming %
\eqref{Gdefhaus+}, we can, for example, ensure that if $ \{
(s_{n},z_{n}) \} $ is a sequence of points in $%
\mathbb{R}
^{d+1}$, $s_{n}\in\lbrack0,T]$, $z_{n}\in\partial D_{s_{n}}^{n}$, $%
\lim_{n\rightarrow\infty}s_{n}=s\in\lbrack0,T]$, $\lim
_{n\rightarrow
\infty}z_{n}=z\in\partial D_{s}$, then
%
\begin{equation}
\lim_{n\rightarrow\infty}h(\Gamma_{s_{n}}^{n} ( z_{n} ) ,\Gamma
_{s} ( z ) )=0. \label{Gamman}
\end{equation}
To see this, we consider, for $ \{ (s_{n},z_{n}) \} $ and $(s,z)$
given as above, $(s_{n},z_{n},\gamma_{s_{n}}^{n})\in G_{T}^{n}$ and $%
(s,z,\gamma_{s})\in G_{T}$. Given $(s_{n},z_{n},\gamma
_{s_{n}}^{n})\in
G_{T}^{n}$, we let $(\hat{s}_{n},\hat{z}_{n},\hat{\gamma}_{\hat{s}%
_{n}}^{n})\in G_{T}$ be a point on $G_{T}$ which minimizes the
distance, as
defined in \eqref{Gdefhaus-}, from $(s_{n},z_{n},\gamma_{s_{n}}^{n})$
to $%
G_{T}$. Then,
%
\begin{eqnarray}\label{Gamman++}
|\gamma_{s_{n}}^{n}-\gamma_{s}| &\leq&E((s_{n},z_{n},\gamma
_{s_{n}}^{n}),(s,z,\gamma_{s})) \nonumber\\
&\leq&E((s_{n},z_{n},\gamma_{s_{n}}^{n}),(\hat{s}_{n},\hat
{z}_{n},\hat{%
\gamma}_{\hat{s}_{n}}^{n}))+E((\hat{s}_{n},\hat{z}_{n},\hat{\gamma
}_{\hat{s}%
_{n}}^{n}),(s,z,\gamma_{s})) \\
&\leq&H(G_{T},G_{T}^{n})+E((\hat{s}_{n},\hat{z}_{n},\hat{\gamma
}_{\hat{s}%
_{n}}^{n}),(s,z,\gamma_{s})).\nonumber
\end{eqnarray}
Hence,
%
\begin{equation}
h(\Gamma_{s_{n}}^{n} ( z_{n} ) ,\Gamma_{s} ( z ) )\leq
H(G_{T},G_{T}^{n})+R_{n},
\end{equation}
where
%
\begin{eqnarray} \label{Gamman++bbb}
R_{n} &=&\max\{A_{n},B_{n}\}, \nonumber \\
A_{n} &=&\sup\{E((\hat{s}_{n},\hat{z}_{n},\hat{\gamma}{}_{\hat{s}%
_{n}}^{n}),\{(s,z,\Gamma_{s}(z))\})\dvtx  \hat{\gamma}{}_{\hat
{s}_{n}}^{n}\in
\Gamma{}_{\hat{s}_{n}}^{n}(\hat{z}_{n})\}, \\
B_{n} &=&\sup\{E(\{(\hat{s}_{n},\hat{z}_{n},\Gamma_{\hat
{s}_{n}}^{n}(\hat{z%
}_{n}))\},(s,z,\gamma_{s}))\dvtx  \gamma_{s}\in\Gamma_{s}(z)\}.\nonumber
\end{eqnarray}
As, by assumption, $G_{T}$ is closed, we can now first conclude that $%
R_{n}\rightarrow0$ as $n\rightarrow\infty$, and then we find, using %
\eqref{Gdefhaus+}, that $h(\Gamma_{s_{n}}^{n} ( z_{n} ) ,\Gamma
_{s} ( z ) )\rightarrow0$ as $n\rightarrow\infty$. This
completes the proof of \eqref{Gamman}. We are now ready to formulate our
convergence result.

\begin{theorem}
\label{SPconvergence} Let $T>0$ and let $D\subset
\mathbb{R}
^{d+1}$ be a time-dependent domain satisfying \eqref{timedep+}. Let
$\Gamma
=\Gamma_{t} ( z ) $ be a closed convex cone of vectors in $%
\mathbb{R}
^{d}$ for every $z\in\partial D_{t}$, $t\in [ 0,T ] $, and assume
that $\Gamma$ satisfies \eqref{cone1} and \eqref{cone2}. Let $%
\{D^{n}\}_{n=1}^{\infty}$ be a sequence of time-dependent domains $%
D^{n}\subset
\mathbb{R}
^{d+1}$ satisfying \eqref{timedep+} and a uniform exterior sphere condition
in time with radius $r_{0}$ in the sense of \eqref{extsphere-}.\vspace*{1pt} Let $%
\{\Gamma^{n}\}_{n=1}^{\infty}=\{\Gamma_{t}^{n} ( z )
\}_{n=1}^{\infty}$ be a sequence of closed convex cones $\Gamma
^{n}=\Gamma
_{t}^{n} ( z ) $ of vectors in $%
\mathbb{R}
^{d}$ for every $z\in\partial D_{t}^{n}$, $t\in\lbrack0,T]$. For all
$%
n\geq1$, $D^{n}$ and $\Gamma^{n}$ satisfy \eqref{crita+} and \eqref{crite+}
for some $0<\rho_{0}<r_{0}$, $\eta_{0}>0$, $a_{n}$, $e_{n}$ and, moreover,
$((0,T)\times
\mathbb{R}
^{d})\setminus\overline{D^{n}}$ has the $(\delta
_{0},h_{0})$-property of
good projections along $\Gamma^{n}$, for some $0<\delta_{0}<\rho
_{0}$, $%
h_{0}>1$. Assume that $\inf_{n}\{a_{n}\}>0$, $\sup_{n}\{e_{n}\}<1$
and %
\eqref{limitzero+} hold. Regarding the convergence $D^{n}\rightarrow
D$ and $%
\Gamma^{n}\rightarrow\Gamma$, assume that
%
\begin{eqnarray} \label{Dn}
\lim_{n\rightarrow\infty}\sup_{t\in [ 0,T ] }h ( {D_{t}^{n}},%
{D_{t}} ) &=&0, \\ \label{dDn}
\lim_{n\rightarrow\infty}\sup_{t\in [ 0,T ] }h ( \partial{%
D_{t}^{n}},\partial{D_{t}} ) &=&0,
\end{eqnarray}
and, with $G_{T}$ and $G_{T}^{n}$ defined as in \eqref{Gtdef}, that%
%
\begin{equation}
G_{T}^{n}\mbox{ converges to }G_{T}\mbox{ in the sense of \eqref{Gdefhaus+}.}
\label{Gamman+}
\end{equation}
Let $w^{n}\in\mathcal{D} ( [ 0,T ] ,%
\mathbb{R}
^{d} ) $ with $w_{0}^{n}\in\overline{D_{0}^{n}}$ and assume that there
exists a solution $ ( x^{n},\lambda^{n} ) $ to the Skorohod
problem for $ ( D^{n},\Gamma^{n},w^{n} ) $ such that $x^{n}\in
\mathcal{D}^{\rho_{0}} ( [ 0,T ] ,%
\mathbb{R}
^{d} ) $ for all $n\geq1$. Assume that $ \{ w^{n} \} $ is
relatively compact in the Skorohod topology and that $ \{ w^{n} \} $
converges to $w\in\mathcal{D} ( [ 0,T ] ,%
\mathbb{R}
^{d} ) $. Then $ \{ ( x^{n},\lambda^{n} ) \} $
converges to $(x,\lambda)\in\mathcal{D} ( [ 0,T ] ,%
\mathbb{R}
^{d} ) \times\mathcal{BV} ( [ 0,T ] ,%
\mathbb{R}
^{d} ) $ and $(x,\lambda)$ is a solution to the Skorohod problem for $%
( D,\Gamma,w ) $ with $x\in\mathcal{D}^{\rho_{0}} ( [
0,T ] ,%
\mathbb{R}
^{d} ) $.
\end{theorem}

\begin{remark}
\label{subtel}We note that the formulation of Theorem \ref{SPconvergence}
contains several subtle points. First, we do not have to assume that the
elements in the sequence $\{\Gamma^{n}\}_{n=1}^{\infty}$ satisfy %
\eqref{cone1}, \eqref{cone2} and \eqref{Gammanpre}. The reason for
this (see
Remark \ref{conenotres}) is that Theorem \ref{compacttheorem} holds, with
constants independent of $n$, for each element in the sequence $ \{
( w^{n},x^{n},\lambda^{n} ) \} $ even without these
assumptions. Second, we only have to impose very modest restrictions on $D$
but, as can be seen in the proof below, we have to assume that $\Gamma
=\Gamma_{t} ( z ) $ satisfies \eqref{cone1} and \eqref{cone2}.
\end{remark}

\begin{pf*}{Proof of Theorem \ref{SPconvergence}}
As $ \{
w^{n} \} $ is relatively compact in the Skorohod topology, we first
note that Theorem \ref{compacttheorem} can be used to conclude the existence
of positive constants $L_{1}$, $L_{2}$, $L_{3}$ and $L_{4}$,
independent of $%
n$, such that%
%
\begin{eqnarray} \label{gg1}
\Vert x^{n} \Vert_{t_{1},t_{2}} &\leq&L_{1} \Vert
w^{n} \Vert_{t_{1},t_{2}}+L_{2}l_{n} ( t_{2}-t_{1} ) , \nonumber
\\[-8pt]
\\[-8pt]
\vert\lambda^{n} \vert_{t_{2}}- \vert\lambda
^{n} \vert_{t_{1}} &\leq&L_{3} \Vert w^{n} \Vert
_{t_{1},t_{2}}+L_{4}l_{n} ( t_{2}-t_{1} ) ,\nonumber
\end{eqnarray}
whenever $0\leq t_{1}\leq t_{2}\leq T$. As $ \{ w^{n} \} $
converges to $w\in\mathcal{D} ( [ 0,T ] ,%
\mathbb{R}
^{d} ) $, we also see, using \eqref{limitzero+} and \eqref{gg1}, that
$%
\{ ( w^{n},x^{n},\lambda^{n},|\lambda^{n}| ) \} $ is
relatively compact in $\mathcal{D} ( [ 0,T ] ,%
\mathbb{R}
^{d} ) \times\mathcal{D} ( [ 0,T ] ,%
\mathbb{R}
^{d} ) \times\mathcal{D} ( [ 0,T ] ,%
\mathbb{R}
^{d} ) \times\mathcal{D} ( [ 0,T ] ,%
\mathbb{R}
_{+} ) $. Furthermore, we know that $x_{t}^{n}\in\overline{D^{n}}$ for
all $t\in\lbrack0,T]$, $n\geq1$. Hence, $ \{ ( x^{n},\lambda
^{n} ) \} $ converges to some $(x,\lambda)\in\mathcal{D} ( %
[ 0,T ] ,%
\mathbb{R}
^{d} ) \times\mathcal{D} ( [ 0,T ] ,%
\mathbb{R}
^{d} ) $. We intend to prove that $ ( x,\lambda ) \in
\mathcal{D} ( [ 0,T ] ,%
\mathbb{R}
^{d} ) \times\mathcal{BV} ( [ 0,T ] ,%
\mathbb{R}
^{d} ) $ solves the Skorohod problem for $ ( D,\Gamma,w ) $
and to do this, we have to prove that
%
\begin{equation}
\lambda\in\mathcal{BV} ( [ 0,T ] ,%
\mathbb{R}
^{d} ) , \label{BV}
\end{equation}
and we have to verify that
%
\begin{equation}\qquad
( D,\Gamma,w ) \mbox{ and } ( x,\lambda ) \mbox{
satisfy properties (\ref{SP1})--(\ref{SP4}) in Definition \ref{skorohodprob}.}
\label{prop}
\end{equation}
We begin by verifying (\ref{SP1}) in Definition \ref{skorohodprob}.
To do
this, we first note, using the convergence properties of the Skorohod
topology and the fact that $ ( x^{n},\lambda^{n} ) $ solves the
Skorohod problem for $ ( D^{n},\Gamma^{n},w^{n} ) $, that
%
\begin{equation}
x_{t}=w_{t}+\lambda_{t}  \label{hra1}
\end{equation}
for all points of continuity and, hence, since $w$, $x$ and $\lambda$
are c%
\`{a}dl\`{a}g functions, that \eqref{hra1} holds for all $t\in [ 0,T%
] $. Hence, to verify (\ref{SP1}) in Definition \ref{skorohodprob}, it
only remains to ensure that $x_{t}\in\overline{D_{t}}$ for all $t\in [
0,T ] $. To do this, we first note, using Proposition 5.3 and Remark 5.4
in Chapter 3 of \cite{EthierKurtz1986}, that there exists a sequence $%
\{ \tilde{t}_{n} \} $ such that%
%
\begin{eqnarray}\label{convTh3}
\tilde{t}_{n}&\rightarrow& t,\qquad x_{\tilde{t}_{n}}^{n}\rightarrow x_{t},\qquad
x_{%
\tilde{t}_{n}^{-}}^{n}\rightarrow x_{t^{-}},\nonumber
\\[-8pt]
\\[-8pt]
 \lambda_{\tilde{t}%
_{n}}^{n}&\rightarrow&\lambda_{t},\qquad \lambda_{\tilde{t}%
_{n}^{-}}^{n}\rightarrow\lambda_{t^{-}}\qquad\mbox{as }n\rightarrow
\infty.\nonumber
\end{eqnarray}
Furthermore, using the triangle inequality, \eqref{limitzero}, \eqref{Dn}
and \eqref{convTh3}, we obtain
%
\begin{eqnarray}\label{bbbb1-}\qquad
d ( x_{t},\overline{D_{t}} ) &\leq& \vert x_{t}-x_{\tilde{t}%
_{n}}^{n} \vert+d ( x_{\tilde{t}_{n}}^{n},\overline{D_{\tilde{t}%
_{n}}^{n}} ) +h ( \overline{D_{\tilde{t}_{n}}^{n}},\overline{D_{%
\tilde{t}_{n}}} ) +h ( \overline{D_{\tilde{t}_{n}}},\overline{D_{t}}
) \nonumber
\\[-8pt]
\\[-8pt]\qquad
&\leq& \vert x_{t}-x_{\tilde{t}_{n}}^{n} \vert+h ( \overline{%
D_{\tilde{t}_{n}}^{n}},\overline{D_{\tilde{t}_{n}}} ) +l(|\tilde{t}%
_{n}-t|)\rightarrow0, \qquad\mbox{as }n\rightarrow\infty.\nonumber
\end{eqnarray}
This proves that $x_{t}\in\overline{D_{t}}$ for all $t\in [ 0,T ]
$ and, hence, we have verified that $(w_{t},x_{t},\lambda_{t})$
satisfies (%
\ref{SP1}). We next prove \eqref{BV}, that is, that $\lambda\in
\mathcal{BV}%
( [ 0,T ] ,%
\mathbb{R}
^{d} ) $. To do this, we use an argument similar to the proof of Theorem
3.1 in \cite{Costantini1992}, but, as described below, our argument is more
subtle due to the fact that we consider sequences $(D^{n},\Gamma^{n},w^{n})$
where, in particular, $D^{n}$ is time-dependent. Recall that with~$R$ as
introduced in \eqref{defR}, we have
%
\begin{equation}
\sup_{n}\sup_{t\in [ 0,T ] } \vert x_{t}^{n} \vert
<R, \qquad\sup_{t\in [ 0,T ] } \vert x_{t} \vert<R.
\label{defR+}
\end{equation}
Let $G^{\Gamma^{n}}$, $G_{t}$ and $G_{t}^{n}$ be as in \eqref{cone2+}
and %
\eqref{Gtdef}. By the prerequisites of Theorem \ref{SPconvergence} (see
Remark \ref{subtel}), we have that $G^{\Gamma}$ is closed. We next
define a
positive measure~$\mu^{n}$ on $ [ 0,T ] \times{B_{R} (
0 ) }\times S_{1} ( 0 ) $ by setting, for every Borel set $%
A\subset [ 0,T ] \times{B_{R} ( 0 ) }\times S_{1} (
0 ) $,%
%
\begin{equation}
\mu^{n} ( A ) =\int_{0}^{T}\chi_{A\cap G_{T}^{n}} (
s,x_{s}^{n},\gamma_{s}^{n} ) \,d \vert\lambda^{n} \vert
_{s}, \label{mundef}
\end{equation}
where $\gamma_{s}^{n}\in\Gamma_{s}^{n,1}(x_{s}^{n})$ is as in (\ref
{SP2}%
)--(\ref{SP4}) for the solution $ ( x^{n},\lambda^{n} ) $ to the
Skorohod problem for $ ( D^{n},\Gamma^{n},w^{n} ) $ and $\chi
_{A\cap G_{T}^{n}}$ is the characteristic functions for the set $A\cap
G_{T}^{n}$. We then first note that%
%
\begin{equation}
\vert\lambda^{n} \vert_{t}=\mu^{n} ( G_{t}^{n} )
 \qquad\mbox{whenever }t\in [ 0,T ] . \label{4.22}
\end{equation}
We also note that the support of $\mu^{n}$ is contained in $G_{T}^{n}$ in
the sense that $\mu^{n}(A)=0$ whenever $A\subset [ 0,T ] \times{%
B_{R} ( 0 ) }\times S_{1} ( 0 ) $ is such that $A\cap
G_{T}^{n}=\emptyset$. Using this, and the fact that (\ref{SP2}) holds
for $%
\lambda^{n}$, we see that \eqref{4.22} implies that
%
\begin{equation}
\lambda_{t}^{n}=\int_{ [ 0,t ] \times{B_{R} ( 0 ) }%
\times S_{1} ( 0 ) }\gamma\,d\mu^{n} ( s,z,\gamma )
 \qquad\mbox{whenever }t\in [ 0,T ] . \label{lambdamurel}
\end{equation}
Next, using \eqref{limitzero+}, \eqref{gg1}, \eqref{4.22} and the
fact that $%
\{ w^{n} \} $ converges to $w\in\mathcal{D} ( [ 0,T%
] ,%
\mathbb{R}
^{d} ) $, we conclude that
%
\begin{equation}
\sup_{n}\mu^{n} ( G_{T}^{n} ) <\infty, \label{laba1}
\end{equation}
which implies that $ \{ \mu^{n} \} $ is a compact set of measures,
on $ [ 0,T ] \times{B_{R} ( 0 ) }\times S_{1} (
0 ) $, in the sense of the weak$\ast$-topology. Therefore, by the
Banach--Alaoglu theorem, we can conclude that $ \{ \mu^{n} \} $
converges in the weak$\ast$-topology to a measure $\mu$ such that
%
\begin{equation}
\mu\bigl( [ 0,T ] \times{B_{R} ( 0 ) }\times S_{1} (
0 ) \bigr)<\infty. \label{boundmeas}
\end{equation}
Moreover, since $ ( x^{n},\lambda^{n} ) $ converges to $ (
x,\lambda ) $ in the sense of the Skorohod topology, we obtain, using (\ref{lambdamurel}), that
%
\begin{equation}
\lambda_{t}=\int_{ [ 0,t ] \times{B_{R} ( 0 ) }\times
S_{1} ( 0 ) }\gamma\,d\mu ( s,z,\gamma )
\label{lambdamurel2}
\end{equation}
for all $t\in [ 0,T ] $ such that $\lambda_{t}=\lambda_{t^{-}}$.
However, as both sides of (\ref{lambdamurel2}) are right continuous,
(\ref%
{lambdamurel2}) holds for all $t\in [ 0,T ] $. Having proved %
\eqref{lambdamurel2}, we see, also using~\eqref{boundmeas}, that
$\lambda$
is of bounded variation and, hence, \eqref{BV} is proved. We next
claim that
%
\begin{equation}
\lambda_{t}=\int_{G_{t}}\gamma\,d\mu ( s,z,\gamma ) ,
\label{lambdamurel2+}
\end{equation}
that is, we claim that the support of the measure $\mu$ is the set $G_{T}$
in the sense that if $A\subset [ 0,T ] \times{B_{R} (
0 ) }\times S_{1} ( 0 ) $ is such that $A\cap G_{T}=\emptyset$%
, then $\mu(A)=0$. To see this, we let $(\hat{s},\hat{z},\hat
{\gamma})\in
\lbrack0,T]\times{B_{R} ( 0 ) }\times S_{1} ( 0 )
\setminus G_{T}$ and we see, as $G_{T}$ is closed, that if we define,
for $%
\eta>0$, $B((\hat{s},\hat{z},\hat{\gamma}),\eta):=\{(s,z,\gamma
)\dvtx E((s,z,\gamma),(\hat{s},\hat{z},\hat{\gamma}))<\eta\}\cap
\lbrack
0,T]\times{B_{R} ( 0 ) }\times S_{1} ( 0 ) $, then there
exists $\eta_{0}>0$ such that $B((\hat{s},\hat{z},\hat{\gamma
}),2\eta
_{0})\cap G_{T}=\emptyset$. Recall that $E$ is the distance function
introduced in \eqref{Gdefhaus-}. Furthermore, the above setup implies
that $%
E(B((\hat{s},\hat{z},\hat{\gamma}),\eta_{0}),G_{T})>\eta_{0}$ and since
%
\begin{equation}
E(B((\hat{s},\hat{z},\hat{\gamma}),\eta_{0}),G_{T}^{n})\geq
E(B((\hat{s},%
\hat{z},\hat{\gamma}),\eta_{0}),G_{T})-H(G_{T},G_{T}^{n}),
\end{equation}
we can use the assumption in \eqref{Gamman+} to conclude that there
exists $%
n_{0}\in
\mathbb{N}
$ such that
%
\begin{equation}
E(B((\hat{s},\hat{z},\hat{\gamma}),\eta_{0}),G_{T}^{n})\geq\eta_{0}/2
\end{equation}
for all $n\geq n_{0}$. In particular, $B((\hat{s},\hat{z},\hat
{\gamma}),\eta
_{0})\cap G_{T}^{n}=\emptyset$ for all $n\geq n_{0}$. Hence, $\mu
^{n}(B((%
\hat{s},\hat{z},\hat{\gamma}),\eta_{0}))=0$, for all $n\geq
n_{0}$, and $%
\mu(B((\hat{s},\hat{z},\hat{\gamma}),\eta_{0}))=0$ by the
weak$\ast$%
-convergence of $\mu^{n}$ to $\mu$. This completes the proof of %
\eqref{lambdamurel2+}. Having proved \eqref{BV} and \eqref
{lambdamurel2+}, we
see that
%
\begin{equation}
\lambda_{t}=\int_{0}^{t}\gamma_{s}\,d \vert\lambda \vert_{s}
\qquad\mbox{whenever }t\in [ 0,T ]   \label{laba2}
\end{equation}
for some $S_{1} ( 0 ) $-valued Borel measurable function $\gamma
_{s}$ and, to prove (\ref{SP2}), we have to prove that $\gamma_{s}\in
\Gamma
_{s}^{1}(x_{s})$ for all $s\in\lbrack0,T]$. To prove this and to
verify (%
\ref{SP4}), we consider the following two cases:
%
\begin{eqnarray} \label{cases1}
\mbox{Case 1.}\quad t\in [ 0,T ] \mbox{ is such that }%
x_{t}-x_{t^{-}}&\neq&0, \nonumber
\\[-8pt]
\\[-8pt]
\mbox{Case 2.}\quad t\in [ 0,T ] \mbox{ is such that }%
x_{t}-x_{t^{-}}&=&0.\nonumber
\end{eqnarray}

%

\textit{Case} 1. Note that Case 1 occurs for an at most countable
set of jump times of~$x$. Moreover, in Case 1 it is enough to prove that
%
\begin{equation}
\lambda_{t}-\lambda_{t^{-}}\neq0\mbox{ implies that }x_{t}\in
\partial
D_{t}\mbox{ and that }\lambda_{t}-\lambda_{t^{-}}\in\Gamma_{t} (
x_{t} ) . \label{hra2}
\end{equation}
We first note, as we are assuming $\lambda_{t}-\lambda_{t^{-}}\neq0$,
that $ \vert\lambda_{\tilde{t}_{n}}^{n}-\lambda_{\tilde{t}%
_{n}^{-}}^{n} \vert>0$ for $n$ sufficiently large; see \eqref{convTh3}
. Furthermore, since $ ( x^{n},\lambda^{n} ) $ solves the Skorohod
problem for $ ( D^{n},\Gamma^{n},w^{n} ) $, we have that%
%
\begin{equation}
x_{\tilde{t}_{n}}^{n}\in\partial D_{\tilde{t}_{n}}^{n},\qquad {\lambda_{%
\tilde{t}_{n}}^{n}-\lambda_{\tilde{t}_{n}^{-}}^{n}}\in\Gamma
_{\tilde{t}%
_{n}}^{n} ( x_{\tilde{t}_{n}}^{n} ) . \label{convTh2}
\end{equation}
Combining \eqref{limitzero+aa}, \eqref{dDn}, \eqref{convTh3} and
(\ref%
{convTh2}), we obtain%
%
\begin{eqnarray}\label{bbb1'}\qquad
d ( x_{t},\partial D_{t} ) &\leq&|x_{t}-x_{\tilde{t}%
_{n}}^{n}|+d(x_{\tilde{t}_{n}}^{n},\partial D_{\tilde{t}_{n}}^{n})+h(%
\partial D_{\tilde{t}_{n}}^{n},\partial D_{\tilde{t}_{n}})+h(\partial
D_{%
\tilde{t}_{n}},\partial D_{t}) \nonumber
\\[-8pt]
\\[-8pt]\qquad
&\leq&|x_{t}-x_{\tilde{t}_{n}}^{n}|+h(\partial D_{\tilde{t}%
_{n}}^{n},\partial D_{\tilde{t}_{n}})+\hat{l}(|\tilde
{t}_{n}-t|)\rightarrow0%
 \qquad\mbox{as }n\rightarrow\infty.\nonumber
\end{eqnarray}
Hence, using \eqref{bbb1'}, we can, since $\partial D_{t}$ is closed,
conclude that
%
\begin{equation}
x_{t}\in\partial D_{t}. \label{existthprop1}
\end{equation}
We next recall that the set $G^{\Gamma}$, defined in \eqref{cone2},
is, by
assumption, closed. Furthermore, arguing as in \eqref{bbb1'}, we first see
that
%
\begin{eqnarray}\label{bbbb2}
d ( \lambda_{t}-\lambda_{t^{-}},\Gamma_{t}(x_{t}) ) &\leq
&|(\lambda_{t}-\lambda_{t^{-}})-(\lambda_{\tilde
{t}_{n}}^{n}-\lambda_{%
\tilde{t}_{n}^{-}}^{n})|+d(\lambda_{\tilde{t}_{n}}^{n}-\lambda
_{\tilde{t}%
_{n}^{-}}^{n},\Gamma_{\tilde{t}_{n}}^{n}(x_{\tilde{t}_{n}}^{n}))
\nonumber\\
&&{}+h(\Gamma_{\tilde{t}_{n}}^{n}(x_{\tilde{t}_{n}}^{n}),\Gamma_{t}(x_{t}))
\\
&\leq&|\lambda_{t}-\lambda_{\tilde{t}_{n}}^{n}|+|\lambda
_{t^{-}}-\lambda
_{\tilde{t}_{n}^{-}}^{n}|+h(\Gamma_{\tilde{t}_{n}}^{n}(x_{\tilde{t}%
_{n}}^{n}),\Gamma_{t}(x_{t})),\nonumber
\end{eqnarray}
and then, letting $n\rightarrow\infty$, it follows, using \eqref{Gamman},
that $d ( \lambda_{t}-\lambda_{t^{-}},\Gamma_{t}(x_{t}) ) =0$.
Applying the fact that the set $G^{\Gamma}$ is closed, we can therefore
conclude that
%
\begin{equation}
\lambda_{t}-\lambda_{t^{-}}\in\Gamma_{t}(x_{t}).
\end{equation}
This concludes the proof of \eqref{hra2} and, hence, we have verified
(\ref%
{SP2})--(\ref{SP4}) in Case~1.

%

\textit{Case} 2. To verify (\ref{SP2})--(\ref{SP4}) in Case 2, we
first see, by combining (\ref{lambdamurel2+}) and \eqref{laba2}, that%
%
\begin{equation}
\int_{0}^{t}\gamma_{s}\,d \vert\lambda \vert
_{s}=\int_{G_{t}}\gamma\,d\mu ( s,z,\gamma )  \qquad\mbox{whenever }%
t\in [ 0,T ] . \label{lambdamurel3}
\end{equation}
We next introduce a measure $\nu$ on $[0,T]$ by setting
%
\begin{equation}
\nu([0,t])=\mu ( G_{t} )  \qquad\mbox{whenever }t\in [ 0,T ] .
\label{nudef}
\end{equation}
Combining (\ref{lambdamurel3}) and (\ref{nudef}), it is clear that
$\nu
([0,t])=0$ implies $ \vert\lambda \vert_{t}=0$, showing that $%
\vert\lambda \vert$ is absolutely continuous with respect to $%
\nu$. To simplify the notation, in the following we let, for $k\in
\mathbb{%
\mathbb{N}
}$,
%
\begin{equation} \qquad
\Omega_{k}= \biggl\{ ( t,z,\gamma ) \in G_{T}\dvtx \inf_{s\in [ 0,T%
] } \bigl( \vert t-s \vert+ ( \vert
z-x_{s} \vert\wedge \vert z-x_{s^{-}} \vert ) \bigr) >%
\frac{1}{k} \biggr\} .
\end{equation}
Then, using Theorem 1.2.1(iii) in \cite{Friedman1982}, the fact that
$\mu
( U ) \leq$\underline{$\lim$}$_{n\rightarrow\infty}\mu
^{n} ( U ) $ for all open sets $U$ and the fact that $x_{t}^{n}$
converges either to $x_{t}$ or $x_{t^{-}}$, we can conclude that
%
\begin{eqnarray}\label{muzero}
&&\mu \bigl( \{ ( t,z,\gamma ) \in G_{T}\dvtx  z\neq x_{t},z\neq
x_{t^{-}} \} \bigr)\nonumber\\
&&\qquad =
\lim_{k\rightarrow\infty}\mu ( \Omega
_{k} ) \leq\lim_{k\rightarrow\infty}\operatorname{\underline {\lim}}\limits_{n\rightarrow
\infty}\mu^{n} ( \Omega_{k} ) \\
&&\qquad\leq\lim_{k\rightarrow\infty}\operatorname{\underline {\lim}}\limits_{n\rightarrow
\infty}
\mu^{n} \biggl( \biggl\{ ( t,z,\gamma ) \in G_{T}\dvtx  \vert
z-x_{t}^{n} \vert>\frac{1}{2k} \biggr\} \biggr) =0.\nonumber
\end{eqnarray}
If $x_{t}=x_{t^{-}}\in D_{t}$, then, since $z\in\partial D_{t}$ for
all $%
( t,z,\gamma ) \in G_{T}$, we deduce that $z\neq x_{t}$ and $%
z\neq x_{t^{-}}$. Hence, using (\ref{nudef}) and (\ref{muzero}), we
first see
that
%
\begin{equation}
\nu ( \{ t\in [ 0,T ] \dvtx x_{t}=x_{t^{-}},x_{t}\in
D_{t} \} ) =0,
\end{equation}
and then, by the absolute continuity of $ \vert\lambda \vert$
with respect to $\nu$, we can conclude that
%
\begin{equation}
\vert\lambda \vert ( \{ t\in [ 0,T ]
\dvtx x_{t}=x_{t^{-}},x_{t}\in D_{t} \} ) =0. \label{milan1}
\end{equation}
In particular, \eqref{milan1} proves (\ref{SP4}). Hence, it only
remains to
prove that $\gamma_{s}$, as defined in \eqref{laba2}, satisfies
$\gamma
_{s}\in\Gamma_{s}^{1}(x_{s})$ for all $s\in\lbrack0,t]$ and $t\in [
0,T ] $. From (\ref{lambdamurel3}) and the fact that (\ref{muzero})
implies that%
%
\begin{equation}
\mu \bigl( \{ ( s,z,\gamma ) \in G_{t}\dvtx x_{s}=x_{s^{-}}\neq
z \} \bigr) =0,
\end{equation}
we deduce that
%
\begin{eqnarray} \label{ddd}
\int_{ \{ s\in [ 0,t ] \dvtx  x_{s}=x_{s^{-}} \} }\gamma
_{s}\,d \vert\lambda \vert_{s}
&=&\int_{ \{ (
s,z,\gamma ) \in G_{t}\dvtx z=x_{s}=x_{s^{-}} \} }\gamma\,d\mu (
s,z,\gamma ) \nonumber\\
&=&\int_{ \{ s\in [ 0,t ] \dvtx x_{s}=x_{s^{-}} \}
}\int_{\Gamma_{s}^{1} ( x_{s} ) }\gamma p ( s,x_{s},d\gamma
) \,d\nu,
\end{eqnarray}
whenever $t\in [ 0,T ] $. Note that the last equality in %
\eqref{ddd} follows from the definition of $\nu$ in (\ref{nudef}).
Here $%
p ( s,x_{s},\cdot ) $ is a measure on the Borel $\sigma$-algebra
of $S_{1} ( 0 ) $, concentrated on $\Gamma_{s}^{1} (
x_{s} ) $ for $d\nu$-almost all $s\in [ 0,T ] $ such that $%
x_{s}=x_{s^{-}}$, and $p ( \cdot,\cdot,A ) $ is a nonnegative
Borel measurable function for every Borel set $A$. Then, since $ \vert
\lambda \vert$ is absolutely continuous with respect to $\nu$, the
Radon--Nikodym theorem asserts the existence of a nonnegative Borel
measurable function $f$ such that%
%
\begin{eqnarray}\label{ddd5}
f ( s ) \gamma_{s}=\int_{\Gamma_{s}^{1} ( x_{s} )
}\gamma p ( s,x_{s},d\gamma ) ,
\nonumber
\\[-8pt]
\\[-8pt]
 \eqntext{d\nu\mbox{-a.e. for all }%
s\in \{ s\in [ 0,t ] \dvtx x_{s}=x_{s^{-}} \} .}
\end{eqnarray}
From the assumption in \eqref{cone1} we deduce that $f$ is strictly
positive. Thus, by the convexity of $\Gamma_{s} ( x_{s} ) $ and the
absolute continuity of $ \vert\lambda \vert$ with respect to $%
\nu$, we conclude that%
%
\begin{equation}
\gamma_{s}\in\Gamma_{s}^{1} ( x_{s} ) ,\qquad d \vert\lambda
\vert\mbox{-a.e. for all }s\in \{ s\in [ 0,t ]
\dvtx x_{s}=x_{s^{-}} \} . \label{ddd3}
\end{equation}
In particular, \eqref{ddd3} verifies the second part in (\ref{SP2}) and,
hence, the proof in Case~2 is also complete.

%

Having completed the proof of Cases 1 and   2, we
conclude that the proofs of~\eqref{BV} and \eqref{prop} are complete. Hence,
to complete the proof of Theorem \ref{SPconvergence}, we only have to ensure
that $x\in\mathcal{D}^{\rho_{0}} ( [0,T],%
\mathbb{R}
) $. However, this follows from the assumption that $x^{n}\in\mathcal{
D}^{\rho_{0}} ( [ 0,T ] ,%
\mathbb{R}
) $ for all $n\geq1$ and from the fact that $x^{n}\rightarrow x$ in
the Skorohod topology.
\end{pf*}

\subsection{Convergence of a sequence of solutions to approximating Skorohod
problems}

Let $T>0$ and let $D\subset
\mathbb{R}
^{d+1}$ be a time-dependent domain satisfying \eqref{timedep+}, %
\eqref{limitzero} and a uniform exterior sphere condition in time with
radium $r_{0}$ in the sense of \eqref{extsphere-}. Let $\Gamma=\Gamma
_{t} ( z ) $ be a closed convex cone of vectors in $%
\mathbb{R}
^{d}$ for every $z\in\partial D_{t}$, $t\in [ 0,T ] $, and assume
that $\Gamma$ satisfies \eqref{cone1}, \eqref{cone2} and \eqref{Gammanpre}.
Assume that \eqref{crita} and \eqref{crite} hold for some $0<\rho
_{0}<r_{0}$%
, $\eta_{0}>0$, $a$ and $e$. Finally, assume that $ ( [ 0,T ]
\times
\mathbb{R}
^{d} ) \setminus\overline{D}$ has the $ ( \delta_{0},h_{0} )
$-property of good projections along~$\Gamma$, for some $0<\delta
_{0}<\rho
_{0}$, $h_{0}>1$ and let $w\in\mathcal{D}^{( {\delta
_{0}}/{4}\wedge
 {\rho_{0}}/({4h_{0}}))} ( [0,T],%
\mathbb{R}
^{d} ) $ with $w_{0}\in\overline{D_{0}}$. The purpose of this
subsection is to construct a sequence of solutions to Skorohod problems
which approximate the Skorohod problem for $ ( D,\Gamma,w ) $.
Based on this sequence, in the next section we conclude the existence
of a
solution $ ( x,\lambda ) $ to the Skorohod problem for $ (
D,\Gamma,w ) $, in the sense of Definition \ref{skorohodprob}, with $%
x\in\mathcal{D}^{\rho_{0}} ( [0,T],%
\mathbb{R}
^{d} ) $. This will then complete the proof of Theorem \ref{Theorem 1}.
To proceed, we let $n\in
\mathbb{N}
$, $n\gg1$, and we let $ \{ \epsilon_{n} \} $ be a sequence of
real numbers which tends to $0$ as $n\rightarrow\infty$. Then, for
each $n$%
, we can find a partition $\Delta_{n}= \{ \tau_{k}^{n} \}
_{k=0}^{N_{n}}$ of the interval $ [ 0,T ] $, that is, $0=\tau
_{0}^{n}<\tau_{1}^{n}<\cdots<\tau_{N_{n}-1}^{n}<\tau
_{N_{n}}^{n}=T$, such
that~\eqref{korlim++}--\eqref{korlim+++} stated below hold. In
particular,%
%
\begin{equation}
\lim_{n\rightarrow\infty}\Delta_{n}^{\ast}=0  \qquad\mbox{where }\Delta
_{n}^{\ast}:=\max_{k\in \{ 1,\ldots,N_{n}-1 \} }\tau
_{k+1}^{n}-\tau
_{k}^{n}, \label{korlim++}
\end{equation}
and, for some $n_{0}\gg1$,%
%
\begin{eqnarray}\label{korlim--lex}
\Vert w \Vert_{\tau_{k}^{n},\tau_{k+1}^{n}}+l ( \Delta
_{n}^{\ast} ) <\min \biggl\{ \frac{\delta_{0}}{2},\frac{\rho_{0}}{%
2h_{0}} \biggr\} \nonumber
\\[-8pt]
\\[-8pt]
\eqntext{\mbox{whenever }n\geq n_{0},k\in \{
0,\ldots,N_{n}-1 \} .}
\end{eqnarray}
Furthermore, we define $w^{\Delta_{n}}=w_{t}^{\Delta_{n}}$, $t\in
\lbrack
0,T]$, as
%
\begin{equation}
w_{t}^{\Delta_{n}}=w_{\tau_{k}^{\Delta_{n}}}  \qquad\mbox{whenever }t\in [
\tau_{k}^{n},\tau_{k+1}^{n} ) , k\in \{
0,\ldots,N_{n}-1 \} , \label{yyy1}
\end{equation}
and $w_{T}^{\Delta_{n}}=w_{T}$, so that
%
\begin{eqnarray}\label{korlim+++}
w^{\Delta_{n}}&\in&\mathcal{D}^{( {\delta_{0}}/{4}\wedge
{\rho_{0}%
}/({4h_{0}}))} ( [ 0,T ] ,\mathbb{R}
^{d} ) , \nonumber
\\[-8pt]
\\[-8pt]
 w_{0}^{\Delta_{n}}&\in&\overline{D_{0}}\quad \mbox{and}\quad d_{%
\mathcal{D}}([0,T],w^{\Delta_{n}},w)\leq\epsilon_{n}.\nonumber
\end{eqnarray}
In particular, $w^{\Delta_{n}}\in\mathcal{D} ( [ 0,T ] ,%
\mathbb{R}
^{d} ) $ is a step function which approximates $w$ in the Skorohod
topology. Given $\Delta_{n}$ and $w^{\Delta_{n}}$, we define
$D^{\Delta
_{n}}$ and $\Gamma^{\Delta_{n}}$ as in \eqref{approxw+}.
Furthermore, to
obtain a more simple notation, from now on we write $w^{n}$, $D^{n}$
and $%
\Gamma^{n}$ for $w^{\Delta_{n}}$, $D^{\Delta_{n}}$ and $\Gamma
^{\Delta
_{n}}$. Then, given $w^{n}$, $D^{n}$ and $\Gamma^{n}$ as above, we next
define a pair of processes $ ( x^{n},\lambda^{n} ) $ as follows.
Let%
%
\begin{equation}
x_{t}^{n}=w_{0},\qquad \lambda_{t}^{n}=0 \qquad\mbox{for }t\in [ 0,\tau
_{1}^{n} ) . \label{approxw++a}
\end{equation}
If $x_{\tau_{k-1}^{n}}^{n}\in\overline{D_{\tau_{k-1}^{n}}^{n}}$ for
some $%
k\in \{ 1,\ldots,N_{n} \} $, then, by the triangle\vspace*{1pt} inequality and~\eqref{korlim--lex},
%
\begin{equation}
d ( x_{\tau_{k-1}^{n}}^{n}+w_{\tau_{k}^{n}}^{n}-w_{\tau
_{k-1}^{n}}^{n},\overline{D_{\tau_{k}^{n}}^{n}} ) \leq
\Vert w^{n}\Vert _{\tau
_{k-1}^{n},\tau_{k}^{n}}+l(\Delta_{n}^{\ast})<\delta_{0}. \label{jj}
\end{equation}
Hence, by the $ ( \delta_{0},h_{0} ) $-property of good
projections, it follows that if\vspace*{-5pt} $x_{\tau_{k-1}^{n}}^{n}+w_{\tau
_{k}^{n}}^{n}-w_{\tau_{k-1}^{n}}^{n}\notin\overline{D_{\tau_{k}^{n}}^{n}}$,
then there exists a point
%
\begin{equation}
\pi_{\partial D_{\tau_{k}^{n}}^{n}}^{\Gamma_{\tau_{k}^{n}}^{n}} (
x_{\tau_{k-1}^{n}}^{n}+w_{\tau_{k}^{n}}^{n}-w_{\tau_{k-1}^{n}}^{n} )
\in\partial{D_{\tau_{k}^{n}}^{n},} \label{jjj}
\end{equation}
which is the projection of $x_{\tau_{k-1}^{n}}^{n}+w_{\tau
_{k}^{n}}^{n}-w_{\tau_{k-1}^{n}}^{n}$ onto $\partial{D_{\tau_{k}^{n}}^{n}}
$ along $\Gamma_{\tau_{k}^{n}}$. Furthermore, if $x_{\tau
_{k-1}^{n}}^{n}+w_{\tau_{k}^{n}}^{n}-w_{\tau_{k-1}^{n}}^{n}\in
\overline{%
D_{\tau_{k}^{n}}^{n}}$, then we let
%
\begin{equation}
\pi_{\partial D_{\tau_{k}^{n}}^{n}}^{\Gamma_{\tau_{k}^{n}}^{n}} (
x_{\tau_{k-1}^{n}}^{n}+w_{\tau_{k}^{n}}^{n}-w_{\tau_{k-1}^{n}}^{n} )
=x_{\tau_{k-1}^{n}}^{n}+w_{\tau_{k}^{n}}^{n}-w_{\tau_{k-1}^{n}}^{n}.
\label{jjjj}
\end{equation}
Based on this argument, we define, whenever $t\in [ \tau_{k}^{n},\tau
_{k+1}^{n} ) $, $k\in \{ 1,\ldots,N_{n}-1 \} $,
%
\begin{eqnarray}\label{approxw+++b}
x_{t}^{n} &=&\pi_{\partial D_{\tau_{k}^{n}}^{n}}^{\Gamma_{\tau
_{k}^{n}}^{n}} ( x_{\tau_{k-1}^{n}}^{n}+w_{\tau_{k}^{n}}^{n}-w_{\tau
_{k-1}^{n}}^{n} ) ,\nonumber
\\[-8pt]
\\[-8pt]
\lambda_{t}^{n} &=&\lambda_{\tau_{k-1}^{n}}^{n}+\bigl(x_{t}^{n}-(x_{\tau
_{k-1}^{n}}^{n}+w_{\tau_{k}^{n}}^{n}-w_{\tau_{k-1}^{n}}^{n})\bigr),\nonumber
\end{eqnarray}
and, finally, we define $x_{T}^{n}$ and $\lambda_{T}^{n}$ using %
\eqref{approxw+++b} by simply setting $k=N_{n}$ in \eqref
{approxw+++b}. Note
that in this way we have $x_{\tau_{k-1}^{n}}^{n}\in\overline{D_{\tau
_{k-1}^{n}}^{n}}$ for all $k\in \{ 1,\ldots,N_{n} \} $. Next, again
using the $ ( \delta_{0},h_{0} ) $-property of good projections,
we see that%
%
\begin{eqnarray}\label{bound3}
\vert x_{\tau_{k}^{n}}^{n}-x_{\tau_{k-1}^{n}}^{n} \vert&\leq
& \vert\pi_{\partial D_{\tau_{k}^{n}}^{n}}^{\Gamma_{\tau
_{k}^{n}}^{n}} ( x_{\tau_{k-1}^{n}}^{n}+w_{\tau_{k}^{n}}^{n}-w_{\tau
_{k-1}^{n}}^{n} ) - ( x_{\tau_{k-1}^{n}}^{n}+w_{\tau
_{k}^{n}}^{n}-w_{\tau_{k-1}^{n}}^{n} ) \vert \nonumber\\
&&{}+ \vert w_{\tau_{k}^{n}}^{n}-w_{\tau_{k-1}^{n}}^{n} \vert
\nonumber\\
&\leq&h_{0}d ( x_{\tau_{k-1}^{n}}^{n}+w_{\tau_{k}^{n}}^{n}-w_{\tau
_{k-1}^{n}}^{n},\overline{D_{\tau_{k}^{n}}^{n}} ) + \vert w_{\tau
_{k}^{n}}^{n}-w_{\tau_{k-1}^{n}}^{n} \vert \\
&\leq&h_{0} \bigl( \Vert w^{n}\Vert _{\tau_{k-1}^{n},\tau_{k}^{n}}+l(\Delta
_{n}^{\ast}) \bigr) +\Vert w^{n}\Vert _{\tau_{k-1}^{n},\tau_{k}^{n}} \nonumber\\
&\leq&h_{0} \biggl( \frac{\rho_{0}}{2h_{0}} \biggr) +\frac{\delta_{0}}{4}%
<\rho_{0}.\nonumber
\end{eqnarray}
Hence, $x^{n}\in\mathcal{D}^{\rho_{0}} ( [ 0,T ] ,%
\mathbb{R}
^{d} ) $. Using this notation, we next prove the following theorem.

\begin{theorem}
\label{SPconvergence+}Let $T>0$, $D\subset
\mathbb{R}
^{d+1}$, $r_{0}$, $\Gamma=\Gamma_{t} ( z ) $, $0<\rho_{0}<r_{0}$%
, $\eta_{0}>0$, $a$, $e$, $\delta_{0}$ and $h_{0}$ be as in the statement
of Theorem \ref{Theorem 1}. Let $w$ be as in the statement of Theorem
\ref%
{Theorem 1} and let $w^{n}$, $D^{n}$, $\Gamma^{n}$, $x^{n}$ and
$\lambda
^{n}$ be defined as above for $n\geq1$. Then $ ( x^{n},\lambda
^{n} ) $ is a solution to the Skorohod problem for $ ( D^{n},\Gamma
^{n},w^{n} ) $ and $x^{n}\in\mathcal{D}^{\rho_{0}} ( [ 0,T%
] ,%
\mathbb{R}
^{d} ) $ for all $n\geq n_{0}$ for some $n_{0}\in
\mathbb{N}
$. Moreover, $ \{ ( x^{n},\lambda^{n} ) \} $ converges
to $ ( x,\lambda ) \in\mathcal{D} ( [ 0,T ] ,%
\mathbb{R}
^{d} ) \times\mathcal{BV} ( [ 0,T ] ,%
\mathbb{R}
^{d} ) $ and $ ( x,\lambda ) $ is a solution to the Skorohod
problem for $ ( D,\Gamma,w ) $. Furthermore, $x\in\mathcal{D}%
^{\rho_{0}} ( [ 0,T ] ,%
\mathbb{R}
^{d} ) $.
\end{theorem}

\begin{remark}
\label{subtel2} Note that for Theorem \ref{SPconvergence+} we, in contrast
to in Theorem \ref{SPconvergence}, also need to assume \eqref
{Gammanpre} in
order to be able complete the proof [see \eqref{bbbb2again} below].
\end{remark}

\begin{pf*}{Proof of Theorem \ref{SPconvergence+}}
$ ( x^{n},\lambda^{n} ) $ is, by
construction, a solution to the Skorohod problem for $ ( D^{n},\Gamma
^{n},w^{n} ) $ and the statement that $x^{n}\in\mathcal{D}^{\rho
_{0}} ( [ 0,T ] ,%
\mathbb{R}
^{d} ) $ for all $n\geq n_{0}$, for some $n_{0}\in
\mathbb{N}
$, is proved in \eqref{bound3}. Next, using Theorem \ref%
{compacttheoremapprox}, we can conclude the existence of some positive
constants $\hat{L}_{1} ( w,T ) $, $\hat{L}_{2} ( w,T ) $, $%
\hat{L}_{3} ( w,T ) $ and $\hat{L}_{4} ( w,T ) $ such that%
%
\begin{eqnarray}\label{lll1+}
\Vert x^{n} \Vert_{t_{1},t_{2}} &\leq&\hat{L}_{1} (
w,T ) \Vert w \Vert_{t_{1},t_{2}}+\hat{L}_{2} (
w,T ) \bigl( l ( t_{2}-t_{1} ) +l ( \Delta_{n}^{\ast
} ) \bigr) ,\nonumber
\\[-8pt]
\\[-8pt]
\quad \vert\lambda^{n} \vert_{t_{2}}- \vert\lambda
^{n} \vert_{t_{1}} &\leq&\hat{L}_{3} ( w,T ) \Vert
w \Vert_{t_{1},t_{2}}+\hat{L}_{4} ( w,T ) \bigl( l (
t_{2}-t_{1} ) +l ( \Delta_{n}^{\ast} ) \bigr) ,\nonumber
\end{eqnarray}
whenever $0\leq t_{1}\leq t_{2}\leq T$. In particular, note that by choosing
$n_{0}$ sufficiently large we can ensure, using \eqref{korlim++}, that
$%
l ( \Delta_{n}^{\ast} ) \leq\rho_{0}/ ( 4 ( K_{2} (
a,e ) +1 ) ) $ and that \eqref{ahah1+} holds for all $n\geq
n_{0}$. Hence, Theorem \ref{compacttheoremapprox} is applicable. Based
on %
\eqref{lll1+}, we can now argue as in the proof of Theorem \ref
{SPconvergence}%
. In particular, as $ \{ w^{n} \} $ converges to $w\in\mathcal{D}%
( [ 0,T ] ,%
\mathbb{R}
^{d} ) $, we see, using \eqref{lll1+}, that $ \{ (
w^{n},x^{n},\lambda^{n}, \vert\lambda \vert^{n} ) \}
$ is relatively compact in $\mathcal{D} ( [ 0,T ] ,%
\mathbb{R}
^{d} ) \times\mathcal{D} ( [ 0,T ] ,%
\mathbb{R}
^{d} ) \times\mathcal{D} ( [ 0,T ] ,%
\mathbb{R}
^{d} ) \times\mathcal{D} ( [ 0,T ] ,%
\mathbb{R}
_{+} ) $. Furthermore, we know that $x_{t}^{n}\in\overline{D_{n}}$ for
all $t\in [ 0,T ] $, $n\geq1$. Hence, $ \{ (
x^{n},\lambda^{n} ) \} $ converges to some $ ( x,\lambda
) \in\mathcal{D} ( [ 0,T ] ,%
\mathbb{R}
^{d} ) \times\mathcal{D} ( [ 0,T ] ,%
\mathbb{R}
^{d} ) $. We intend to prove that $ ( x,\lambda ) \in
\mathcal{D} ( [ 0,T ] ,%
\mathbb{R}
^{d} ) \times\mathcal{BV} ( [ 0,T ] ,%
\mathbb{R}
^{d} ) $ solves the Skorohod problem for $ ( D,\Gamma,w ) $
and, to do this, we have to prove that
%
\begin{equation}
\lambda\in\mathcal{BV} ( [ 0,T ] ,%
\mathbb{R}
^{d} ) , \label{BV+}
\end{equation}
and we have to verify
%
\begin{equation}
\mbox{properties \eqref{SP1}--\eqref{SP4} in Definition \ref{skorohodprob}.}
\label{prop+}
\end{equation}
The proof of \eqref{BV+} and \eqref{prop+} follows along the lines of the
proof of \eqref{BV} and~\eqref{prop} in Theorem \ref{SPconvergence}
and we
shall only outline the main differences between the proofs. To start with,
the statements in \eqref{hra1} and \eqref{convTh3} remain true.
However, the
argument in \eqref{bbbb1-} has to be changed. In particular, in this
case we
see, using~\eqref{convTh3} and \eqref{korlim++}, that%
%
\begin{eqnarray}\label{bbbb1-+}
d ( x_{t},\overline{D_{t}} ) &\leq&|x_{t}-x_{\tilde{t}_{n}}^{n}|+h(%
\overline{D_{\tilde{t}_{n}}^{n}},\overline{D_{\tilde
{t}_{n}}})\nonumber
\\[-8pt]
\\[-8pt]
&&{}+l(|\tilde{t}%
_{n}-t|)+l(\Delta_{n}^{\ast})\rightarrow0  \qquad\mbox{as }n\rightarrow
\infty,\nonumber
\end{eqnarray}
which, since $\overline{D_{t}}$ is closed, proves that $x_{t}\in
\overline{%
D_{t}}$, for all $t\in [ 0,T ] $. Hence, we have verified that $%
( w_{t},x_{t},\lambda_{t} ) $ satisfies \eqref{SP1}. Next,
arguing as in \eqref{defR+}--\eqref{laba2}, using \eqref{lll1+} to
conclude %
\eqref{laba1}, we can conclude that \eqref{BV+} holds and that
%
\begin{equation}
\lambda_{t}=\int_{0}^{t}\gamma_{s}\,d \vert\lambda \vert
_{s}=\int_{G_{t}}\gamma\,d\mu ( s,z,\gamma )   \qquad\mbox{whenever }%
t\in [ 0,T ]
\end{equation}
for some $S_{1} ( 0 ) $-valued Borel measurable function $\gamma
_{s}$. Hence, to prove \eqref{SP2}, we again have to prove that
$\gamma
_{s}\in\Gamma_{s}^{1} ( x_{s} ) $ for all $t\in [ 0,T ]
$. As in the proof of Theorem~\ref{SPconvergence}, we consider Case 1 and
Case 2. In fact, Case 2 can be handled exactly as in the proof of
Theorem %
\ref{SPconvergence} and hence shall only discuss the proof of Case 1. To
prove Case 1, we first see that the statements in \eqref{hra2} and %
\eqref{convTh2} can be repeated and, arguing as in \eqref{bbbb1-+},
we obtain%
%
\begin{eqnarray}
d ( x_{t},\partial D_{t} ) &\leq&|x_{t}-x_{\tilde{t}%
_{n}}^{n}|+h(\partial D_{\tilde{t}_{n}}^{n},\partial D_{\tilde
{t}_{n}})\nonumber
\\[-8pt]
\\[-8pt]
&&{}+\hat{%
l}(|\tilde{t}_{n}-t|)+\hat{l}(\Delta_{n}^{\ast})\rightarrow0  \qquad\mbox
{as }%
n\rightarrow\infty.\nonumber
\end{eqnarray}
Hence, since $\partial D_{t}$ is closed, we can conclude that $x_{t}\in
\partial D_{t}$. To proceed, we deduce as in \eqref{bbbb2} that%
%
\begin{equation}\qquad
d ( \lambda_{t}-\lambda_{t^{-}},\Gamma_{t}(x_{t}) ) \leq
|\lambda_{t}-\lambda_{\tilde{t}_{n}}^{n}|+|\lambda_{t^{-}}-\lambda
_{%
\tilde{t}_{n}^{-}}^{n}|+h(\Gamma_{\tilde{t}_{n}}^{n}(x_{\tilde{t}%
_{n}}^{n}),\Gamma_{t}(x_{t})). \label{bbbb2again}
\end{equation}
Obviously the first two terms on the right-hand side in \eqref{bbbb2again}
tend to zero as $n\rightarrow\infty$. Concerning the third term, we first
note that there exists, for $n$ large enough, some integer $k ( n )
$ such that
%
\begin{equation}
\Gamma_{\tilde{t}_{n}}^{n}(x_{\tilde{t}_{n}}^{n})=\Gamma_{\tau
_{k(n)}^{n}}^{n}\bigl(x_{\tau_{k(n)}^{n}}^{n}\bigr). \label{hgf}
\end{equation}
Hence, as $|\tilde{t}_{n}-\tau_{k(n)}^{n}|\leq l(\Delta_{n}^{\ast
})$, we
can conclude, using \eqref{convTh3}, that $\tau_{k(n)}^{n}\rightarrow
t$ as
$n\rightarrow\infty$. Moreover, as $x_{\tau
_{k(n)}^{n}}^{n}=x_{^{\tilde{t}%
_{n}}}^{n}$, we can also use \eqref{convTh3} to conclude that\vspace*{1pt} $x_{\tau
_{k(n)}^{n}}^{n}\rightarrow x_{t}$ as $n\rightarrow\infty$. In particular,
based on these conclusions, it follows from \eqref{Gammanpre} that
also the
third term on the right-hand side in \eqref{bbbb2again} tends to zero
as $%
n\rightarrow\infty$. Hence, having proved that $d ( \lambda
_{t}-\lambda_{t^{-}},\Gamma_{t} ( x_{t} ) ) =0$, the proof
of Case 1 can now be completed as in Theorem \ref{SPconvergence}.

%

Having completed the proof of Cases 1 and   2, we can conclude
that the proofs of \eqref{BV+} and \eqref{prop+} are complete. Hence, to
complete the proof of Theorem \ref{SPconvergence+}, we only have to ensure
that $x\in\mathcal{D}^{\rho_{0}} ( [ 0,T ] ,%
\mathbb{R}
^{d} ) $. However, again this follows from the fact that $x^{n}\in
\mathcal{D}^{\rho_{0}} ( [ 0,T ] ,%
\mathbb{R}
^{d} ) $ for all $n\geq n_{0}$ and from the fact that $x^{n}\rightarrow
x$ in the Skorohod topology.
\end{pf*}

\section{\texorpdfstring{Proof of Theorems \protect\ref{Theorem 1}, \protect\ref
{Theorem 2} and \protect\ref{Theorem 3}}{Proof of Theorems 1.2, 1.3 and 1.9}} \label{sectskor}

\mbox{}

\begin{pf*}{Proof of Theorem \ref{Theorem 1}}
Theorem \ref{Theorem 1}
now follows immediately from Theorem \ref{SPconvergence+}.
\end{pf*}
%

\begin{pf*}{Proof of Theorem \ref{Theorem 2}} Using Theorem \ref%
{compacttheorem} and \eqref{limitzero}, we see that
%
\begin{eqnarray}
\lim_{t_{2}\rightarrow t_{1}} \vert x_{t_{2}}-x_{t_{1}} \vert&\leq&
\lim_{t_{2}\rightarrow t_{1}} \Vert x \Vert_{t_{1},t_{2}}\nonumber\\
 &\leq
&\lim_{t_{2}\rightarrow t_{1}}\bigl(L_{1} \Vert w \Vert
_{t_{1},t_{2}}+L_{2}l ( |t_{2}-t_{1}| ) \bigr) \\
&\leq&0+L_{2}\lim_{t_{2}\rightarrow t_{1}}l ( |t_{2}-t_{1}| ) =0.\nonumber
\end{eqnarray}
This proves that $x$ is continuous.
\end{pf*}
%

\begin{pf*}{Proof of Theorem \ref{Theorem 3}} Let $W$ be a $m$%
-dimensional Wiener process on a filtered probability space $ ( \Omega,
\mathcal{F,} \{ \mathcal{F}_{t} \} ,P ) $ and in the following let
$W_{t}$, $t\in\lbrack0,T]$, be a continuous path of $W$. We
define, for $n\in
\mathbb{N}
$, $n\gg1$, $k\in \{ 0,1,\ldots,2^{n}-1 \} $,
%
\begin{eqnarray}\label{dsds}
D_{t}^{n}&=&D_{kT/2^{n}},\nonumber
\\[-8pt]
\\[-8pt]
 \Gamma_{t}^{n}&=&\Gamma_{kT/2^{n}} \qquad
\mbox{whenever }%
t\in [ kT/2^{n}, ( k+1 ) T/2^{n} ) ,\nonumber
\end{eqnarray}
and $D_{T}^{n}=D_{T}, \Gamma_{T}^{n}=\Gamma_{T}$. Furthermore, we
recursively define three processes $X^{n}=X_{t}^{n}$, $Z^{n}=Z_{t}^{n}$
and $%
\Lambda^{n}=\Lambda_{t}^{n}$, for $t\in\lbrack0,T]$, in the following
way. Let
%
\begin{equation}
X_{0}^{n}=\hat{z},\qquad Z_{0}^{n}=\hat{z},\qquad \Lambda_{0}^{n}=0,
\label{dsds1}
\end{equation}
and let, for $k\in \{ 0,1,\ldots,2^{n}-1 \} $,
%
\begin{eqnarray}\label{dsds2}
Z_{(k+1)T/2^{n}}^{n}
&=&Z_{kT/2^{n}}^{n}+\frac{T}{2^{n}}b (
kT/2^{n},X_{kT/2^{n}}^{n} ) \nonumber \\
&&{}+\sigma ( kT/2^{n},X_{kT/2^{n}}^{n} ) \bigl( W_{ (
k+1 ) T/2^{n}}-W_{kT/2^{n}} \bigr) , \\
X_{ ( k+1 ) T/2^{n}}^{n}
&=&\pi_{\partial
D_{(k+1)T/2^{n}}}^{\Gamma_{(k+1)T/2^{n}}} \bigl( X_{kT/2^{n}}^{n}+Z_{ (
k+1 ) T/2^{n}}^{n}-Z_{kT/2^{n}}^{n} \bigr) .\nonumber
\end{eqnarray}
We here have to make sure that $X_{ ( k+1 ) T/2^{n}}^{n}$ is well
defined. To do this, we note that, either\vspace*{1pt} $X_{kT/2^{n}}^{n}+Z_{ (
k+1 ) T/2^{n}}^{n}-Z_{kT/2^{n}}^{n}\in\overline{D_{(k+1)T/2^{n}}^{n}}$
or $X_{kT/2^{n}}^{n}+Z_{ ( k+1 ) T/2^{n}}^{n}-Z_{kT/2^{n}}^{n}\in
\mathbb{R}
^{d}\setminus\overline{D_{(k+1)T/2^{n}}^{n}}$. In the first case we
identify the projection with the\vspace*{1pt} point itself, whereas, in the latter case,
we have to assert the existence of appropriate projections onto
$\partial
D_{(k+1)T/2^{n}}^{n}$. However, assuming $X_{kT/2^{n}}^{n}\in
D_{kT/2^{n}}^{n}$, we see that
%
\begin{eqnarray}\label{dsds4}
&&d\bigl(X_{kT/2^{n}}^{n}+Z_{ ( k+1 )
T/2^{n}}^{n}-Z_{kT/2^{n}}^{n},D_{(k+1)T/2^{n}}^{n}\bigr) \nonumber\\
&&\qquad\leq l(T/2^{n})+\bigl|Z_{ ( k+1 ) T/2^{n}}^{n}-Z_{kT/2^{n}}^{n}\bigr|
\\
&&\qquad\leq l(T/2^{n})+({T}/{2^{n}})\Bigl(\sup_{\overline{D}} \vert b \vert
\Bigr)+\Bigl(\sup_{\overline{D}} \Vert\sigma \Vert\Bigr)\bigl|W_{ ( k+1 )
T/2^{n}}-W_{kT/2^{n}}\bigr|.\nonumber
\end{eqnarray}
Hence, since $W_{t}$ is a continuous path, there must exist some
$n_{0}\in
\mathbb{%
\mathbb{N}
}$ such that
%
\begin{equation}
d\bigl(X_{kT/2^{n}}^{n}+Z_{ ( k+1 )
T/2^{n}}^{n}-Z_{kT/2^{n}}^{n},D_{(k+1)T/2^{n}}^{n}\bigr)<\delta_{0},
\label{dsds5}
\end{equation}
whenever $n\geq n_{0}$ and $k\in \{ 0,1,\ldots,2^{n}-1 \} $. By %
\eqref{dsds5} and the $ ( \delta_{0},h_{0} ) $-property of good
projections, it follows that the projection $\pi_{\partial
D_{(k+1)T/2^{n}}}^{\Gamma_{(k+1)T/2^{n}}}(X_{kT/2^{n}}^{n}+Z_{ (
k+1 ) T/2^{n}}^{n}-Z_{kT/2^{n}}^{n})$ is well defined, for $n\geq n_{0}$
, whenever $X_{kT/2^{n}}^{n}+Z_{ ( k+1 )
T/2^{n}}^{n}-Z_{kT/2^{n}}^{n}\in
\mathbb{R}
^{d}\setminus\overline{D_{(k+1)T/2^{n}}^{n}}$. Furthermore, using the
definition of $Z_{ ( k+1 ) T/2^{n}}^{n}$ in \eqref{dsds2}, we also
see that
%
\begin{eqnarray}\label{dsds3+}
\bigl|Z_{ ( k+1 ) T/2^{n}}^{n}-Z_{kT/2^{n}}^{n}\bigr| &\leq&l(T/2^{n})+({T}/{%
2^{n}})\Bigl(\sup_{\overline{D}} \vert b \vert\Bigr) \nonumber
\\[-8pt]
\\[-8pt]
&&{}+\Bigl(\sup_{\overline{D}} \Vert\sigma \Vert\Bigr)\bigl|W_{ ( k+1 )
T/2^{n}}-W_{kT/2^{n}}\bigr|,\nonumber
\end{eqnarray}
and, hence, once more using that $W_{t}$ is a continuous path, we can ensure
that
%
\begin{eqnarray}\label{dsds4ha1}
\mathrm{(i)} &&\hspace*{5pt}Z^{n}\in\mathcal{D}^{( {\delta_{0}}/{4}\wedge {\rho
_{0}}/({4h_{0}}))} ( [ 0,T ] ,%
\mathbb{R}
^{d} ) , \nonumber
\\[-8pt]
\\[-8pt]
\mathrm{(ii)} &&\hspace*{5pt}h_{0} \bigl( \Vert Z^{n}\Vert _{kT/2^{n},(k+1)T/2^{n}}+l(\Delta
_{n}^{\ast}) \bigr) +\Vert Z^{n}\Vert _{kT/2^{n},(k+1)T/2^{n}}\leq\rho_{0},\nonumber
\end{eqnarray}
whenever $n\geq n_{0}$ and $k\in \{ 0,1,\ldots,2^{n}-1 \} $. We next
let, for $k\in \{ 0,1,\ldots,2^{n}-1 \} $,
%
\begin{equation}
\quad \Lambda_{ ( k+1 ) T/2^{n}}^{n}=\Lambda_{kT/2^{n}}^{n}+X_{ (
k+1 ) /2^{n}}^{n}-X_{kT/2^{n}}^{n}-Z_{ ( k+1 )
T/2^{n}}^{n}+Z_{kT/2^{n}}^{n}. \label{dsds6}
\end{equation}
Finally, we define, for $kT/2^{n}\leq t< ( k+1 ) T/2^{n}$, $k\in
\{ 0,1,\ldots,2^{n}-1 \} $,
%
\begin{equation}
X_{t}^{n}=X_{kT/2^{n}}^{n},\qquad Z_{t}^{n}=Z_{kT/2^{n}}^{n},\qquad \Lambda
_{t}^{n}=\Lambda_{kT/2^{n}}^{n}.
\end{equation}
Then, by arguing as in the proof of \eqref{bound3}, using (i) and (ii)
in \eqref{dsds4ha1}, we can conclude that%
%
\begin{equation}
X^{n}\in\mathcal{D}^{\rho_{0}} ( [ 0,T ] ,%
\mathbb{R}
^{d} )  \qquad \mbox{whenever }n\geq n_{0}. \label{boliden}
\end{equation}
Furthermore, using the definitions above, it is clear that
%
\begin{equation}
Z_{t}^{n}=\hat{z}+\int_{0}^{t}b ( s,X_{s}^{n} )
\,ds+\int_{0}^{t}\sigma ( s,X_{s}^{n} ) \,dW_{s}-\varepsilon
^{n} ( t ) , \label{RSDE2}
\end{equation}
where%
%
\begin{equation}\qquad
\sup_{t\in [ 0,T ] } \vert\varepsilon^{n} ( t )
\vert
\leq\Bigl(\sup_{\overline{D}} \vert b \vert\Bigr)\frac{T}{2^{n}}%
+\Bigl(\sup_{\overline{D}} \Vert\sigma \Vert\Bigr)\mathop{\mathop{\sup}_{0\leq s\leq t\leq T ,}}_{\vert
s-t \vert\leq {T}/{2^{n}}} \vert
W_{t}-W_{s} \vert. \label{RSDE1}
\end{equation}
By construction, $ ( X^{n},\Lambda^{n} ) $ solves the Skorohod
problem for $ ( D^{n},\Gamma^{n},Z^{n} ) $ and using Theorem \ref%
{compacttheoremapprox}, we can conclude that there exist positive
constants $%
\hat{L}_{1} ( Z,T ) $, $\hat{L}_{2} ( Z,T ) $, $\hat{L}%
_{3} ( Z,T ) $ and $\hat{L}_{4} ( Z,T ) $, independent of $%
n$, for $n\geq n_{0}$, such that
%
\begin{eqnarray}\label{jjja}\qquad
\Vert X^{n} \Vert_{t_{1},t_{2}} &\leq&\hat{L}_{1} (
Z,T ) \Vert Z \Vert_{t_{1},t_{2}}+\hat{L}_{2} (
Z,T ) \bigl(l ( t_{2}-t_{1} ) +l(T/2^{n})\bigr), \nonumber
\\[-8pt]
\\[-8pt]\qquad
|\Lambda^{n}|_{t_{2}}-|\Lambda^{n}|_{t_{1}} &\leq&\hat{L}_{3} (
Z,T ) \Vert Z \Vert_{t_{1},t_{2}}+\hat{L}_{4} (
Z,T ) \bigl(l ( t_{2}-t_{1} ) +l(T/2^{n})\bigr),\nonumber
\end{eqnarray}
whenever $0\leq t_{1}\leq t_{2}\leq T$ and $n\geq n_{0}$. Hence, the
sequence $\{ ( Z^{n},X^{n},\Lambda^{n} ) \}$ is relatively compact
in the Skorohod topology and we can conclude, by the construction of $Z^{n}$
and (\ref{RSDE2})--(\ref{RSDE1}), that also the sequence $ \{ (
W,Z^{n},X^{n},\break\Lambda^{n},\varepsilon^{n} ) \} $ is relatively
compact in the Skorohod topology. In fact, as this argument can be repeated
for each continuous path of $W_{t}$, it follows, by considering convergent
subsequences if necessary, that the sequence of vector valued processes
$%
\{ ( W,Z^{n},X^{n},\Lambda^{n},\varepsilon^{n} ) \} $
defined on $ ( \Omega,\mathcal{F},P ) $ converges in law to a
stochastic process $ ( W,Z,X,\Lambda,0 ) $. Furthermore, using the
Skorohod representation theorem (see, e.g., \cite
{Billingsley1999} and
\cite{EthierKurtz1986}), there exists a complete probability space
$(\tilde{%
\Omega},\mathcal{\tilde{F}},\tilde{P})$ and versions $\{(\tilde
{W}^{n},%
\tilde{Z}^{n},\tilde{X}^{n},\tilde{\Lambda}^{n},\tilde{\varepsilon
}^{n})\}$
and $(\tilde{W},\tilde{Z},\tilde{X},\tilde{\Lambda},0)$ of $ \{
(W,Z^{n},X^{n},\Lambda^{n},\varepsilon^{n}) \} $ and $(W,Z,X,\Lambda
,0)$ on $(\tilde{\Omega},\mathcal{\tilde{F}},\tilde{P})$, such
that $\{(%
\tilde{W}^{n},\tilde{Z}^{n},\break \tilde{X}^{n},\tilde{\Lambda
}^{n},\tilde{%
\varepsilon}^{n})\}$ converges to $(\tilde{W},\tilde{Z},\tilde
{X},\tilde{%
\Lambda},0)$ $\tilde{P}$-almost surely. Moreover, using Theorem \ref%
{SPconvergence+} and the fact that $(\tilde{X}^{n},\tilde{\Lambda}^{n})$
solves, $\tilde{P}$-almost surely, the Skorohod problem for
$(D^{n},\Gamma
^{n},\tilde{Z}^{n})$, it follows that $(\tilde{X},\tilde{\Lambda})$
solves, $%
\tilde{P}$-almost surely, the Skorohod problem for $(D,\Gamma,\tilde{Z})$.
In particular, $(\tilde{X},\tilde{\Lambda})\in\mathcal{D} ( [ 0,T%
] ,%
\mathbb{R}
^{d} ) \times\mathcal{BV} ( [ 0,T ] ,%
\mathbb{R}
^{d} ) $ and
%
\begin{eqnarray}\label{RSDEth1a}
\tilde{X}_{t} &=&\hat{z}+\tilde{Z}_{t}+\tilde{\Lambda}_{t},
\\\label{RSDEth1a+}
\tilde{\Lambda}_{t} &=&\int_{0}^{t}\gamma_{s}\,d|\tilde{\Lambda}|_{s},\qquad
\gamma_{s}\in\Gamma_{s}(\tilde{X}_{s})\cap S_{1}(0),\  d|\tilde
{\Lambda}|%
\mbox{-a.e.}, \\ \label{RSDEth2a}
\tilde{X}_{t} &\in&
\overline{D_{t}},\qquad d|\tilde{\Lambda}|(\{t\in [
0,T ] \dvtx \tilde{X}_{t}\in D_{t}\})=0
\end{eqnarray}
holds $\tilde{P}$-almost surely whenever $t\in [ 0,T ] $. We next
want to verify that
%
\begin{equation}
\tilde{Z}_{t}=\int_{0}^{t}b(s,\tilde{X}_{s})\,ds+\int_{0}^{t}\sigma
(s,%
\tilde{X}_{s})\,d\tilde{W}_{s} \label{hahaga}
\end{equation}
holds $\tilde{P}$-almost surely whenever $t\in\lbrack0,T]$. Indeed,
using %
\eqref{RSDE2} and \eqref{RSDE1}, we can, following \cite{Costantini1992},
simply quote Theorem 2.2 in \cite{KurtzProtter1991}, which in our case
states that\vspace*{1.5pt} since $\{(\sigma(s,\tilde{X}_{s}^{n}),\tilde
{W}_{s}^{n})\}$
converges to $(\sigma(s,\tilde{X}_{s}),\tilde{W}_{s})$ $\tilde{P}$-almost
surely whenever $s\in\lbrack0,T]$, %
%
\begin{eqnarray} \label{hahaga+++}
\int_{0}^{t}\sigma(s,\tilde{X}_{s}^{n})\,d\tilde{W}_{s}^{n}\quad\mbox{converges
to}\quad \int_{0}^{t}\sigma(s,\tilde{X}_{s})\,d\tilde{W}_{s},\nonumber
\\[-8pt]
\\[-8pt]
 \eqntext{\tilde{P}%
\mbox{-almost surely,}}
\end{eqnarray}
whenever $t\in [ 0,T ] $, as $n\rightarrow\infty$. This proves %
\eqref{hahaga}. Now let $\mathcal{\tilde{F}}_{t}$ and $\mathcal
{\tilde{F}}%
_{t}^{n}$ be the $\sigma$-algebras generated by $\{\tilde
{W}_{s}\dvtx s\leq t\}$
and $\{\tilde{W}_{s}^{n}\dvtx s\leq t\}$, respectively. We next prove that $
\tilde{W}$ is a $m$-dimensional Wiener process on the filtered probability
space $(\tilde{\Omega},\mathcal{\tilde{F},}\{\mathcal{\tilde
{F}}_{t}\},%
\tilde{P})$. To obtain this, we first note that the $\sigma$-algebra
generated by $\{\tilde{W}_{s}^{n}-\tilde{W}_{t}^{n}\dvtx s\geq t\}$, for
$t\in
\lbrack0,T]$, is independent of $\mathcal{\tilde{F}}_{t}^{n}$. Furthermore,
since $\tilde{W}^{n}\rightarrow\tilde{W}$ $\tilde{P}$-almost
surely, it
follows that the $\sigma$-algebra generated by $\{\tilde
{W}_{s}-\tilde{W}%
_{t}\dvtx s\geq t\}$, for $t\in\lbrack0,T]$, is independent of $\mathcal
{\tilde{%
F}}_{t}$. In particular, $\{\tilde{W}_{t}\dvtx t\in\lbrack0,T]\}$ is a
martingale with respect to $\{\mathcal{\tilde{F}}_{t}\dvtx t\in\lbrack
0,T]\}$
and $\tilde{P}$. Now let $\tilde{W}_{t}^{n}=(\tilde
{W}_{t}^{n,1},\ldots,\tilde{W%
}_{t}^{n,m})$ and $\tilde{W}_{t}=(\tilde{W}_{t}^{1},\ldots,\tilde
{W}_{t}^{m})$.
Then, using essentially the same argument as in \eqref{hahaga+++}, we also
see that $\tilde{W}^{n,i}\tilde{W}^{n,j}\rightarrow\tilde
{W}^{i}\tilde{W}%
^{j}$, $\tilde{P}$-almost surely, for all $i,j\in\{1,\ldots,m\}$,
and, hence,
as above, we can conclude that $\{\tilde{W}_{t}^{i}\tilde
{W}_{t}^{j}-\delta
_{ij}t\dvtx t\in\lbrack0,T]\}$, with $\delta_{ij}$ being the Kronecker delta,
is a martingale with respect to $\{\mathcal{\tilde{F}}_{t}\dvtx t\in
\lbrack
0,T]\}$ and $\tilde{P}$. By the L\'{e}vy characterization of Wiener
processes (see, e.g., Theorem II.6.1 in \cite
{IkedaWatanabe1989}), we
can thus conclude that $\tilde{W}$ is a $m$-dimensional Wiener process
on $(%
\tilde{\Omega},\mathcal{\tilde{F},}\{\mathcal{\tilde{F}}_{t}\}
,\tilde{P})$.
To finally conclude that $(\tilde{X},\tilde{\Lambda})$ is a weak
solution in the sense of Definition \ref{weaksol}, and hence to
complete the
proof of Theorem \ref{Theorem 3}, it only remains to prove the
existence of
a version of $\tilde{X}$ on $(\tilde{\Omega},\mathcal{\tilde
{F}},\tilde{P})$%
, denoted $\hat{X}$, such that $\hat{X}\in\mathcal{C}( [ 0,T ] ,%
\mathbb{R}
^{d})$, $\tilde{P}$-almost surely. However, using standard arguments, we
first note that there exists a version of $\tilde{Z}$ on $(\tilde
{\Omega},%
\mathcal{\tilde{F}},\tilde{P})$, denoted $\hat{Z}$, such that $\hat
{Z}\in
\mathcal{C}( [ 0,T ] ,%
\mathbb{R}
^{d})$, $\tilde{P}$-almost surely, and such that $(\tilde{X},\tilde
{\Lambda}%
) $ solves, $\tilde{P}$-almost surely, the Skorohod problem for
$(D,\Gamma,%
\hat{Z})$. Furthermore, by \eqref{boliden}, it is clear that $\tilde
{X}\in
\mathcal{D}^{\rho_{0}} ( [ 0,T ] ,%
\mathbb{R}
^{d} ) $. Hence, by Theorem \ref{Theorem 2}, $\tilde{X}$ is
continuous $%
\tilde{P}$-almost surely. This completes the proof of Theorem \ref
{Theorem 3}.
\end{pf*}

\begin{appendix}\label{appx}

\section*{Appendix: Geometry of time-dependent domains}
\setcounter{equation}{0}

Concerning the $(\delta_{0},h_{0})$-property of good projections along
$%
\Gamma$, the following result follows immediately from Theorem 4.1 in
\cite{Costantini1992}.

\begin{lemma}
\label{costprop1} Let $T>0$ and let $D\subset
\mathbb{R}
^{d+1}$ be a time-dependent domain satisfying \eqref{timedep+} and a uniform
exterior sphere condition in time with radius $r_{0}$ in the sense of %
\eqref{extsphere-}. Let $\Gamma=\Gamma_{t} ( z ) $ be a closed
convex cone of vectors in $%
\mathbb{R}
^{d}$ for every $z\in\partial D_{t}$, $t\in\lbrack0,T]$ and assume
that $%
\Gamma$ satisfies \eqref{cone1} and \eqref{cone2}. Assume that there exists
a continuous map $Q\dvtx G^{N}\rightarrow
\mathbb{R}
^{d}$ such that
%
\begin{eqnarray}\label{Qdef}
Q ( t,z,N_{t} ( z ) ) &=&\Gamma_{t} ( z ) \qquad
\mbox{for all }z\in\partial D_{t},t\in [ 0,T ] ,\nonumber
\\
Q ( t,z,\lambda v ) &=&\lambda Q ( t,z,v )\\
    \eqntext{\mbox{for all
}\lambda\geq0, v\in N_{t} ( z ) , z\in\partial
D_{t}, t\in [ 0,T ] .}
\end{eqnarray}
Moreover, assume that
%
\begin{eqnarray}\label{Qdef1}
\sup_{t\in\lbrack0,T], z\in\partial D_{t}}\max_{v\in N_{t}^{1} (
z ) } \vert Q ( t,z,v ) \vert&:=& \Vert
Q \Vert<\infty,\nonumber
\\[-8pt]
\\[-8pt]
\inf_{t\in\lbrack0,T], z\in\partial D_{t}}\min_{v\in N_{t}^{1} (
z ) }v\cdot Q ( t,z,v ) &:=&q>0.\nonumber
\end{eqnarray}
Let
%
\begin{eqnarray}\label{nolldef}
\delta_{0}&:=&r_{0} \bigl( 1-\sqrt{1- ( q/ \Vert Q \Vert )
^{2}} \bigr) ,\nonumber
\\[-8pt]
\\[-8pt] h_{0}&:=&\frac{q/ \Vert Q \Vert}{1-\sqrt{%
1- ( q/ \Vert Q \Vert ) ^{2}}}.\nonumber
\end{eqnarray}
Then $([0,T]\times
\mathbb{R}
^{d})\setminus\overline{D}$ has the $(\delta_{0},h_{0})$-property of good
projections along $\Gamma$.
\end{lemma}

Note that in Lemma \ref{costprop1} we have that $q/ \Vert Q \Vert
<1 $, $0<\delta_{0}<r_{0}$ and $h_{0}>1$ by construction. To continue,
given $T>0$ and $D$ as above, we say that $D$ is a $\mathcal
{H}_{1+\alpha}$%
-domain if we can find a $\rho>0$ such that, for all $z_{0}\in
\partial
D_{t_{0}}$, $t_{0}\in\lbrack0,T]$, there exists a function $\psi (
t,z ) $, $\psi\in\mathcal{H}_{1+\alpha} ( C_{\rho} (
t_{0},z_{0} ) ) $, with the properties%
%
\begin{eqnarray}\label{H1alpha}
D\cap C_{\rho} ( t_{0},z_{0} ) &=& \{ \psi ( t,z )
>0 \} \cap C_{\rho} ( t_{0},z_{0} ) , \nonumber\\
\partial D\cap C_{\rho} ( t_{0},z_{0} ) &=& \{ \psi (
t,z ) =0 \} \cap C_{\rho} ( t_{0},z_{0} ) , \\
\inf_{ ( t,z ) \in\partial D\cap C_{\rho} (
t_{0},z_{0} ) } \vert\nabla_{z}\psi ( t,z ) \vert
&>&0\nonumber
\end{eqnarray}
for all $ ( t,z ) \in ( 0,T ) \times
\mathbb{R}
^{d}$.

\begin{lemma}
\label{compactprop} Let $T>0$ and let $D\subset
\mathbb{R}
^{d+1}$ be a time-dependent domain satisfying \eqref{timedep+} and a uniform
exterior sphere condition in time with radius $r_{0}$ in the sense of %
\eqref{extsphere-}. Assume, in addition, that $D$ is a $\mathcal
{H}_{1+\alpha
} $-domain for some $\alpha\in ( 0,1 ] $. Let $\Gamma=\Gamma
_{t}(z)$ be, for every $z\in\partial D_{t}$, $t\in\lbrack0,T]$, a closed
convex cone of vectors in $%
\mathbb{R}
^{d}$ with the specific form $\{\lambda\gamma_{t} ( z ) \dvtx
\lambda>0\}$, for some $S_{1} ( 0 ) $-valued function $\gamma
_{t} ( z ) $ which is uniformly continuous, in both space and time,
and satisfies
%
\begin{equation}
\beta=\inf_{t\in [ 0,T ] }\inf_{z\in\partial D_{t}} \langle
\gamma_{t} ( z ) ,n_{t} ( z ) \rangle>0.
\label{condbeta}
\end{equation}
Then $D$ satisfies \eqref{crita} and \eqref{crite} and
%
\begin{equation}
l ( r ) \leq Lr^{\tilde{\alpha}}  \qquad\mbox{whenever }r\in [ 0,T%
]  \label{Holder}
\end{equation}
for some $0<L<\infty$ and with $\tilde{\alpha}=(1+\alpha)/2\in (
0,1%
] $.
\end{lemma}

\begin{pf} By the uniform continuity of $\gamma_{t} (
z ) $ in space and time, it is clear that the variation of $\gamma
_{t} ( z ) $ can be made arbitrarily small on temporal
neighborhoods. Following Proposition 2.5 in \cite{Costantini1992}, it
therefore immediately follows that criteria~(\ref{crita}) and (\ref{crite})
are satisfied for some $0<\rho_{0}<r_{0}$ and $\eta_{0}>0$. Hence, it
remains to prove (\ref{Holder}) and we note that it suffices to prove
(\ref%
{Holder}) for small values of $r$. Let $z\in
\mathbb{R}
^{d}$ be arbitrary and let $z_{s}\in\partial D_{s}$ be such that $%
\vert z-z_{s} \vert=d ( z,\partial D_{s} ) $. We now
claim that
%
\begin{equation}
z-z_{s}\parallel n_{s} ( z_{s} ) \label{2}
\end{equation}
or, in other words, that $z-z_{s}$ and $n_{s} ( z_{s} ) $ are
parallel. To prove this claim, we can assume, without loss of generality,
that $z_{s}=0$. As $D$ is a time-dependent domain of class $\mathcal{H}
_{1+\alpha}$, we can assume the existence of a function $\psi$, with
property (\ref{H1alpha}), such that $\{\psi(s,y)=0\dvtx  y\in\overline
{B_{\rho
}(0)}\cap\partial D_{s}\}$, where $B_{\rho} ( 0 ) $ is a
(spatial) neighborhood of the origin with the radius $\rho$ as given in
the definition of $\mathcal{H}_{1+\alpha}$-domains. We consider the
minimization problem
%
\begin{equation}
\min_{y\in\overline{B}_{\rho}(0)\cap\partial D_{s}}|z-y|^{2}.
\label{minprob}
\end{equation}
Then, as the minimum in (\ref{minprob}) is realized at the origin, we see
that $z=\lambda\nabla\psi(s,0)$ for some Lagrange multiplier
$\lambda$.
Obviously, this proves (\ref{2}). Next, by Lemma~\ref
{llhatequivalence}, we
see that, for $r<\rho\wedge r_{0}$,
%
\begin{equation}
l ( r ) =\mathop{\mathop{\sup}_{ s,t\in[ 0,T]}}_{ \vert
s-t \vert\leq r}\sup_{z\in\overline{D_{s}}}d ( z,D_{t} )
=\mathop{\mathop{\sup}_{ 0\leq s\leq t\leq T}}_{ \vert s-t \vert
\leq r}%
\vert z_{t}-z_{s} \vert \label{lrep}
\end{equation}
for some $z_{s}\in\partial D_{s}$, $z_{t}\in\partial D_{t}$ such that
%
\begin{equation}
z_{t}-z_{s}\parallel n_{s} ( z_{s} )   \label{para}
\end{equation}
and $ \vert z_{t}-z_{s} \vert<\rho\wedge r_{0}$. Furthermore,
employing once more the fact that $D$ is a time-dependent domain of
class $%
\mathcal{H}_{1+\alpha}$, we conclude the existence of a function
$\psi$,
with the property $\psi ( s,z_{s} ) =\psi ( t,z_{t} ) =0$%
, which is continuously differentiable in space. Taylor expanding $\psi
$ up
to the first order in spatial coordinates, we obtain, by the $\mathcal
{H}%
_{1+\alpha}$-regularity of $\psi$,%
%
\begin{eqnarray}\label{Taylor}
\psi ( t,z_{t} ) -\psi ( s,z_{s} ) &=&\psi (
t,z_{t} ) -\psi ( s,z_{t} ) +\psi ( s,z_{t} ) -\psi
( s,z_{s} ) \nonumber\\
&=&\psi ( t,z_{t} ) -\psi ( s,z_{t} ) + \langle
z_{t}-z_{s},\nabla_{z}\psi ( s,z_{s} ) \rangle \\
&&{}+\mathcal{O} ( \vert z_{t}-z_{s} \vert^{1+\alpha} ) .\nonumber
\end{eqnarray}
As $n_{s} ( z_{s} ) =\frac{\nabla_{z}\psi (
s,z_{s} ) }{ \vert\nabla_{z}\psi ( s,z_{s} )
\vert}$, it follows from (\ref{lrep})--(\ref{Taylor}) that
%
\begin{eqnarray}
l ( r ) &=&\mathop{\mathop{\sup}_{ 0\leq s\leq t\leq T }}_{ \vert
s-t \vert\leq r} \vert z_{t}-z_{s} \vert=\mathop{\mathop{\sup}_{ %
0\leq s\leq t\leq T}}_{ \vert s-t \vert\leq r}\frac{%
| \langle z_{t}-z_{s},\nabla_{z}\psi ( s,z_{s} )
\rangle|}{ \vert\nabla_{z}\psi ( s,z_{s} )
\vert}\nonumber
\\[-8pt]
\\[-8pt]
&\leq&\mathop{\mathop{\sup}_{ 0\leq s\leq t\leq T}}_{ \vert s-t
\vert
\leq r}\frac{ \vert\psi ( t,z_{t} ) -\psi (
s,z_{t} ) \vert}{ \vert\nabla_{z}\psi (
s,z_{s} ) \vert}+\mathop{\mathop{\sup}_{ 0\leq s\leq t\leq T}}_{ %
\vert s-t \vert\leq r}\frac{\mathcal{O} ( \vert
z_{t}-z_{s} \vert^{1+\alpha} ) }{ \vert\nabla_{z}\psi
( s,z_{s} ) \vert}.\nonumber
\end{eqnarray}
Furthermore, for small $r$, we can assume, without loss of generality, that
for some $\delta>0$ independent of $ ( s,z_{s} ) $ we have $%
\vert\nabla_{z}\psi ( s,z_{s} ) \vert\geq\delta$.
Hence, by the $\mathcal{H}_{1+\alpha}$-regularity of $\psi$ and the
definition of $l$, we see that
%
\begin{equation}
l ( r ) \leq c\delta^{-1}\bigl(r^{ ( \alpha+1 )
/2}+l(r)^{1+\alpha}\bigr),
\end{equation}
provided $r<\rho\wedge r_{0}$. Finally, since $l(r)\rightarrow0$ as $%
r\rightarrow0$, there exists $\epsilon>0$ such that if $r\leq
\epsilon$,
then $c\delta^{-1}l(r)^{\alpha}\leq1/2$. Combining these facts, we
conclude that
%
\begin{equation}
l ( r ) \leq Lr^{ ( \alpha+1 ) /2}
\end{equation}
for some constant $L$. Hence, (\ref{Holder}) holds with $\tilde
{\alpha}%
= ( \alpha+1 ) /2$.
\end{pf}

\begin{remark}
In the following let $\mathcal{C}_{b}^{1}$ and $\mathcal{C}_{b}^{2}$ be
spaces containing all functions with bounded derivatives up to orders one
and two, respectively. Consider\vspace*{1pt} a bounded spatial domain $\Omega
\subset
\mathbb{R}
^{d}$ which is $\mathcal{C}_{b}^{1}$-smooth and satisfies a uniform exterior
sphere condition. Moreover, assume that the cone of directions of reflection
$\Gamma ( z ) $ has the specific form $\{\lambda\gamma (
z ) \dvtx  \lambda>0\}$, for some $S_{1} ( 0 ) $-valued function $%
\gamma ( z ) $ which is continuous and satisfies%
%
\begin{equation}
\beta:=\inf_{z\in\partial\Omega} \langle\gamma ( z )
,n ( z ) \rangle>0.
\end{equation}
Then Proposition 2.5 in \cite{Costantini1992} states that the
time-independent counterparts of criteria \eqref{crita} and \eqref{crite}
are satisfied. Furthermore, Theorem 4.5 in \cite{Costantini1992}
states that
the time-independent counterparts of \eqref{crita} and \eqref{crite} are
also satisfied for piecewise $\mathcal{C}_{b}^{1}$-smooth domains
$\Omega$
if the function $\gamma ( z ) $ is uniformly continuous on each
face of $\partial\Omega$ and satisfies some nondegeneracy and consistency
criteria. Finally, we also mention Theorem 4.6 in \cite{Costantini1992}
which states that unique projections may be found if~$\Omega$ is a
piecewise $\mathcal{C}_{b}^{2}$-smooth domain and if $\gamma ( z )
$ is Lipschitz continuous on each face of $\partial\Omega$ and satisfies
some other nondegeneracy and consistency criteria.
\end{remark}
\end{appendix}

\section*{Acknowledgments}

The authors would like to express their sincere gratitude to several
anonymous referees for their valuable comments and suggestions, which have
resulted in considerable improvement of the results and presentation of
this article.

%

\printaddresses


\begin{thebibliography}{00}

\bibitem{AtarBudhirajaRamanan2008}
%
\begin{barticle}[mr]
\bauthor{\bsnm{Atar},~\bfnm{Rami}\binits{R.}},
\bauthor{\bsnm{Budhiraja},~\bfnm{Amarjit}\binits{A.}} \AND
\bauthor{\bsnm{Ramanan},~\bfnm{Kavita}\binits{K.}}
(\byear{2008}).
\btitle{Deterministic and stochastic differential inclusions with multiple
surfaces of discontinuity}.
\bjournal{Probab. Theory Related Fields}
\bvolume{142}
\bpages{249--283}.
\bid{doi={10.1007/s00440-007-0104-z}, mr={2413272}}
\end{barticle}
%
\endbibitem

\bibitem{AtarDupuis1999}
%
\begin{barticle}[mr]
\bauthor{\bsnm{Atar},~\bfnm{Rami}\binits{R.}} \AND
\bauthor{\bsnm{Dupuis},~\bfnm{Paul}\binits{P.}}
(\byear{1999}).
\btitle{Large deviations and queueing networks: Methods for rate function
identification}.
\bjournal{Stochastic Process. Appl.}
\bvolume{84}
\bpages{255--296}.
\bid{doi={10.1016/S0304-4149(99)00051-4}, mr={1719274}}
\end{barticle}
%
\endbibitem

\bibitem{AtarDupuis2002}
%
\begin{barticle}[mr]
\bauthor{\bsnm{Atar},~\bfnm{Rami}\binits{R.}} \AND
\bauthor{\bsnm{Dupuis},~\bfnm{Paul}\binits{P.}}
(\byear{2002}).
\btitle{A differential game with constrained dynamics and viscosity solutions
of a related {HJB} equation}.
\bjournal{Nonlinear Anal.}
\bvolume{51}
\bpages{1105--1130}.
\bid{doi={10.1016/S0362-546X(01)00134-1}, mr={1926618}}
\end{barticle}
%
\endbibitem

\bibitem{BassHsu2000}
%
\begin{barticle}[mr]
\bauthor{\bsnm{Bass},~\bfnm{Richard~F.}\binits{R.~F.}} \AND
\bauthor{\bsnm{Hsu},~\bfnm{Elton~P.}\binits{E.~P.}}
(\byear{2000}).
\btitle{Pathwise uniqueness for reflecting {B}rownian motion in {E}uclidean
domains}.
\bjournal{Probab. Theory Related Fields}
\bvolume{117}
\bpages{183--200}.
\bid{doi={10.1007/s004400050003}, mr={1771660}}
\end{barticle}
%
\endbibitem

\bibitem{BenjaminiChenRodhe2004}
%
\begin{barticle}[mr]
\bauthor{\bsnm{Benjamini},~\bfnm{Itai}\binits{I.}},
\bauthor{\bsnm{Chen},~\bfnm{Zhen-Qing}\binits{Z.-Q.}} \AND
\bauthor{\bsnm{Rohde},~\bfnm{Steffen}\binits{S.}}
(\byear{2004}).
\btitle{Boundary trace of reflecting {B}rownian motions}.
\bjournal{Probab. Theory Related Fields}
\bvolume{129}
\bpages{1--17}.
\bid{doi={10.1007/s00440-003-0318-7}, mr={2052860}}
\end{barticle}
%
\endbibitem

\bibitem{BernardElKharroubi1991}
%
\begin{barticle}[mr]
\bauthor{\bsnm{Bernard},~\bfnm{Alain}\binits{A.}} \AND
\bauthor{\bparticle{el~}\bsnm{Kharroubi},~\bfnm{Ahmed}\binits{A.}}
(\byear{1991}).
\btitle{R\'egulations d\'eterministes et stochastiques dans le premier
``orthant'' de {${\mathbb R}^n$}}.
\bjournal{Stochastics Stochastics Rep.}
\bvolume{34}
\bpages{149--167}.
\bid{mr={1124833}}
\end{barticle}
%
\endbibitem

\bibitem{Billingsley1999}
%
\begin{bbook}[mr]
\bauthor{\bsnm{Billingsley},~\bfnm{Patrick}\binits{P.}}
(\byear{1999}).
\btitle{Convergence of Probability Measures},
\bedition{2nd} ed.
\bpublisher{Wiley}, \baddress{New York}.
\bid{doi={10.1002/9780470316962}, mr={1700749}}
\end{bbook}
%
\endbibitem

\bibitem{Borkowski2007}
%
\begin{bincollection}[mr]
\bauthor{\bsnm{Borkowski},~\bfnm{Dariusz}\binits{D.}}
(\byear{2007}).
\btitle{Chromaticity denoising using solution to the {S}korokhod problem}.
In \bbooktitle{Image Processing Based on Partial Differential
Equations}
(\beditor{X.-C. Tai et al.}, eds.)
\bpages{149--161}.
\bpublisher{Springer}, \baddress{Berlin}.
\bid{doi={10.1007/978-3-540-33267-1_9}, mr={2424226}}
\end{bincollection}
%
\endbibitem

\bibitem{BurdzyChenSylvester2003}
%
\begin{barticle}[mr]
\bauthor{\bsnm{Burdzy},~\bfnm{Krzysztof}\binits{K.}},
\bauthor{\bsnm{Chen},~\bfnm{Zhen-Qing}\binits{Z.-Q.}} \AND
\bauthor{\bsnm{Sylvester},~\bfnm{John}\binits{J.}}
(\byear{2003}).
\btitle{The heat equation and reflected {B}rownian motion in time-dependent
domains. {II}. {S}ingularities of solutions}.
\bjournal{J. Funct. Anal.}
\bvolume{204}
\bpages{1--34}.
\bid{doi={10.1016/S0022-1236(03)00128-9}, mr={2004743}}
\end{barticle}
%
\endbibitem

\bibitem{BurdzyChenSylvester2004AP}
%
\begin{barticle}[mr]
\bauthor{\bsnm{Burdzy},~\bfnm{Krzysztof}\binits{K.}},
\bauthor{\bsnm{Chen},~\bfnm{Zhen-Qing}\binits{Z.-Q.}} \AND
\bauthor{\bsnm{Sylvester},~\bfnm{John}\binits{J.}}
(\byear{2004}).
\btitle{The heat equation and reflected {B}rownian motion in time-dependent
domains}.
\bjournal{Ann. Probab.}
\bvolume{32}
\bpages{775--804}.
\bid{doi={10.1214/aop/1079021464}, mr={2039943}}%
\end{barticle}
%
\endbibitem

\bibitem{BurdzyChenSylvester2004JMAA}
%
\begin{barticle}[mr]
\bauthor{\bsnm{Burdzy},~\bfnm{Chris}\binits{C.}},
\bauthor{\bsnm{Chen},~\bfnm{Zhen-Qing}\binits{Z.-Q.}} \AND
\bauthor{\bsnm{Sylvester},~\bfnm{John}\binits{J.}}
(\byear{2004}).
\btitle{The heat equation in time dependent domains with insulated boundaries}.
\bjournal{J. Math. Anal. Appl.}
\bvolume{294}
\bpages{581--595}.
\bid{doi={10.1016/j.jmaa.2004.02.032}, mr={2061344}}
\end{barticle}
%
\endbibitem

\bibitem{BurdzyKangRamanan2009}
%
\begin{barticle}[mr]
\bauthor{\bsnm{Burdzy},~\bfnm{Krzysztof}\binits{K.}},
\bauthor{\bsnm{Kang},~\bfnm{Weining}\binits{W.}} \AND
\bauthor{\bsnm{Ramanan},~\bfnm{Kavita}\binits{K.}}
(\byear{2009}).
\btitle{The {S}korokhod problem in a time-dependent interval}.
\bjournal{Stochastic Process. Appl.}
\bvolume{119}
\bpages{428--452}.
\bid{doi={10.1016/j.spa.2008.03.001}, mr={2493998}}
\end{barticle}
%
\endbibitem

\bibitem{BurdzyNualart2002}
%
\begin{barticle}[mr]
\bauthor{\bsnm{Burdzy},~\bfnm{Krzysztof}\binits{K.}} \AND
\bauthor{\bsnm{Nualart},~\bfnm{David}\binits{D.}}
(\byear{2002}).
\btitle{Brownian motion reflected on {B}rownian motion}.
\bjournal{Probab. Theory Related Fields}
\bvolume{122}
\bpages{471--493}.
\bid{doi={10.1007/s004400100165}, mr={1902187}}
\end{barticle}
%
\endbibitem

\bibitem{BurdzyToby1995}
%
\begin{barticle}[mr]
\bauthor{\bsnm{Burdzy},~\bfnm{Krzysztof}\binits{K.}} \AND
\bauthor{\bsnm{Toby},~\bfnm{Ellen}\binits{E.}}
(\byear{1995}).
\btitle{A {S}korohod-type lemma and a decomposition of reflected {B}rownian
motion}.
\bjournal{Ann. Probab.}
\bvolume{23}
\bpages{586--604}.
\bid{mr={1334162}}
\end{barticle}
%
\endbibitem

\bibitem{Costantini1992}
%
\begin{barticle}[mr]
\bauthor{\bsnm{Costantini},~\bfnm{C.}\binits{C.}}
(\byear{1992}).
\btitle{The {S}korohod oblique reflection problem in domains with
corners and
application to stochastic differential equations}.
\bjournal{Probab. Theory Related Fields}
\bvolume{91}
\bpages{43--70}.
\bid{doi={10.1007/BF01194489}, mr={1142761}}
\end{barticle}
%
\endbibitem

\bibitem{CostantiniGobetKaroui2006}
%
\begin{barticle}[mr]
\bauthor{\bsnm{Costantini},~\bfnm{Cristina}\binits{C.}},
\bauthor{\bsnm{Gobet},~\bfnm{Emmanuel}\binits{E.}} \AND
\bauthor{\bsnm{El~Karoui},~\bfnm{Nicole}\binits{N.}}
(\byear{2006}).
\btitle{Boundary sensitivities for diffusion processes in time dependent
domains}.
\bjournal{Appl. Math. Optim.}
\bvolume{54}
\bpages{159--187}.
\bid{doi={10.1007/s00245-006-0863-4}, mr={2239532}}
\end{barticle}
%
\endbibitem

\bibitem{DupuisIshii1991}
%
\begin{barticle}[mr]
\bauthor{\bsnm{Dupuis},~\bfnm{Paul}\binits{P.}} \AND
\bauthor{\bsnm{Ishii},~\bfnm{Hitoshi}\binits{H.}}
(\byear{1991}).
\btitle{On {L}ipschitz continuity of the solution mapping to the {S}korokhod
problem, with applications}.
\bjournal{Stochastics Stochastics Rep.}
\bvolume{35}
\bpages{31--62}.
\bid{mr={1110990}}%
\end{barticle}
%
\endbibitem

\bibitem{DupuisIshii1993}
%
\begin{barticle}[mr]
\bauthor{\bsnm{Dupuis},~\bfnm{Paul}\binits{P.}} \AND
\bauthor{\bsnm{Ishii},~\bfnm{Hitoshi}\binits{H.}}
(\byear{1993}).
\btitle{S{DE}s with oblique reflection on nonsmooth domains}.
\bjournal{Ann. Probab.}
\bvolume{21}
\bpages{554--580}.
\bid{mr={1207237}}
\end{barticle}
%
\endbibitem

\bibitem{DupuisIshii2008}
%
\begin{barticle}[mr]
\bauthor{\bsnm{Dupuis},~\bfnm{Paul}\binits{P.}} \AND
\bauthor{\bsnm{Ishii},~\bfnm{Hitoshi}\binits{H.}}
(\byear{2008}).
\btitle{Correction: ``{SDE}s with oblique reflection on nonsmooth
domains.''}
\bjournal{Ann. Probab.}
\bvolume{36}
\bpages{1992--1997}.
\bid{doi={10.1214/07-AOP374}, mr={2440929}}
\end{barticle}
%
\endbibitem

\bibitem{DupuisRamanan1998}
%
\begin{barticle}[mr]
\bauthor{\bsnm{Dupuis},~\bfnm{Paul}\binits{P.}} \AND
\bauthor{\bsnm{Ramanan},~\bfnm{Kavita}\binits{K.}}
(\byear{1998}).
\btitle{A {S}korokhod problem formulation and large deviation analysis
of a
processor sharing model}.
\bjournal{Queueing Systems Theory Appl.}
\bvolume{28}
\bpages{109--124}.
\bid{mr={1628485}}
\end{barticle}
%
\endbibitem

\bibitem{DupuisRamanan1999b}
%
\begin{barticle}[vtex]
\bauthor{\bsnm{Dupuis},~\bfnm{Paul}\binits{P.}} \AND
\bauthor{\bsnm{Ramanan},~\bfnm{Kavita}\binits{K.}}
(\byear{1999}).
\btitle{Convex duality and the {S}korokhod problem. {II}}.
\bjournal{Probab. Theory Related Fields}
\bvolume{115}
\bpages{197--236}.
\end{barticle}
%
\endbibitem

\bibitem{DupuisRamanan1999a}
%
\begin{barticle}[vtex]
\bauthor{\bsnm{Dupuis},~\bfnm{Paul}\binits{P.}} \AND
\bauthor{\bsnm{Ramanan},~\bfnm{Kavita}\binits{K.}}
(\byear{1999}).
\btitle{Convex duality and the {S}korokhod problem. {I}}.
\bjournal{Probab. Theory Related Fields}
\bvolume{115}
\bpages{153--195}.
\end{barticle}
%
\endbibitem

\bibitem{DupuisRamanan2000a}
%
\begin{barticle}[mr]
\bauthor{\bsnm{Dupuis},~\bfnm{Paul}\binits{P.}} \AND
\bauthor{\bsnm{Ramanan},~\bfnm{Kavita}\binits{K.}}
(\byear{2000}).
\btitle{An explicit formula for the solution of certain optimal control
problems on domains with corners}.
\bjournal{Teory Probab. Math. Statist.}
\bvolume{63}
\bpages{33--49}.
\bid{mr={1870773}}%
\end{barticle}
%
\endbibitem

\bibitem{DupuisRamanan2000b}
%
\begin{barticle}[mr]
\bauthor{\bsnm{Dupuis},~\bfnm{Paul}\binits{P.}} \AND
\bauthor{\bsnm{Ramanan},~\bfnm{Kavita}\binits{K.}}
(\byear{2000}).
\btitle{A multiclass feedback queueing network with a regular {S}korokhod
problem}.
\bjournal{Queueing Systems Theory Appl.}
\bvolume{36}
\bpages{327--349}.
\bid{doi={10.1023/A:1011037419624}, mr={1823974}}
\end{barticle}
%
\endbibitem

\bibitem{DupuisWilliams1994}
%
\begin{barticle}[mr]
\bauthor{\bsnm{Dupuis},~\bfnm{Paul}\binits{P.}} \AND
\bauthor{\bsnm{Williams},~\bfnm{Ruth~J.}\binits{R.~J.}}
(\byear{1994}).
\btitle{Lyapunov functions for semimartingale reflecting {B}rownian motions}.
\bjournal{Ann. Probab.}
\bvolume{22}
\bpages{680--702}.
\bid{mr={1288127}}
\end{barticle}
%
\endbibitem

\bibitem{ElKarouiKapoudjianPardouxPengQuenez1997}
%
\begin{barticle}[mr]
\bauthor{\bsnm{El~Karoui},~\bfnm{N.}\binits{N.}},
\bauthor{\bsnm{Kapoudjian},~\bfnm{C.}\binits{C.}},
\bauthor{\bsnm{Pardoux},~\bfnm{E.}\binits{E.}},
\bauthor{\bsnm{Peng},~\bfnm{S.}\binits{S.}} \AND
\bauthor{\bsnm{Quenez},~\bfnm{M.~C.}\binits{M.~C.}}
(\byear{1997}).
\btitle{Reflected solutions of backward {SDE}'s, and related obstacle problems
for {PDE}'s}.
\bjournal{Ann. Probab.}
\bvolume{25}
\bpages{702--737}.
\bid{doi={10.1214/aop/1024404416}, mr={1434123}}
\end{barticle}
%
\endbibitem

\bibitem{ElKarouiKaratzas1991b}
%
\begin{barticle}[mr]
\bauthor{\bsnm{El~Karoui},~\bfnm{N.}\binits{N.}} \AND
\bauthor{\bsnm{Karatzas},~\bfnm{I.}\binits{I.}}
(\byear{1991}).
\btitle{Correction: ``{A} new approach to the {S}korohod problem, and its
applications''}.
\bjournal{Stochastics Stochastics Rep.}
\bvolume{36}
\bpages{265}.
\bid{mr={1128498}}
\end{barticle}
%
\endbibitem

\bibitem{ElKarouiKaratzas1991a}
%
\begin{barticle}[mr]
\bauthor{\bsnm{El~Karoui},~\bfnm{Nicole}\binits{N.}} \AND
\bauthor{\bsnm{Karatzas},~\bfnm{Ioannis}\binits{I.}}
(\byear{1991}).
\btitle{A new approach to the {S}korohod problem, and its applications}.
\bjournal{Stochastics Stochastics Rep.}
\bvolume{34}
\bpages{57--82}.
\bid{mr={1104422}}
\end{barticle}
%
\endbibitem

\bibitem{ElKharroubiBenTaharYaacoubi2002}
%
\begin{barticle}[mr]
\bauthor{\bsnm{El~Kharroubi},~\bfnm{Ahmed}\binits{A.}},
\bauthor{\bsnm{Ben~Tahar},~\bfnm{Abdelghani}\binits{A.}} \AND
\bauthor{\bsnm{Yaacoubi},~\bfnm{Abdelhak}\binits{A.}}
(\byear{2002}).
\btitle{On the stability of the linear {S}korohod problem in an orthant}.
\bjournal{Math. Methods Oper. Res.}
\bvolume{56}
\bpages{243--258}.
\bid{doi={10.1007/s001860200210}, mr={1938213}}
\end{barticle}
%
\endbibitem

\bibitem{EthierKurtz1986}
%
\begin{bbook}[mr]
\bauthor{\bsnm{Ethier},~\bfnm{Stewart~N.}\binits{S.~N.}} \AND
\bauthor{\bsnm{Kurtz},~\bfnm{Thomas~G.}\binits{T.~G.}}
(\byear{1986}).
\btitle{Markov Processes: Characterization and Convergence}.
\bpublisher{Wiley}, \baddress{New York}.
\bid{doi={10.1002/9780470316658}, mr={838085}}
\end{bbook}
%
\endbibitem

\bibitem{Frankowska1985}
%
\begin{barticle}[mr]
\bauthor{\bsnm{Frankowska},~\bfnm{Halina}\binits{H.}}
(\byear{1985}).
\btitle{A viability approach to the {S}korohod problem}.
\bjournal{Stochastics Stochastics Rep.}
\bvolume{14}
\bpages{227--244}.
\bid{mr={800245}}
\end{barticle}
%
\endbibitem

\bibitem{Friedman1982}
%
\begin{bbook}[mr]
\bauthor{\bsnm{Friedman},~\bfnm{Avner}\binits{A.}}
(\byear{1982}).
\btitle{Foundations of Modern Analysis}.
\bpublisher{Dover}, \baddress{New York}.
\bid{mr={663003}}
\end{bbook}
%
\endbibitem

\bibitem{HarrisonReiman1981}
%
\begin{barticle}[mr]
\bauthor{\bsnm{Harrison},~\bfnm{J.~Michael}\binits{J.~M.}} \AND
\bauthor{\bsnm{Reiman},~\bfnm{Martin~I.}\binits{M.~I.}}
(\byear{1981}).
\btitle{Reflected {B}rownian motion on an orthant}.
\bjournal{Ann. Probab.}
\bvolume{9}
\bpages{302--308}.
\bid{mr={606992}}
\end{barticle}
%
\endbibitem

\bibitem{IkedaWatanabe1989}
%
\begin{bbook}[mr]
\bauthor{\bsnm{Ikeda},~\bfnm{Nobuyuki}\binits{N.}} \AND
\bauthor{\bsnm{Watanabe},~\bfnm{Shinzo}\binits{S.}}
(\byear{1989}).
\btitle{Stochastic Differential Equations and Diffusion Processes},
\bedition{2nd} ed.
\bseries{North-Holland Math. Library}
\bvolume{24}.
\bpublisher{North-Holland}, \baddress{Amsterdam}.
\bid{mr={1011252}}
\end{bbook}
%
\endbibitem

\bibitem{KangWilliams2007}
%
\begin{barticle}[mr]
\bauthor{\bsnm{Kang},~\bfnm{W.}\binits{W.}} \AND
\bauthor{\bsnm{Williams},~\bfnm{R.~J.}\binits{R.~J.}}
(\byear{2007}).
\btitle{An invariance principle for semimartingale reflecting {B}rownian
motions in domains with piecewise smooth boundaries}.
\bjournal{Ann. Appl. Probab.}
\bvolume{17}
\bpages{741--779}.
\bid{doi={10.1214/105051606000000899}, mr={2308342}}
\end{barticle}
%
\endbibitem

\bibitem{KonstantopoulosAnantharam1995}
%
\begin{bmisc}[author]
\bauthor{\bsnm{Konstantopoulos},~\bfnm{Takis}\binits{T.}} \AND
\bauthor{\bsnm{Anantharam},~\bfnm{Venkat}\binits{V.}}
(\byear{1995}).
\bhowpublished{An
optimal flow control scheme that regulates the burstiness of traffic subject
to delay constraints. \textit{IEEE/ACM Transactions on Networking}}
\bvolume{3}
\bpages{423--432}.
\end{bmisc}
%
\endbibitem

\bibitem{Kruk2000}
%
\begin{barticle}[mr]
\bauthor{\bsnm{Kruk},~\bfnm{Lukasz}\binits{L.}}
(\byear{2000}).
\btitle{Optimal policies for {$n$}-dimensional singular stochastic control
problems. {I}. {T}he {S}korokhod problem}.
\bjournal{SIAM J. Control Optim.}
\bvolume{38}
\bpages{1603--1622 (electronic)}.
\bid{doi={10.1137/S0363012998347535}, mr={1766432}}
\end{barticle}
%
\endbibitem

\bibitem{KrukLehoczkyRamananShreve2007}
%
\begin{barticle}[mr]
\bauthor{\bsnm{Kruk},~\bfnm{Lukasz}\binits{L.}},
\bauthor{\bsnm{Lehoczky},~\bfnm{John}\binits{J.}},
\bauthor{\bsnm{Ramanan},~\bfnm{Kavita}\binits{K.}} \AND
\bauthor{\bsnm{Shreve},~\bfnm{Steven}\binits{S.}}
(\byear{2007}).
\btitle{An explicit formula for the {S}korokhod map on {$[0,a]$}}.
\bjournal{Ann. Probab.}
\bvolume{35}
\bpages{1740--1768}.
\bid{doi={10.1214/009117906000000890}, mr={2349573}}
\end{barticle}
%
\endbibitem

\bibitem{KurtzProtter1991}
%
\begin{barticle}[mr]
\bauthor{\bsnm{Kurtz},~\bfnm{Thomas~G.}\binits{T.~G.}} \AND
\bauthor{\bsnm{Protter},~\bfnm{Philip}\binits{P.}}
(\byear{1991}).
\btitle{Weak limit theorems for stochastic integrals and stochastic
differential equations}.
\bjournal{Ann. Probab.}
\bvolume{19}
\bpages{1035--1070}.
\bid{mr={1112406}}
\end{barticle}
%
\endbibitem

\bibitem{Kushner2001}
%
\begin{bbook}[mr]
\bauthor{\bsnm{Kushner},~\bfnm{Harold~J.}\binits{H.~J.}}
(\byear{2001}).
\btitle{Heavy Traffic Analysis of Controlled Queueing and Communication
Networks}.
\bseries{Appl. Math.}
\bvolume{47}.
\bpublisher{Springer}, \baddress{New York}.
\bid{mr={1834938}}
\end{bbook}
%
\endbibitem

\bibitem{Lieberman1996}
%
\begin{bbook}[mr]
\bauthor{\bsnm{Lieberman},~\bfnm{Gary~M.}\binits{G.~M.}}
(\byear{1996}).
\btitle{Second Order Parabolic Differential Equations}.
\bpublisher{World Scientific}, \baddress{River Edge, NJ}.
\bid{mr={1465184}}
\end{bbook}
%
\endbibitem

\bibitem{LionsSznitman1984}
%
\begin{bmisc}[author]
\bauthor{\bsnm{Lions},~\bfnm{Pierre-Louis}\binits{P.-L.}}
\AND
\bauthor{\bsnm{Sznitman},~\bfnm{Alain-Sol}\binits{A.-S.}}
(\byear{1984}).
\bhowpublished{Stochastic
differential equations with reflecting boundary conditions. {\em
Comm. Pure Appl. Math.}}
\bvolume{37}
\bpages{511--537}.
\end{bmisc}
%
\endbibitem

\bibitem{MandelbaumMassey1995}
%
\begin{barticle}[mr]
\bauthor{\bsnm{Mandelbaum},~\bfnm{Avi}\binits{A.}} \AND
\bauthor{\bsnm{Massey},~\bfnm{William~A.}\binits{W.~A.}}
(\byear{1995}).
\btitle{Strong approximations for time-dependent queues}.
\bjournal{Math. Oper. Res.}
\bvolume{20}
\bpages{33--64}.
\bid{doi={10.1287/moor.20.1.33}, mr={1320446}}
\end{barticle}
%
\endbibitem

\bibitem{Marin-RubioReal2004}
%
\begin{barticle}[mr]
\bauthor{\bsnm{Mar{\'{\i}}n-Rubio},~\bfnm{P.}\binits{P.}} \AND
\bauthor{\bsnm{Real},~\bfnm{J.}\binits{J.}}
(\byear{2004}).
\btitle{Some results on stochastic differential equations with reflecting
boundary conditions}.
\bjournal{J. Theoret. Probab.}
\bvolume{17}
\bpages{705--716}.
\bid{doi={10.1023/B:JOTP.0000040295.09922.85}, mr={2091557}}
\end{barticle}
%
\endbibitem

\bibitem{NystromOnskog2009c}
%
\begin{bmisc}[author]
\bauthor{\bsnm{Nystr\"{o}m},~\bfnm{Kaj}\binits{K.}} \AND
\bauthor{\bsnm{\"{O}nskog},~\bfnm{Thomas}\binits{T.}}
(\byear{2010}).
\bhowpublished{Weak approximation of
obliquely reflected diffusions in time dependent domains.  \textit{J. Comput. Math.}}
\bvolume{28}
\bpages{579--605}.
\end{bmisc}
%
\endbibitem

\bibitem{Ramanan2006}
%
\begin{barticle}[mr]
\bauthor{\bsnm{Ramanan},~\bfnm{Kavita}\binits{K.}}
(\byear{2006}).
\btitle{Reflected diffusions defined via the extended {S}korokhod map}.
\bjournal{Electron. J. Probab.}
\bvolume{11}
\bpages{934--992}
(electronic).
\bid{mr={2261058}}
\end{barticle}
%
\endbibitem

\bibitem{RamananReiman2003}
%
\begin{barticle}[mr]
\bauthor{\bsnm{Ramanan},~\bfnm{Kavita}\binits{K.}} \AND
\bauthor{\bsnm{Reiman},~\bfnm{Martin~I.}\binits{M.~I.}}
(\byear{2003}).
\btitle{Fluid and heavy traffic diffusion limits for a generalized processor
sharing model}.
\bjournal{Ann. Appl. Probab.}
\bvolume{13}
\bpages{100--139}.
\bid{doi={10.1214/aoap/1042765664}, mr={1951995}}
\end{barticle}
%
\endbibitem

\bibitem{RamananReiman2008}
%
\begin{barticle}[mr]
\bauthor{\bsnm{Ramanan},~\bfnm{Kavita}\binits{K.}} \AND
\bauthor{\bsnm{Reiman},~\bfnm{Martin~I.}\binits{M.~I.}}
(\byear{2008}).
\btitle{The heavy traffic limit of an unbalanced generalized processor sharing
model}.
\bjournal{Ann. Appl. Probab.}
\bvolume{18}
\bpages{22--58}.
\bid{doi={10.1214/07-AAP438}, mr={2380890}}
\end{barticle}
%
\endbibitem

\bibitem{Ramasubramanian2000}
%
\begin{barticle}[mr]
\bauthor{\bsnm{Ramasubramanian},~\bfnm{S.}\binits{S.}}
(\byear{2000}).
\btitle{A subsidy-surplus model and the {S}korokhod problem in an orthant}.
\bjournal{Math. Oper. Res.}
\bvolume{25}
\bpages{509--538}.
\bid{doi={10.1287/moor.25.3.509.12215}, mr={1855180}}
\end{barticle}
%
\endbibitem

\bibitem{Ramasubramanian2006}
%
\begin{barticle}[mr]
\bauthor{\bsnm{Ramasubramanian},~\bfnm{S.}\binits{S.}}
(\byear{2006}).
\btitle{An insurance network: {N}ash equilibrium}.
\bjournal{Insurance Math. Econom.}
\bvolume{38}
\bpages{374--390}.
\bid{doi={10.1016/j.insmatheco.2005.10.005}, mr={2212535}}
\end{barticle}
%
\endbibitem

\bibitem{Reiman1984}
%
\begin{barticle}[mr]
\bauthor{\bsnm{Reiman},~\bfnm{Martin~I.}\binits{M.~I.}}
(\byear{1984}).
\btitle{Open queueing networks in heavy traffic}.
\bjournal{Math. Oper. Res.}
\bvolume{9}
\bpages{441--458}.
\bid{doi={10.1287/moor.9.3.441}, mr={757317}}
\end{barticle}
%
\endbibitem

\bibitem{Robert2003}
%
\begin{bbook}[mr]
\bauthor{\bsnm{Robert},~\bfnm{Philippe}\binits{P.}}
(\byear{2003}).
\btitle{Stochastic Networks and Queues},
\bedition{French} ed.
\bseries{Appl. Math.}
\bvolume{52}.
\bpublisher{Springer}, \baddress{Berlin}.
\bid{mr={1996883}}
\end{bbook}
%
\endbibitem

\bibitem{Saisho1987}
%
\begin{barticle}[mr]
\bauthor{\bsnm{Saisho},~\bfnm{Yasumasa}\binits{Y.}}
(\byear{1987}).
\btitle{Stochastic differential equations for multidimensional domain with
reflecting boundary}.
\bjournal{Probab. Theory Related Fields}
\bvolume{74}
\bpages{455--477}.
\bid{doi={10.1007/BF00699100}, mr={873889}}
\end{barticle}
%
\endbibitem

\bibitem{Saisho1988}
%
\begin{bincollection}[mr]
\bauthor{\bsnm{Saisho},~\bfnm{Yasumasa}\binits{Y.}}
(\byear{1988}).
\btitle{Mutually repelling particles of {$m$} types}.
In \bbooktitle{Probability Theory and Mathematical Statistics
({K}yoto, \textit{1986})}.
\bseries{Lecture Notes in Math.}
\bvolume{1299}
\bpages{444--453}.
\bpublisher{Springer}, \baddress{Berlin}.
\bid{doi={10.1007/BFb0078503}, mr={936019}}
\end{bincollection}
%
\endbibitem

\bibitem{Saisho1991}
%
\begin{barticle}[mr]
\bauthor{\bsnm{Saisho},~\bfnm{Yasumasa}\binits{Y.}}
(\byear{1991}).
\btitle{On the equation describing the random motion of mutually reflecting
molecules}.
\bjournal{Proc. Japan Acad. Ser. A Math. Sci.}
\bvolume{67}
\bpages{293--298}.
\bid{mr={1151341}}
\end{barticle}
%
\endbibitem

\bibitem{Saisho1994}
%
\begin{barticle}[mr]
\bauthor{\bsnm{Saisho},~\bfnm{Yasumasa}\binits{Y.}}
(\byear{1994}).
\btitle{A model of the random motion of mutually reflecting molecules
in {$\mathbb R^ d$}}.
\bjournal{Kumamoto J. Math.}
\bvolume{7}
\bpages{95--123}.
\bid{mr={1273971}}
\end{barticle}
%
\endbibitem

\bibitem{Skorohod1961a}
%
\begin{barticle}[mr]
\bauthor{\bsnm{Skorohod},~\bfnm{A.~V.}\binits{A.~V.}}
(\byear{1961}).
\btitle{Stochastic equations for diffusion processes with a boundary}.
\bjournal{Teor. Verojatnost. i Primenen.}
\bvolume{6}
\bpages{287--298}.
\bid{mr={0145598}}
\end{barticle}
%
\endbibitem

\bibitem{SonerShreve1989}
%
\begin{barticle}[mr]
\bauthor{\bsnm{Soner},~\bfnm{H.~Mete}\binits{H.~M.}} \AND
\bauthor{\bsnm{Shreve},~\bfnm{Steven~E.}\binits{S.~E.}}
(\byear{1989}).
\btitle{Regularity of the value function for a two-dimensional singular
stochastic control problem}.
\bjournal{SIAM J. Control Optim.}
\bvolume{27}
\bpages{876--907}.
\bid{doi={10.1137/0327047}, mr={1001925}}
\end{barticle}
%
\endbibitem

\bibitem{SoucaliucWerner2002}
%
\begin{barticle}[mr]
\bauthor{\bsnm{Soucaliuc},~\bfnm{Florin}\binits{F.}} \AND
\bauthor{\bsnm{Werner},~\bfnm{Wendelin}\binits{W.}}
(\byear{2002}).
\btitle{A note on reflecting {B}rownian motions}.
\bjournal{Electron. Comm. Probab.}
\bvolume{7}
\bpages{117--122 (electronic)}.
\bid{mr={1917545}}
\end{barticle}
%
\endbibitem

\bibitem{Taksar1992}
%
\begin{barticle}[mr]
\bauthor{\bsnm{Taksar},~\bfnm{M.~I.}\binits{M.~I.}}
(\byear{1992}).
\btitle{Skorohod problems with nonsmooth boundary conditions}.
\bjournal{J. Comput. Appl. Math.}
\bvolume{40}
\bpages{233--251}.
\bid{doi={10.1016/0377-0427(92)90108-A}, mr={1170903}}
\end{barticle}
%
\endbibitem

\bibitem{Tanaka1979}
%
\begin{barticle}[mr]
\bauthor{\bsnm{Tanaka},~\bfnm{Hiroshi}\binits{H.}}
(\byear{1979}).
\btitle{Stochastic differential equations with reflecting boundary
condition in
convex regions}.
\bjournal{Hiroshima Math. J.}
\bvolume{9}
\bpages{163--177}.
\bid{mr={529332}}
\end{barticle}
%
\endbibitem

\end{thebibliography}
\end{document}